\documentclass[12pt,sans]{article}

\usepackage[margin=1in]{geometry}
\usepackage{microtype}
\usepackage{url}
\usepackage[colorlinks = true,
            linkcolor = blue,
            urlcolor  = blue,
            citecolor = blue,
            anchorcolor = blue]{hyperref}
\usepackage[utf8]{inputenc} 
\usepackage[T1]{fontenc}    
\usepackage{lmodern}
\usepackage{mathrsfs}
\usepackage[sans]{dsfont}   
\usepackage{booktabs}       
\usepackage{nicefrac}       
\usepackage{amsmath,amsthm,amssymb,amsfonts}   
\usepackage{latexsym}
\usepackage{optidef}
\usepackage{graphicx}
\usepackage{caption}
\usepackage{subcaption}
\usepackage{epstopdf}
\usepackage{cleveref}
\usepackage{color}
\usepackage{algorithm,algpseudocode}
\usepackage{tcolorbox}
\usepackage{xcolor} 
\newtcolorbox{tbox}[3][]{%
colframe=#2,colback=#2!10,coltitle=#2!20!black,title={#3},#1}
\usepackage{accents}
\usepackage{enumitem}
\usepackage{tikz-cd} 
\tikzset{smalltext/.style={"\textup{\small #1}" description}}
\usetikzlibrary{fit,backgrounds,arrows.meta}
\usepackage{pgfplots}
\pgfplotsset{compat=1.18}
\usepackage{silence}
\usepackage{autonum} 

\makeatletter
\newcommand{\@tvsp}{\mathchoice{{}\mkern-4.5mu}{{}\mkern-4.5mu}{{}\mkern-2.5mu}{}}
\makeatother

\newcommand{\esp}{\mathbb{E}}

\newcommand{\prob}{\mathbb{P}}
\newcommand{\probn}{\mathbf{P}}

\newcommand{\unit}{\mathbf{1}}

\newcommand{\re}{\mathbb{R}}
\DeclareMathOperator{\tr}{tr}
\DeclareMathOperator{\erank}{r}

\DeclareMathOperator{\interior}{int}
\DeclareMathOperator{\ri}{ri}

\DeclareMathOperator{\co}{conv}
\DeclareMathOperator{\ext}{ext}
\DeclareMathOperator{\KL}{KL}
\DeclareMathOperator{\E}{E}
\DeclareMathOperator{\dom}{dom}

\DeclareMathOperator{\Count}{Count}
\DeclareMathOperator{\IQG}{IQG}

\newcommand{\bfSigma}{\boldsymbol{\Sigma}}

\newcommand{\bfA}{\mathbf{A}}

\newcommand{\bfC}{\mathbf{C}}

\newcommand{\bfI}{\mathbf{I}}
\newcommand{\bfM}{\mathbf{M}}

\newcommand{\bfW}{\mathbf{W}}

\newcommand{\bfZ}{\mathbf{Z}}

\newcommand{\be}{\boldsymbol e}

\newcommand{\bg}{\boldsymbol g}
\newcommand{\bh}{\boldsymbol h}
\newcommand{\bx}{\boldsymbol x}

\newcommand{\bu}{\boldsymbol u}
\newcommand{\bv}{\boldsymbol v}

\newcommand{\bz}{\boldsymbol z}
\newcommand{\bzero}{\mathbf{0}}

\newcommand{\bmmu}{\boldsymbol\mu}
\newcommand{\bnu}{\boldsymbol\nu}

\newcommand{\bxi}{\boldsymbol \xi}

\newcommand{\calB}{\mathcal B}
\newcommand{\calC}{\mathcal C}

\newcommand{\calE}{\mathcal E}
\newcommand{\calF}{\mathcal F}

\newcommand{\calM}{\mathcal M}
\newcommand{\calN}{\mathcal N}

\newcommand{\calR}{\mathcal R}
\newcommand{\calS}{\mathcal S}

\newcommand{\calZ}{\mathcal Z}

\newcommand{\sfT}{\mathsf{T}}

\newcommand{\mbS}{\mathbb{S}}
\newcommand{\mbN}{\mathbb{N}}

\newcommand{\mbX}{\mathbb{X}}

\newcommand{\mcE}{\mathscr{E}}

\DeclareMathOperator*{\argmin}{argmin}

\theoremstyle{plain}
\newtheorem{theorem}{Theorem}
\newtheorem{proposition}[theorem]{Proposition}
\newtheorem{lemma}[theorem]{Lemma}
\newtheorem{corollary}[theorem]{Corollary}
\newtheorem{assumption}[theorem]{Assumption}

\theoremstyle{definition}
\newtheorem{definition}{Definition}

\theoremstyle{remark}

\newtheorem{remark}{Remark}

\numberwithin{equation}{section}

\title{The optimization landscape of robust high-dimensional mean estimation: stability, moment adaptation, and mirror descent reweighting}
\date{}

\author{
Roberto I. Oliveira \and
Zoraida F. Rico \and
Philip Thompson
}

\newcommand{\Addresses}{{
  \bigskip
  \footnotesize

Roberto I.~Oliveira, \textsc{IMPA, Brazil}\par\nopagebreak
  \textit{E-mail address}: \texttt{rbimfo@impa.br} 

\and 

Zoraida F.~Rico, \textsc{Bocconi University, Italy}\par\nopagebreak
  \textit{E-mail address}: \texttt{fr.zoraida@gmail.com} 

\and

  P.~Thompson, \textsc{FGV EMAp, School of Applied Mathematics, Brazil}\par\nopagebreak
  \textit{E-mail address}: \texttt{philip.thompson@fgv.br}

  \medskip
}}

\begin{document}

\maketitle
\Addresses

\begin{abstract}
We study efficient robust estimation of a high-dimensional mean under strong contamination when the clean distribution has only $p\ge2$ moments, bounded by $\nu_p$. We give a unified proof of stability for every $p\ge2$ with the optimal dimension-free sub-Gaussian rate, linking the reweighting, trimming and truncation techniques of the outlier and heavy-tailed robust literatures through concentration of counting and truncation processes. Building on stability and spectral signatures, we quantify the optimization landscape over sample weights that remove at most a fraction $\varepsilon_{n,k}=k/n$ of the mass, where the trimming parameter $k$ slightly exceeds the number of outliers. An unobserved convex spectral loss has inexact quadratic growth, with curvature $\frac{1-\varepsilon_{n,k}}{\varepsilon_{n,k}}$, relative to a mean-induced Bregman divergence, and the observable spectral loss supplies a dimension-free, inexact, relatively weakly star-convex subgradient oracle for it. This yields a Lyapunov analysis of mirror descent reweighting for $\varepsilon_{n,k}<\frac13$ that covers Hedge, the projected subgradient method, a new mean-induced mirror descent, and a centered multiplicative weights update. Each method uses one top-eigenvector computation per step and attains the optimal sub-Gaussian rate without knowledge of $(p,\nu_p)$. Mean-induced mirror descent further achieves faster optimization rates and last-iterate guarantees, at the cost of a more expensive projection step in higher dimensions.
\end{abstract}






\section{Introduction}

Robust statistics, initiated by Huber~\cite{1964Huber} (see also~\cite{2011Huber:Ronchetti}), studies estimators whose guarantees degrade gracefully when a small fraction $\varepsilon$ of the sample consists of outliers, caused, e.g., by measurement errors or corrupted data. One of its most basic tasks, used as a subroutine in, e.g., classification, regression and stochastic optimization, is to estimate the mean of a $d$-dimensional distribution from such a contaminated sample. The recent literature \cite{2016Diakonikolas:Kamath:Kane:Li:Moitra:Stewart,2020Diakonikolas:Kane:Pensia,2023Diakonikolas:Kane} adopts the following contamination model.

\begin{assumption}[Strong contamination with bounded $p$-th moments]
\label{assump:finite:covariance:contamination}
Let $\bx_1,\ldots,\bx_n$ be i.i.d.\ copies of a random vector $\bx\in\re^d$ with
distribution $\probn$, mean $\bmmu:=\esp[\bx]$ and covariance matrix
$\bfSigma\neq\bzero$. Assume that, for some $p\ge2$,
\[
    \nu_p := \sup_{\|\bv\|_2=1}
    \Big(\esp\big[|\langle\bx-\bmmu,\bv\rangle|^p\big]\Big)^{1/p} < \infty .
\]
Let $\varepsilon\in[0,1/2)$. The observed sample $\tilde\bx_1,\ldots,\tilde\bx_n$
consists of arbitrary random vectors in $\re^d$, possibly depending on
$\bx_1,\ldots,\bx_n$, such that
$
    \#\{i\in[n] : \tilde\bx_i\neq \bx_i\} \le \varepsilon n.
$
\end{assumption}
Our main goal in this paper is to understand the optimization aspects of estimating $\bmmu$ in the above model. Here $\nu_p^p$ is the largest directional $p$-th moment of $\bx$ about $\bmmu$, and $\varepsilon$ is the \emph{contamination fraction}. In particular, $\nu_2^2=\|\bfSigma\|$, the operator norm of $\bfSigma$ Assumption~\ref{assump:finite:covariance:contamination} describes the \emph{strong contamination model}: an adversary may inspect the clean sample and then replace up to an $\varepsilon$-fraction of it
arbitrarily~\cite{2019Diakonikolas:Kamath:Kane:Li:Moitra:Stewart,2021Lugosi:Mendelson}. This model is more general than Huber's  model~\cite{1964Huber}, in which the outliers are drawn independently of the inliers. Modern robust estimation seeks estimators that are computationally efficient and achieve the optimal rate in the dimension, $n$ and~$\varepsilon$. Notably, the sharp dependence on the dimension is not through $d$ but through the \emph{effective rank} $\tr(\bfSigma)/\|\bfSigma\|$, called \emph{dimension-free} rates \cite{2021Lugosi:Mendelson,2017Koltchinskii:Lounici}.

For $d=1$, mean estimation under moment assumptions is well understood: the trimmed mean is efficient and attains the optimal rate~\cite{2026Oliveira:Orenstein:Rico}. In high dimensions, however, there is a trade-off: the multivariate trimmed mean of~\cite{2021Lugosi:Mendelson} attains the optimal dimension-free rate but has no known polynomial-time implementation, while efficient coordinate-wise estimators have errors that grow with the dimension. For Gaussian data, the Tukey median exhibits the same trade-off: it is minimax optimal~\cite[Thm.~2.1]{2018Chen:Gao:Ren} but NP-hard to compute~\cite{2006Bernholt}.

This trade-off has been overcome in the past decade. Starting with the seminal works~\cite{2016Diakonikolas:Kamath:Kane:Li:Moitra:Stewart,2016Lai:Rao:Vempala}, a line of research has developed polynomial-time, and even nearly linear-time~\cite{2019Cheng:Diakonikolas:Ge}, estimators whose corruption error is dimension-independent and nearly matches the information-theoretic limit for broad classes of distributions. Many of them rely on \emph{iterative reweighting}~\cite{2018Steinhardt:Charikar:Valiant,2019Dong:Hopkins:Li,2020Cheng:Diakonikolas:Ge:Soltanolkotabi,2020Hopkins:Li:Zhang,2022Zhu:Jiao:Steinhardt,2022Dalalyan:Minasyan}, a soft version of the original \emph{filtering} method~\cite{2016Diakonikolas:Kamath:Kane:Li:Moitra:Stewart,2018Diakonikolas:Kamath:Kane:Li:Moitra:Stewart,2017Diakonikolas:Kamath:Kane:Li:Moitra:Stewart,2023Diakonikolas:Kane}.

Let us formalize the optimization problem behind many of these works \cite[Section 1.1]{2020Cheng:Diakonikolas:Ge:Soltanolkotabi}. Let $0<k<n$ be an integer, the \emph{trimming hyperparameter},\footnote{By
Lemma~\ref{lemma:capped:extreme} in Section \ref{ss:capped:simplex:spechtahedron}, every $w\in\Delta_{n,k}$ is a convex combination of uniform distributions on $n-k$ points, i.e., of weightings that trim $k$ points.} let $\varepsilon_{n,k}:=k/n$ be a proxy for $\varepsilon$, and let $\Delta_n$ be the probability simplex of $\re^n$. For $\omega\in\Delta_n$, the set $\Delta_{n,k}(\omega):=\{w\in\Delta_n:w_i\le\frac{\omega_i}{1-\varepsilon_{n,k}},\ i\in[n]\}$
consists of the distributions obtained from $\omega$ by removing an
$\varepsilon_{n,k}$-fraction of its mass and renormalizing~\cite{2022Zhu:Jiao:Steinhardt},~\cite[Chap.~2]{2023Diakonikolas:Kane}.
With $\unit_{[n]}\in\re^n$ the all-ones vector,
\begin{align}
    \Delta_{n,k}:=\Delta_{n,k}(\unit_{[n]}/n)
    =\Big\{w\in\Delta_n:w_i\le\frac{1}{n-k},\ i\in[n]\Big\}.
\end{align}
This set is sometimes called the \emph{capped
simplex} in the optimization literature~\cite{2015Wang:Lu,2021Ang:Ma:Liu:Huang:Wang}. 

For $w\in\Delta_n$ and $\bnu\in\re^d$, define the weighted mean, the weighted second moment about $\bnu$, and the weighted covariance:
\begin{align}
    \bmmu_w:=\sum_{i=1}^nw_i\tilde\bx_i,~~~
    \bfSigma_w(\bnu):=\sum_{i=1}^nw_i(\tilde\bx_i-\bnu)(\tilde\bx_i-\bnu)^\top,~~~
    \bfSigma_w:=\bfSigma_w(\bmmu_w).
\end{align}
Let $\mbX\in\re^{n\times d}$ be the data matrix with rows $\tilde\bx_1^\top,\ldots,\tilde\bx_n^\top$, so that $\bmmu_w=\mbX^\top w$, and let $M_{\mbX}:=\max_{i\in[n]}\|\tilde\bx_i\|_2$.

Following the framework of~\cite{2016Diakonikolas:Kamath:Kane:Li:Moitra:Stewart}, \cite[Section~1.1]{2020Cheng:Diakonikolas:Ge:Soltanolkotabi} and follow-up works (e.g.,~\cite{2022Zhu:Jiao:Steinhardt,2020Hopkins:Li:Zhang}) approximately solve, by iterative reweighting,
\begin{align}
    \min_{w\in\Delta_{n,k}}\|\bfSigma_w\|,
    \label{problem:surrogate:spectral:noncentered}
\end{align}
where $\|\cdot\|$ is the spectral norm and, typically, $k\ge\varepsilon n$. This problem is nonconvex because the centering $\bmmu_w$ depends on $w$.

Despite this non-convexity, simple reweighting schemes yield efficient and accurate estimators, thanks to two ingredients introduced in~\cite{2016Diakonikolas:Kamath:Kane:Li:Moitra:Stewart} and refined in many follow-up works. Informally:
\begin{quote}
\textbf{Stability.} \emph{Every reweighting of the clean sample that removes at most an $\varepsilon$-fraction of its mass has a weighted mean close to $\bmmu$ in Euclidean norm and a weighted second moment about $\bmmu$ close to $\bfSigma$ in spectral norm.}
\end{quote}
Closely related notions appear under other names, such as \emph{good} sets, data points or subsamples~\cite{2016Diakonikolas:Kamath:Kane:Li:Moitra:Stewart,2019Li:lecture:notes,2019Cheng:Diakonikolas:Ge,2019Dong:Hopkins:Li,2022Zhu:Jiao:Steinhardt},
and \emph{resilience}~\cite{2018Steinhardt:Charikar:Valiant,2022Zhu:Jiao:Steinhardt:generalized_resilience},
which constrains only the weighted means~\cite[Exercise~2.2]{2023Diakonikolas:Kane}.

The second is an algorithmic certificate based on the spectral norm of the weighted covariance matrix, called a \emph{certificate for the empirical mean} in~\cite{2016Diakonikolas:Kamath:Kane:Li:Moitra:Stewart},\cite[Lemma~2.6]{2023Diakonikolas:Kane} and related to the \emph{spectral signatures} of~\cite{2017Diakonikolas:Kamath:Kane:Li:Moitra:Stewart,2019Li:lecture:notes,2019Dong:Hopkins:Li,2022Zhu:Jiao:Steinhardt:generalized_resilience,2022Zhu:Jiao:Steinhardt,2022Dalalyan:Minasyan,2020Hopkins:Li:Zhang}. See also \cite{2018Tran:Li:Madry}. Informally:
\begin{quote}
\textbf{Spectral signature.} \emph{Suppose the clean sample is stable. If a reweighting of the corrupted sample removes at most an $\varepsilon$-fraction of its mass and its weighted covariance has small spectral norm, then its weighted mean is close to $\bmmu$ in Euclidean distance.}
\end{quote}

Together, these ingredients show that approximate minimizers
of~$\|\bfSigma_w\|$ over~$\Delta_{n,k}$ yield accurate estimates of~$\bmmu$.
Algorithms exploit this by reweighting iteratively and computing at each step an approximate top eigenvector~\cite{2016Diakonikolas:Kamath:Kane:Li:Moitra:Stewart,2018Diakonikolas:Kamath:Kane:Li:Moitra:Stewart,2017Diakonikolas:Kamath:Kane:Li:Moitra:Stewart,2019Dong:Hopkins:Li,2020Cheng:Diakonikolas:Ge:Soltanolkotabi,2020Hopkins:Li:Zhang,2022Zhu:Jiao:Steinhardt,2022Dalalyan:Minasyan,2023Diakonikolas:Kane},
or solving a semidefinite
program (SDP)~\cite{2016Diakonikolas:Kamath:Kane:Li:Moitra:Stewart,2018Steinhardt:Charikar:Valiant,2019Cheng:Diakonikolas:Ge}. We discuss both ingredients in Section~\ref{s:stability:spectral:signature}.

Despite this progress, important questions remain open, to the best of our knowledge. The first concerns stability. Most of the works above guarantee their error bounds only with constant probability, rather than attaining the \emph{sub-Gaussian rates} of heavy-tailed mean estimation~\cite{2019Lugosi:Mendelson}; see the
survey~\cite{2019Lugosi:Mendelson:survey}. Diakonikolas, Kane and Pensia~\cite{2020Diakonikolas:Kane:Pensia} were among the first to establish stability with sub-Gaussian rates. Their rates, however, are not dimension-free for all $p\ge2$, and in some regimes they are only near-optimal; see Section~\ref{s:related:work} for a detailed discussion. This motivates the question:
\begin{quote}
  \textbf{Question 1.} \emph{Is there a unified proof of stability, with the optimal dimension-free sub-Gaussian rate and the optimal dependence on $\varepsilon$, namely, for given probability of failure $\delta\in(0,1)$, 
  \begin{align}
    r = r(n,\delta,\varepsilon) := \sqrt{\frac{\tr(\bfSigma)}{n}}
    +\nu_2\sqrt{\frac{\log(3/\delta)}{n}}
    +\inf_{p\ge2}\nu_p\,\varepsilon^{1-\frac1p},\label{eq:subgaussian:rate}
  \end{align}
    valid for every $p\ge2$ with $\nu_p<\infty$?}
\end{quote}

Relatedly, current stability-based reweighting algorithms typically require knowledge of $(p,\nu_p)$~\cite{2020Diakonikolas:Kane:Pensia} and are often tailored to a distribution class fixed in
advance~\cite{2019Cheng:Diakonikolas:Ge,2020Hopkins:Li:Zhang}. This motivates the following question.
\begin{quote}
  \textbf{Question 2.} \emph{Is there an efficient stability-based estimator that, without knowledge of $(p,\nu_p)$, attains the optimal rate $r(n,\delta,\varepsilon)$ for every $p\ge2$ with $\nu_p<\infty$?}
\end{quote}
Regarding this question, the polynomial-time estimators known to attain $r(n,\delta,\varepsilon)$ for $p>2$ require $\bfSigma$ to be known or isotropic, or to be approximated a priori by some $\hat\bfSigma$ with $\|\hat\bfSigma-\bfSigma\|$ small \cite{2022Dalalyan:Minasyan,2020Diakonikolas:Kane:Pensia,2020Hopkins:Li:Zhang,2025Bateni:Minasyan:Dalalyan}. In fact, these previous results are somewhat weaker than indicated in the preceding sentence; see Section \ref{s:related:work} for details. Sum-of-squares methods achieve optimal $r(n,\delta,\varepsilon)$ at the cost of the stronger condition of \emph{certifiable hypercontractivity} and $n\gtrsim d^{O(p)}$ samples~\cite{2018Kothari:Steinhardt:Steurer,2018Hopkins:Li,2022Kothari:Manohar:Zhang}. For general unknown $\bfSigma$, filtering and reweighting methods attain only $r(n,\delta,\varepsilon)+\sqrt{\varepsilon\|\bfSigma\|}$, which is optimal only for $p=2$, requiring $n\gtrsim d$. Statistical query lower bounds~\cite{2017Diakonikolas:Kane:Stewart}, \cite[Exercise~8.8]{2023Diakonikolas:Kane} suggest that no efficient statistical query algorithm with sample size subquadratic in $d$ can improve on this rate, even for distributions with bounded higher moments. There is also worst-case complexity-theoretic evidence that improving on this rate for distributions with bounded higher moments may be computationally intractable~\cite{2019Hopkins:Li}.

A second point is the practical appeal of simple iterative methods. The seminal works~\cite{2016Diakonikolas:Kamath:Kane:Li:Moitra:Stewart,2016Lai:Rao:Vempala} relied on convex relaxations or intricate iterative spectral procedures, which motivated simpler filtering methods~\cite{2023Diakonikolas:Kane} based on iterative reweighting, analyzed either through regret minimization, using multiplicative weights~\cite{2012Arora:Hazan:Kale} over the capped simplex~\cite{2019Dong:Hopkins:Li,2020Hopkins:Li:Zhang,2022Zhu:Jiao:Steinhardt},  or as first-order methods on a non-convex objective~\cite{2020Cheng:Diakonikolas:Ge:Soltanolkotabi,2022Zhu:Jiao:Steinhardt}. As argued in~\cite{2020Cheng:Diakonikolas:Ge:Soltanolkotabi}, such methods deserve study in their own right, being easy to implement and widely used in practice. We are interested in the statistical precision of iterative reweighting with only an approximate top eigenvector per step, rather than a semidefinite program or a full eigendecomposition~\cite{2016Lai:Rao:Vempala,2025Bateni:Minasyan:Dalalyan}.

For first-order methods, the relevant object is the \emph{optimization landscape} of robust mean estimation. Cheng, Diakonikolas, Ge and Soltanolkotabi~\cite{2020Cheng:Diakonikolas:Ge:Soltanolkotabi} showed that, despite non-convexity, any approximate stationary point of~\eqref{problem:surrogate:spectral:noncentered} yields an accurate estimate of $\bmmu$. Zhu, Jiao and Steinhardt~\cite{2022Zhu:Jiao:Steinhardt} simplified this analysis via quasi-gradients~\cite{2015Hazan:Levy:ShalevShwartz}, showing that such points are also approximate global minimizers~\cite[Thm.~2.1 and Lemma~2.2]{2022Zhu:Jiao:Steinhardt}, but only for $\varepsilon<1/3$: for $\varepsilon\ge1/3$, stationary points can be arbitrarily far from optimal~\cite[Thm.~1.1]{2022Zhu:Jiao:Steinhardt}. This motivates the following question.
\begin{quote}
  \textbf{Question 3.} \emph{Is there a precise quantitative description of the optimization landscape of robust mean estimation that bounds the error $\|\bmmu_w-\bmmu\|_2$ by standard optimization metrics, such as the optimality gap? Can such a description guide the design and analysis of efficient robust mean estimators?}
\end{quote}

\subsection{Main results}\label{s:main:results}

We present our results along the questions of the introduction: stability with optimal sub-Gaussian rates for all $p\ge2$ (Questions~1--2), motivated by~\cite{2020Diakonikolas:Kane:Pensia,2020Hopkins:Li:Zhang}, and the optimization landscape and first-order methods (Question~3), motivated by~\cite{2020Cheng:Diakonikolas:Ge:Soltanolkotabi,2022Zhu:Jiao:Steinhardt,2020Hopkins:Li:Zhang}.

\paragraph{A hidden convex problem.}
Instead of solving~\eqref{problem:surrogate:spectral:noncentered}, we target the unobserved convex problem
\begin{align}
  F^*:=\min_{w\in\Delta_{n,k}}\Big\{F(w):=\lambda_1(\bfSigma_w(\bmmu)-\bfSigma)\Big\},
  \label{eq:landscape:problem}
\end{align}
where $\lambda_1(\cdot)$ is the largest eigenvalue. $F$ is convex since $w\mapsto\bfSigma_w(\bmmu)$ is linear, but depends on the unknown $(\bmmu,\bfSigma)$. Replacing $\bmmu$ by $\bmmu_w$ and $\bfSigma$ by a symmetric positive semidefinite surrogate $\hat\bfSigma$ gives the observable function
\begin{align}
  \hat F(w) := \lambda_1(\bfSigma_w-\hat\bfSigma),\qquad w\in\Delta_{n,k},
  \label{eq:surrogate:spectral:function}
\end{align}
which equals the objective $\|\bfSigma_w\|$ of~\eqref{problem:surrogate:spectral:noncentered} when $\hat\bfSigma=\bzero$. Unlike~\cite{2020Cheng:Diakonikolas:Ge:Soltanolkotabi,2022Zhu:Jiao:Steinhardt,2020Hopkins:Li:Zhang}, which treat the objective $\|\bfSigma_w\|$ of~\eqref{problem:surrogate:spectral:noncentered}, we use $\hat F$ only as an inexact subgradient oracle for $F$. The surrogate enters our bounds only through $\tau:=\lambda_1(\hat\bfSigma-\bfSigma)+\lambda_1(\bfSigma-\hat\bfSigma)\ge0$, which vanishes when $\hat\bfSigma-\bfSigma$ is a multiple of $\bfI_d$.

Two obstacles arise: $F$ is unobserved, and $\hat F$ is nonconvex. $\hat F$ is weakly convex~\cite{1998Rockafellar:Wets}, but with a Euclidean modulus which grows with $n$, $d$ and the size of the outliers~\cite{2020Cheng:Diakonikolas:Ge:Soltanolkotabi}. We overcome both by using \emph{relative} regularity notions from optimization~\cite{2017Bauschke:Bolte:Teboulle,2018Lu:Freund:Nesterov,2019Lu}, measured in the (regularized) \emph{mean-induced divergence}
\begin{align}
  D_{\mbX,\alpha}(\omega\|w) := \tfrac12\|\mbX^\top(\omega-w)\|_2^2+\tfrac\alpha2\|\omega-w\|_2^2,\qquad \alpha\ge0.\label{eq:mean:induced:div}
\end{align}
The name reflects the identity $\mbX^\top(\omega-w)=\bmmu_\omega-\bmmu_w$.

\paragraph{Stability and landscape.}
Given any $\delta\in(0,1)$ and $\eta\in(0,1]$, let $k:=\lceil(1+\eta)\varepsilon n\rceil+\lceil\log(3/\delta)\rceil$. Our bounds hold up to factors depending only on $\eta$. 

\begin{theorem}[Stability and landscape, informal]\label{thm:informal:landscape}
Under Assumption~\ref{assump:finite:covariance:contamination}, if $\varepsilon_{n,k}=\frac{\lceil(1+\eta)\varepsilon n\rceil+\lceil\log(3/\delta)\rceil}{n}<\frac13$, then with probability at least $1-\delta$:
\begin{enumerate}
  \item[\rm (a)] the sample is \emph{favorable} (Definition~\ref{def:favorable} in Section \ref{s:stability:spectral:signature}): the spectral signature holds; any minimizer $w^*$ of~\eqref{eq:landscape:problem} is a good weight, i.e., $F(w^*)\lesssim r^2/\varepsilon_{n,k}$; and $F\gtrsim -r^2/\varepsilon_{n,k}$ uniformly on $\Delta_{n,k}$. In particular, the sample is \emph{stable} with radius $\lesssim r$;
  \item[\rm (b)] for every $\alpha\ge0$, $F$ has \emph{inexact quadratic growth} on $\Delta_{n,k}$ relative to $D_{\mbX,\alpha}$ (Definition \ref{def:relative:iqg} in Section \ref{ss:relative:regularity}), with curvature $m=\frac{1-\varepsilon_{n,k}}{\varepsilon_{n,k}}$ and radius $\lesssim r+\sqrt\alpha\,\frac{\sqrt{2k}}{n-k}$. For $\alpha=0$ or $\alpha\asymp\nu_2^2/\varepsilon_{n,k}$, the radius is $\lesssim r$. 
  \item[\rm (c)] $\partial\hat F$ is an inexact \emph{weakly star-convex} oracle for $F$ on $\Delta_{n,k}$ relative to $D_{\mbX,0}$, with dimension-free parameters $(4,O(r),\tau)$, and a \emph{Lipschitz} inexact oracle relative to the Kullback--Leibler divergence $D_{\KL}$, the Euclidean divergence $D_{\E}$ and $D_{\mbX,\alpha}$ with Lipschitz constants  $\frac52M_{\mbX}^2$, $\frac52\sqrt n\,M_{\mbX}^2$ and $M_{\mbX}\sqrt{4+nM_{\mbX}^2/(4\alpha)}$ for $\alpha>0$, respectively (Definition \ref{def:inexact:oracle} in Section \ref{ss:relative:regularity}). 
\end{enumerate}
\end{theorem}
For the formal statement, see Theorem \ref{thm:favorability:subgaussian:rate} in Section \ref{ss:favorable}, Theorem \ref{thm:stability:subgaussian:rate} in Section \ref{ss:stability} and Propositions \ref{prop:IQG}-\ref{prop:weakly:convex:oracle} and Lemma \ref{lemma:lipschitz} in Section \ref{s:landscape}.

Item~(a) gives a unified proof of stability for every $p\ge2$ (Question~1). The uniform lower bound in~(a) is, to our knowledge, new in algorithmic robust statistics, reflecting that lower tails of nonnegative processes are easier to control than upper tails. Items~(b)--(c) answer Question~3: combined, they yield the Lyapunov dissipation function of Definition~\ref{def:lyapunov:Phi:m:r} of Section \ref{s:mirror:descent}, which drives all convergence proofs;  see Theorem~\ref{thm:informal:rates}.

\paragraph{Algorithms.}
Both families below use the score $s^t\in\re^n$,
\begin{align}
  s_i^t:=\langle\bv^t,[(\tilde\bx_i-\bmmu_{w^t})^{\otimes2}-\hat\bfSigma]\bv^t\rangle~~~(i\in[n]),\label{eq:score}
\end{align}
where $\bv^t$ is a top unit eigenvector of $\bfSigma_{w^t}-\hat\bfSigma$, so each step needs one top-eigenvector computation. A direct computation gives $s^t+c^t\unit_{[n]}\in\partial\hat F(w^t)$ with $c^t:=\langle\bv^t,\hat\bfSigma\bv^t\rangle-\langle\bv^t,\bmmu_{w^t}\rangle^2$. Both start at $w^1=\unit_{[n]}/n$, with step sizes $\gamma_t>0$, and neither uses $(p,\nu_p)$.

\emph{Mirror descent} with reference function $h$ (Section~\ref{s:bregman:geometry}) computes
\begin{align}
    w^{t+1}:=\argmin_{w\in\Delta_{n,k}}
    \Big\{\langle g^t,w\rangle+\gamma_t^{-1}D_h(w\|w^t)\Big\},
    \qquad g^t\in\partial\hat F(w^t);
    \label{eq:md:update}
\end{align}
since $\langle\unit_{[n]},w\rangle=1$ on $\Delta_n$, one may take $g^t=s^t$. We consider $D_{\KL}$ (entropic mirror descent, or Hedge), $D_{\E}$ (the projected subgradient method) and $D_{\mbX,\alpha}$ with $\alpha>0$ (mean-induced mirror descent). The last is matched to the quadratic growth of Theorem~\ref{thm:informal:landscape}(b), which yields faster rates (Table~\ref{tab:rates}).

The \emph{multiplicative weights update} (MWU), a variant of~\cite{2020Hopkins:Li:Zhang} with costs centered at $\hat\bfSigma$, computes
\begin{align}
    \bar w^{t+1}_i:=\frac{w_i^t\big(1-\frac{\gamma_t}{{M_\infty^2}} s_i^t\big)}
    {\sum_{j=1}^nw_j^t\big(1-\frac{\gamma_t}{{M_\infty^2}}s_j^t\big)},\quad i\in[n],
    \qquad
    w^{t+1}:=\argmin_{w\in\Delta_{n,k}}D_{\KL}(w\|\bar w^{t+1}),
    \label{eq:mwu}
\end{align}
with $\gamma_t\in(0,1]$ and ${M_\infty^2}:=2\max\{4M_{\mbX}^2,\|\hat\bfSigma\|\}\ge2\sup_t\|s^t\|_\infty$, so that $\frac{\gamma_t}{{M_\infty^2}}|s_i^t|\le\frac12$. Unlike mirror descent, MWU is not invariant to shifts of $s^t$ by multiples of $\unit_{[n]}$.

\paragraph{Statistical and optimization rates.} The next theorem summarizes the statistical and optimization rates of our four algorithms. For the formal statement see Corollaries \ref{cor:rates:md}-\ref{cor:rates:mwu} in Section \ref{s:stat:opt:rates}.

\begin{theorem}[Statistical and optimization rates, informal.]\label{thm:informal:rates}
On the event of Theorem~\ref{thm:informal:landscape}, the averaged iterates of all four methods (with $\alpha\lesssim\nu_2^2/\varepsilon_{n,k}$ for mean-induced mirror descent) satisfy the bounds of Table~\ref{tab:rates}. In particular, if $\varepsilon_{n,k}$ is bounded away from $\frac13$, their statistical error is of order $r(n,\delta,\varepsilon)+\sqrt{\varepsilon_{n,k}\,\tau}$, attained after the number of iterations in Table~\ref{tab:rates} and without knowledge of $(p,\nu_p)$. The surrogate $\hat\bfSigma$ enters only through $\tau$, which vanishes whenever $\hat\bfSigma-\bfSigma$ is a multiple of $\bfI_d$. Hence the optimal rate $r(n,\delta,\varepsilon)$ of Question~2 is attained for $\hat\bfSigma=\bfSigma$ (known covariance) and for $\hat\bfSigma=\bzero$ when $\bfSigma=\sigma^2\bfI_d$ (isotropic, $\sigma$ unknown), while for general unknown $\bfSigma$ and $\hat\bfSigma=\bzero$, $\tau\le\|\bfSigma\|$ gives $r(n,\delta,\varepsilon)+\sqrt{\varepsilon_{n,k}\|\bfSigma\|}$.
\end{theorem}
Theorem \ref{thm:informal:rates} answers Question 2 in that our algorithms are adaptive to $(p,\nu_p)$; see Section \ref{s:related:work} for a detailed comparision with prior work. The last iterate of mean-induced mirror descent obeys the same bounds, except that its surrogate term is the larger $\sqrt\tau\,\varepsilon_{n,k}^{1/4}$ (Table~\ref{tab:rates}). The condition $\varepsilon_{n,k}<\frac13$ is sharp for our Lyapunov dissipation argument: it is exactly $m=\frac{1-\varepsilon_{n,k}}{\varepsilon_{n,k}} > 2$, and it matches the threshold of~\cite{2022Zhu:Jiao:Steinhardt}. We have observed this threshold in our numerical experiments for all four methods; see Section \ref{s:simulations}.

\begin{table}[h]
\centering
\caption{Rates for $\varepsilon_{n,k}\asymp\varepsilon+\frac{\log(3/\delta)}{n}$ bounded away from $\frac13$ (MD: mirror descent). Here $r=r(n,\delta,\varepsilon)$ is the rate~\eqref{eq:subgaussian:rate}; for the mean-induced geometry, $\alpha\asymp\nu_2^2/\varepsilon_{n,k}$ and $M\asymp M_{\mbX}(1+\sqrt k\,M_{\mbX}/\nu_2)$. For all averaged iterates, the statistical part of the spectral bound is of order $r^2/\varepsilon_{n,k}+\tau$.}
\label{tab:rates}
\small
\setlength{\tabcolsep}{4pt}
\begin{tabular}{@{}lccc@{}}
\toprule
 & \multicolumn{2}{c}{Mean error} & Spectral gap\\
\cmidrule(lr){2-3}\cmidrule(l){4-4}
Method & statistical & optimization & optimization\\
\midrule
Entropic/Euclidean MD & $r+\sqrt{\tau\varepsilon_{n,k}}$ & $M_{\mbX}\,\varepsilon_{n,k}^{3/4}T^{-1/4}$ & $M_{\mbX}^2\sqrt{\varepsilon_{n,k}/T}$\\
Mean-ind.\ MD (window average) & $r+\sqrt{\tau\varepsilon_{n,k}}$ & $M\varepsilon_{n,k}^{3/4}T^{-1/2}$ & $M^2\sqrt{\varepsilon_{n,k}}\,T^{-1}$\\
Mean-ind.\ MD (last iterate) & $r+\sqrt\tau\,\varepsilon_{n,k}^{1/4}$ & $M\sqrt{\varepsilon_{n,k}}\,T^{-1/2}$ & ---\\
MWU & $r+\sqrt{\tau\varepsilon_{n,k}}$ & ${M_\infty}\,\varepsilon_{n,k}^{3/4}T^{-1/4}$ & ${M_\infty^2}\sqrt{\varepsilon_{n,k}/T}$\\
\bottomrule
\end{tabular}
\end{table}

Entropic and Euclidean MD and MWU attain the rates $T^{-1/4}$ for the mean and $T^{-1/2}$ for the spectral gap, with the same dependence on $\varepsilon_{n,k}$ and on $M_{\mbX}$ (or $\sqrt{{M_\infty^2}}$).\footnote{Although the Lipschitz constant and the Bregman diameter of $\Delta_{n,k}$ differ between entropic and Euclidean MD, the product of the former with the square root of the latter is of order $M_{\mbX}^2\sqrt{\varepsilon_{n,k}}$ in both, so the rates coincide: on $\Delta_{n,k}$, the $\{n/\log n\}^{1/2}$ advantage usually attributed to the entropic geometry on $\Delta_n$ disappears; see Corollary \ref{cor:rates:md} in Section \ref{s:stat:opt:rates}.} Mean-induced MD improves these to $T^{-1/2}$ and $T^{-1}$ and also controls the last iterate, at three costs: $M$ grows with $\sqrt k/\nu_2$; the tuning $\alpha\asymp\nu_2^2/\varepsilon_{n,k}$ requires knowledge of $\nu_2$; and each step solves a strongly convex quadratic program over $\Delta_{n,k}$ with matrix $\mbX\mbX^\top+\alpha\bfI_n$. Its mean optimization error is smaller than that of entropic MD once $T\gtrsim(1+\sqrt k\,M_{\mbX}/\nu_2)^4$.

\begin{remark}[Approximate top eigenvector computation]\label{rem:approx:eigenvector}
  Suppose that, for all $t\in\mathbb{N}$, $\hat\bfSigma$ is such that $\bfSigma_{w^t}-\hat\bfSigma\succeq\mathbf{0}$ and $\bv^t\in\mbS^{d-1}$ is only a $\chi$-approximate top eigenvector, i.e., $\langle\bv^t,[\bfSigma_{w^t}-\hat\bfSigma]\bv^t\rangle\ge(1-\chi)\lambda_1(\bfSigma_{w^t}-\hat\bfSigma)$ for some $\chi\in[0,1)$. This $\bv^t$ can be computed, e.g., by a randomized power method. Define $\tau_\chi:=\tau+\chi(\lambda_d(\bfSigma-\hat\bfSigma))_+ + O(\chi)\,r^2/\varepsilon_{n,k}$. In that case:
  \begin{itemize}
    \item the required threshold is $\varepsilon_{n,k}<\frac{1-\chi}{3-2\chi}$. We do not know if this dependence on $\chi$ is an artifact of the analysis;
    \item in Table~\ref{tab:rates} and Theorem~\ref{thm:informal:rates}, $\tau$ is replaced by $\tau_\chi$; the spectral gap bounds acquire an overall factor $(1-\chi)^{-1}$ and the mean-error bounds a factor $(1-\chi)^{-1/2}$. The last iterate of mean-induced mirror descent picks up only an absolute constant;
    \item The term $\chi\,r_1^2/\varepsilon_{n,k}$ in $\tau_\chi$ only inflates $r$ by a factor $1+\sqrt\chi$; for $\hat\bfSigma=\bzero$, the bound $r+\sqrt{\varepsilon_{n,k}\|\bfSigma\|}$ for general unknown $\bfSigma$ persists, while for $\bfSigma=\sigma^2\bfI_d$ (where $\tau=0$) the error is $r+\sigma\sqrt{\chi\,\varepsilon_{n,k}}$, which is $\lesssim r$ as soon as $\chi\lesssim r^2/(\varepsilon_{n,k}\nu_2^2)$; for $\hat\bfSigma=\bfSigma$ (known covariance), $\tau_\chi=O(\chi)\,r^2/\varepsilon_{n,k}$ and the surrogate term $\sqrt{\varepsilon_{n,k}\tau_\chi}=\sqrt{\chi}\,O(r)$ is of the same order as the statistical term.
  \end{itemize}
See Remark \ref{rem:approx:eigenvector:details} in Section \ref{s:stat:opt:rates} for further details. 
\end{remark}

The constants $M_{\mbX}=\max_i\|\tilde\bx_i\|_2$ and ${M_\infty}$ may be arbitrarily large under contamination, but they affect only the optimization errors. As observed in~\cite[Section~4.3.2]{2022Zhu:Jiao:Steinhardt}, via a regret analysis with quasi-subgradients, no pruning is needed for statistical validity.  This helps in practice when the outliers do not inflate $M_{\mbX}$ much. In the worst case $M_{\mbX}$ may be exponential in the input size. A standard pruning step~\cite{2020Hopkins:Li:Zhang,2020Diakonikolas:Kane:Pensia} restores a polynomial-time guarantee and reduces the optimization errors of all methods, at the cost of increasing the contamination rate.

\begin{remark}[From a high-probability landscape to statistical precision]
At a high level, our analysis follows a template common in high-dimensional estimation: with high probability, the sample satisfies deterministic conditions; these give the optimization problem a favorable landscape; and this landscape drives a first-order method to statistical precision. Our contributions are threefold.
\begin{enumerate}
  \item \emph{Probability.} Counting and truncated empirical processes have recently emerged as a central tool for analyzing trimmed-mean estimators \cite{2019Oliveira:Orenstein,2021Lugosi:Mendelson,2024Oliveira:Rico,2026Abdalla:Zhivotovskiy,2025Oliveira:Resende}, and they also underlie algorithmic robust estimation (see, e.g., \cite{2020Diakonikolas:Kane:Pensia}). We sharpen the analysis of these processes and connect it to stability conditions, spectral signatures, and reweighting algorithms. See Sections \ref{s:robust:averages:counting}-\ref{s:stability:spectral:signature}.
  \item \emph{Landscape.} Building on these probabilistic results, we establish improved stability and spectral signature conditions, which we call \emph{favorability}, at the optimal sub-Gaussian rate $r(n,\delta,\varepsilon)$ for every $p\ge2$. We further show that, whenever $\varepsilon_{n,k}<\tfrac13$, a subgradient surrogate computable from the data yields weaker, inexact versions of the standard conditions for convergence of first-order methods. These conditions hold with respect to the geometry on the capped simplex induced by weighted means, with dimension-free parameters. See Sections \ref{s:stability:spectral:signature}-\ref{s:landscape}.
  \item \emph{Optimization.} We show that these weaker conditions suffice for mirror descent and MWU, run over the capped simplex, to converge down to the statistical precision whenever $\varepsilon_{n,k}<\tfrac13$. This is in the spirit of results for several otherwise quite different high-dimensional estimation problems (in several cases non-convex), where benign local or restricted geometry, rather than global strong convexity, drives convergence. Examples include sparse regression \cite{2012Negahban:RavikumarWainwright:Yu}, low-rank recovery \cite{2021Charisopoulos:Chen:Davis:Diaz:Ding:Drusvyatskiy}, phase retrieval \cite{2019Duchi:Ruan,2020Davis:Drusvyatskiy:Paquette}, and blind deconvolution \cite{2021Charisopoulos:Davis:Diaz:Drusvyatskiy}. See Sections \ref{s:landscape}-\ref{s:stat:opt:rates}. 
\end{enumerate}
\end{remark}

\subsection{Related work and organization}\label{s:related:work}

\paragraph{Stability and statistical error with optimal sub-Gaussian rate.} Following~\cite{2016Diakonikolas:Kamath:Kane:Li:Moitra:Stewart,2017Diakonikolas:Kamath:Kane:Li:Moitra:Stewart,2018Diakonikolas:Kamath:Kane:Li:Moitra:Stewart,2018Steinhardt:Charikar:Valiant}, Diakonikolas, Kane and Pensia~\cite{2020Diakonikolas:Kane:Pensia} were among the first to prove stability with sub-Gaussian rates; see also~\cite{2022Zhu:Jiao:Steinhardt:generalized_resilience}. For $p=2$, their rates are dimension-free and near-optimal, up to a factor $\sqrt{\log\erank(\bfSigma)}$ in $\sqrt{\tr(\bfSigma)/n}$ that is removed after median-of-means (MoM) preprocessing. For $p\ge4$, they require isotropy, have a suboptimal $\sqrt{d\log d/n}$ term instead of $\sqrt{d/n} = \sqrt{\tr(\bfSigma)/n}$, and carry $\nu_4$ instead of $\nu_2$ in front of $\sqrt{\log(1/\delta)/n}$. Hopkins, Li and Zhang~\cite{2020Hopkins:Li:Zhang} relate stability for $p=2$ to MoM estimators~\cite{2019Lugosi:Mendelson}, and optimal dimension-free rates are known for Gaussian~\cite{2022Dalalyan:Minasyan} and, up to logarithmic factors, sub-Gaussian~\cite{2025Bateni:Minasyan:Dalalyan} distributions. However, it is known that MoM does not achieve optimal rates in terms of $(p,\nu_p)$ for $p>2$ \cite{2025deJuan:Mazuelas}, and existing proofs treat ranges of $p$ separately, leave $2<p<4$ uncovered, and lose logarithmic factors without MoM. These limitations carry over to the statistical rates of stability-based algorithms~\cite{2020Diakonikolas:Kane:Pensia,2020Hopkins:Li:Zhang,2022Dalalyan:Minasyan,2025Bateni:Minasyan:Dalalyan}.

We prove stability for every $p\ge2$ with the optimal rate $r(n,\delta,\varepsilon)$, using techniques that connect uniform control of trimmed means with the concentration of counting and truncation processes~\cite{2019Oliveira:Orenstein,2021Lugosi:Mendelson,2026Abdalla:Zhivotovskiy,2024Oliveira:Rico}; these arguments also avoid the randomized rounding steps of~\cite{2022Depersin:Lecue,2020Hopkins:Li:Zhang} (Sections~\ref{s:robust:averages:counting}--\ref{s:counting:truncation:in:probability}). Combined with Theorem~\ref{thm:informal:landscape}, they connect the outlier-robust and heavy-tailed literatures~\cite{2020Hopkins:Li:Zhang} through algorithmic reweighting, and yield, to our knowledge, the first reweighting algorithms that attain, for all $p\ge2$ and adaptively to $(p,\nu_p)$, the rate $r(n,\delta,\varepsilon)$ when $\bfSigma$ is known, or unknown but isotropic, and $r(n,\delta,\varepsilon)+\sqrt{\tau\varepsilon}$ otherwise.

\paragraph{Optimization landscape and reweighting algorithms.} Among reweighting methods that avoid large SDPs and full SVDs, the closest to ours are the original filtering algorithms~\cite{2016Diakonikolas:Kamath:Kane:Li:Moitra:Stewart,2017Diakonikolas:Kamath:Kane:Li:Moitra:Stewart,2023Diakonikolas:Kane}; the related method of~\cite{2018Steinhardt:Charikar:Valiant} also reweights iteratively, but solves an SDP at each step. Dong, Hopkins and Li~\cite{2019Dong:Hopkins:Li} analyzed matrix multiplicative weights via regret; with sketched approximate SVDs, theirs remains the fastest method, running in time $\tilde O(nd)$ for $p=2$, but it is suboptimal in terms of $(p,\nu_p)$ for $p>2$~\cite{2020Hopkins:Li:Zhang}. Cheng, Diakonikolas, Ge and Soltanolkotabi~\cite{2020Cheng:Diakonikolas:Ge:Soltanolkotabi} first cast reweighting as a projected subgradient method, and Zhu, Jiao and Steinhardt~\cite{2022Zhu:Jiao:Steinhardt} analyzed it via regret with quasi-subgradient scores. Hopkins, Li and Zhang~\cite{2020Hopkins:Li:Zhang} simplified both the regret analysis of MWU and the projected-subgradient analysis, refining~\cite{2019Dong:Hopkins:Li,2020Cheng:Diakonikolas:Ge:Soltanolkotabi}. None of these methods is adaptive to $(p,\nu_p)$, and their runtime guarantees rely on knowing it, so we do not compare runtimes. To our knowledge, our mirror descent analysis, which builds on the relative optimization landscape of Theorem~\ref{thm:informal:landscape} and covers Hedge and mean-induced mirror descent, is new in this context. More broadly, the underlying convergence theory, for convex objectives with inexact relative quadratic growth and inexact relative subgradient oracles, valid when the curvature satisfies $m>2$, also appears to be new in optimization.

\paragraph{Additional related work.} We also review work outside the scope of the questions raised in the introduction. A large body of work in robust statistics addresses heavy-tailed mean estimation; see, for instance,~\cite{2019Lugosi:Mendelson} and the survey~\cite{2019Lugosi:Mendelson:survey}. This line of work has motivated efficient estimators based on approaches other than iterative reweighting, including~\cite{2020Hopkins,2019Cherapanamjeri:Flammarion:Bartlett,2022Depersin:Lecue,2020Lei:Luh:Venkat:Zhang}. For completeness, we also mention estimators based on the sum-of-squares method~\cite{2018Kothari:Steinhardt:Steurer,2018Hopkins:Li,2022Kothari:Manohar:Zhang}, which are beyond our scope since they require solving large SDPs.

\paragraph{Organization.} Section~\ref{s:preliminaries} collects notation and tools on the capped simplex, the spectrahedron, Bregman geometry and relative regularity. Sections~\ref{s:robust:averages:counting} and~\ref{s:counting:truncation:in:probability} relate reweighting, trimming and truncation, and prove the concentration bounds for counting and truncation processes. Section~\ref{s:stability:spectral:signature} introduces favorable samples and derives the spectral signature and stability with the optimal sub-Gaussian rate. Section~\ref{s:landscape} establishes the relative optimization landscape of Theorem~\ref{thm:informal:landscape}. Sections~\ref{s:mirror:descent} and~\ref{s:mwu} develop the convergence theory for mirror descent and the multiplicative weights update, and Section~\ref{s:stat:opt:rates} combines these results into the statistical and optimization rates of Theorem~\ref{thm:informal:rates}. Section~\ref{s:simulations} reports numerical experiments, and the Appendix contains the remaining proofs.

\section{Preliminaries}\label{s:preliminaries}

This section fixes notation and collects the tools used throughout: the capped simplex and the spectrahedron (Section~\ref{ss:capped:simplex:spechtahedron}), Bregman geometry (Section~\ref{s:bregman:geometry}), and relative regularity of inexact oracles and quadratic growth (Section~\ref{ss:relative:regularity}). Proofs not given in the text are in the Appendix.

\paragraph{Additional notation.}
For $S\subseteq[n]$, $\unit_S\in\re^n$ is its indicator vector. For a set $C$, $\co(C)$, $\ri(C)$ and $\ext(C)$ denote its convex hull, relative interior and extreme points. For $a,b\in\re$, $a\wedge b:=\min\{a,b\}$, $a\vee b:=\max\{a,b\}$ and $a_+:=a\vee0$. For matrices, $\langle\bfM,\bfW\rangle:=\tr(\bfM^\top\bfW)$, and $\|\bfM\|_F$ and $\|\bfM\|$ are the Frobenius and spectral norms, both equal to $\|\cdot\|_2$ on vectors. We write $\bfI_n$ for the identity, $\mbS^{d-1}$ for the unit sphere of $\re^d$, $\calS^d\supseteq\calS_+^d$ for the symmetric and positive semidefinite $d\times d$ matrices, $\bu^{\otimes2}:=\bu\bu^\top$, and $\lambda_1(\bfA)\ge\cdots\ge\lambda_d(\bfA)$ for the eigenvalues of $\bfA\in\calS^d$. For $w,\omega\in\Delta_n$, $\mathrm{TV}(w,\omega):=\max_{A\subseteq[n]}|\sum_{i\in A}(w_i-\omega_i)|=\sum_i(w_i-\omega_i)_+=\frac12\|w-\omega\|_1$. The \emph{spectrahedron}, or set of $d\times d$ \emph{density matrices}, is $\calM_d:=\{\bfM\in\calS_+^d:\tr(\bfM)=1\}$, a standard tool in algorithmic robust statistics~\cite{2018Steinhardt:Charikar:Valiant,2019Cheng:Diakonikolas:Ge,2020Diakonikolas:Kane:Pensia,2020Hopkins:Li:Zhang}, online learning~\cite{2012Arora:Hazan:Kale} and convex optimization~\cite{2004Nemirovski,2005BenTal:Nemirovski,2011Juditsky:Nemirovski:Tauvel}. For $\bfM\in\calM_d$, $\bfC_{\bfM}\in\calS_+^d$ denotes its positive semidefinite square root. Finally, for $(w,\bfM)\in\Delta_n\times\calM_d$, let
$f(w,\bfM):=\langle\bfSigma_w(\bmmu)-\bfSigma,\bfM\rangle$ and $\hat f(w,\bfM):=\langle\bfSigma_w-\hat\bfSigma,\bfM\rangle$. By Lemma~\ref{lemma:calM:d:convex:hull} in the following section, $F(w)=\max_{\bfM\in\calM_d}f(w,\bfM)$ and $\hat F(w)=\max_{\bfM\in\calM_d}\hat f(w,\bfM)$.

\subsection{The capped simplex and the spectrahedron}\label{ss:capped:simplex:spechtahedron}

This subsection collects mostly standard facts about the capped simplex $\Delta_{n,k}$ and the spectrahedron $\calM_d$ used in the landscape analysis. Lemma~\ref{lemma:capped:extreme} is folklore~\cite{2016Diakonikolas:Kamath:Kane:Li:Moitra:Stewart,2019Diakonikolas:Kamath:Kane:Li:Moitra:Stewart,2019Cheng:Diakonikolas:Ge}, proved in~\cite[Lemma~1]{2022Dalalyan:Minasyan}; others appear, explicitly or implicitly, in the literature. Lemma \ref{lemma:capped:diameter} is included in case of the Euclidean divergence. Except for Lemma \ref{lemma:capped:extreme}, we prove them in Appendices~\ref{app:capped:tv}--\ref{app:calM:d:convex:hull} and~\ref{ss:capped:KL} for completeness. 

\begin{lemma}[Vertices]\label{lemma:capped:extreme}
The extreme points of $\Delta_{n,k}$ are the uniform distributions on $n-k$ points, and $\Delta_{n,k}$ is their convex hull. Precisely, 
$
  \ext\Delta_{n,k}=\{\frac{\unit_S}{n-k}:S\subseteq[n],\ |S|=n-k\}
$
and 
$
  \Delta_{n,k}=\co\ext\Delta_{n,k}.
$
\end{lemma}

\begin{lemma}[Closeness to uniform]\label{lemma:capped:tv}
For every $w\in\Delta_{n,k}$, $\mathrm{TV}(w,\unit_{[n]}/n)\le\varepsilon_{n,k}$,
with equality at the vertices $w=\unit_S/(n-k)$.
\end{lemma}

\begin{lemma}[Weighted versus uniform averages]\label{lemma:weight:dif}
For every $w\in\Delta_{n,k}$ and $a\in\re^n$,
$
    \sum_{i=1}^n(w_i-\frac1n)a_i
    \le\varepsilon_{n,k}(\max_ia_i-\min_ia_i)
    \le2\,\varepsilon_{n,k}\max_i|a_i|,
$
and, if $a\ge0$, the left-hand side is at most $\varepsilon_{n,k}\max_ia_i$.
\end{lemma}

\begin{lemma}[Euclidean diameter of $\Delta_{n,k}$]\label{lemma:capped:diameter}
For integers $0\le k<n$, we have 
$
  \max_{w,w'\in\Delta_{n,k}}\|w-w'\|_2=\frac{\sqrt{2\min\{k,\,n-k\}}}{n-k}
$
and
$
\max_{w\in\Delta_{n,k}}\|w-\nicefrac{\unit_{[n]}}{n}\|_2 = \sqrt{\nicefrac{k}{n(n-k)}}.
$
\end{lemma}

\begin{lemma}[Entropic diameter of $\Delta_{n,k}$]\label{lemma:capped:kl}
For every $w\in\Delta_{n,k}$,
$
    D_{\KL}(w\|\nicefrac{\unit_{[n]}}{n})
    \le \log\frac{1}{1-\varepsilon_{n,k}}
    \le\frac{\varepsilon_{n,k}}{1-\varepsilon_{n,k}},
$
with equality in the first inequality at the vertices $w=\unit_S/(n-k)$,
$|S|=n-k$.
\end{lemma}

The following standard lemma is proved in Appendix~\ref{app:calM:d:convex:hull} for completeness.
\begin{lemma}[Rank-one decomposition]\label{lemma:calM:d:convex:hull}
$\calM_d=\co\{\bu^{\otimes2}:\bu\in\mbS^{d-1}\}$. Consequently,
$\lambda_1(\bfA)=\max_{\bfM\in\calM_d}\langle\bfA,\bfM\rangle$ and $\lambda_d(\bfA)=\min_{\bfM\in\calM_d}\langle\bfA,\bfM\rangle$ for every
$\bfA\in\calS^d$.
\end{lemma}

\subsection{Notation for random samples and estimator-related quantities}\label{sub:p:sample}

For convenience, we collect here the notation related to Assumption~\ref{assump:finite:covariance:contamination} and to our weighted-average estimators. Let $\bx_1,\ldots,\bx_n$ be i.i.d.\ copies of a random vector $\bx\in\re^d$ with mean $\bmmu:=\esp[\bx]$ and covariance matrix $\bfSigma\neq\bzero$. The contamination level is a parameter $\varepsilon\in[0,1/2)$. The observed sample $\tilde\bx_1,\ldots,\tilde\bx_n\in\re^d$ is arbitrary except that $\#\{i\in[n]:\tilde\bx_i\neq\bx_i\}\le\varepsilon n$. We write $\bxi:=\bx-\bmmu$, $\bxi_i:=\bx_i-\bmmu$, and $\tilde\bxi_i:=\tilde\bx_i-\bmmu$ for the vectors centered at $\bmmu$.

Our estimators of $\bmmu$ are weighted averages of the observed sample, $\bmmu_w:=\sum_{i=1}^n w_i\tilde\bx_i$ for $w\in\Delta_{n,k}$. Let $\mbX\in\re^{n\times d}$ be the data matrix with rows $\tilde\bx_1^\top,\ldots,\tilde\bx_n^\top$, so that $\bmmu_w=\mbX^\top w$, and let $M_{\mbX}:=\max_{i\in[n]}\|\tilde\bx_i\|_2$. We also use the weighted covariance-type matrices $\bfSigma_w(\bnu):=\sum_{i=1}^n w_i(\tilde\bx_i-\bnu)(\tilde\bx_i-\bnu)^\top$ and $\bfSigma_w:=\bfSigma_w(\bmmu_w)$.

\subsection{Bregman geometry}\label{s:bregman:geometry}

Let $h:\re^n\to\re\cup\{+\infty\}$ be a proper, closed, convex function, differentiable on $\interior(\dom(h))\neq\emptyset$, where
$\dom(h):=\{x:h(x)<+\infty\}$. Following~\cite{2018Lu:Freund:Nesterov}, we call
$h$ a \emph{reference function}; its \emph{Bregman
divergence}~\cite{1967Bregman} is
$
    D_h(x\|x') := h(x)-h(x')-\langle\nabla h(x'),x-x'\rangle,
$
for $x\in\dom(h)$ and $x'\in\interior(\dom(h))$. Let  $\Delta\subseteq\re^n$ be a nonempty convex set with $\Delta\cap\interior(\dom(h))\neq\emptyset$. If, for some norm $\|\cdot\|$, $D_h(x\|x')\ge\frac12\|x-x'\|^2$ for all $x\in\Delta\cap\dom(h)$ and  $x'\in\Delta\cap\interior(\dom(h))$, we call $h$ a \emph{distance-generating function} on $\Delta$ with respect to $\|\cdot\|$~\cite{1983Nemirovskij:Yudin,2012Lan}.

We use three reference functions on $\Delta_n$. Recall that $\bmmu_w=\mbX^\top w$.

\paragraph{Entropy.}
The \emph{negative entropy}
$h_{\KL}(w):=\sum_{i=1}^nw_i\log w_i$ for $w\in\re^n_+$ (with $0\log0:=0$), and $h_{\KL}(w):=+\infty$ otherwise, has  $\Delta_n\cap\interior(\dom(h_{\KL}))=\ri(\Delta_n)$. Its Bregman divergence is the Kullback--Leibler divergence,
\begin{align}
    D_{\KL}(\omega\|w)=\sum_{i=1}^n\omega_i\log\frac{\omega_i}{w_i},
    \qquad \omega\in\Delta_n,\ w\in\ri(\Delta_n).
\end{align}
By Pinsker's inequality, $h_{\KL}$ is a distance-generating function on $\Delta_n$ with respect to $\|\cdot\|_1$.

\paragraph{Euclidean and mean-induced divergences.}
The reference functions $h_{\E}(w):=\frac12\|w\|_2^2$ and $h_{\mbX}(w):=\frac12\|\bmmu_w\|_2^2=\frac12\|\mbX^\top w\|_2^2$ are finite on $\re^n$, so every point of $\Delta_n$, including its boundary, is admissible. Their Bregman divergences are 
\begin{align}
    D_{\E}(\omega\|w)=\tfrac12\|\omega-w\|_2^2,
    \qquad
    D_{\mbX}(\omega\|w)=\tfrac12\|\bmmu_\omega-\bmmu_w\|_2^2=\tfrac12\|\mbX^\top(\omega-w)\|_2^2.
\end{align}

We call $D_{\mbX}$ the \emph{mean-induced divergence}. For $\alpha\ge0$, the \emph{regularized mean-induced divergence}
$
    D_{\mbX,\alpha}(\omega\|w)
    :=D_{\mbX}(\omega\|w)+\alpha D_{\E}(\omega\|w)
    =\tfrac12\|\omega-w\|_{\mbX,\alpha}^2
$
is induced by $h_{\mbX,\alpha}:=h_{\mbX}+\alpha h_{\E}$, where
$\|u\|_{\mbX,\alpha}^2:=\|\mbX^\top u\|_2^2+\alpha\|u\|_2^2
=\langle u,(\mbX\mbX^\top+\alpha\bfI_n)u\rangle$. For $\alpha>0$, $\|\cdot\|_{\mbX,\alpha}$ is a norm, so $h_{\mbX,\alpha}$ is a distance-generating function on $\Delta_n$ with respect to it.

\subsection{Relative regularity}\label{ss:relative:regularity}

Throughout, $h$ and $\Delta$ are as in Section~\ref{s:bregman:geometry}. In this subsection only, $F$ and $\hat F$ denote general functions, which may differ from those of the Introduction. Let $F:\re^n\to\re\cup\{+\infty\}$ be a proper closed function to be minimized over $\Delta$. For $w\in\dom(F)$, the \emph{Fr\'echet subdifferential} $\partial F(w)$ is the set of $g\in\re^n$ such that $F(w')\ge F(w)+\langle g,w'-w\rangle+o(\|w'-w\|)$ as $w'\to w$; for $w\notin\dom(F)$, we set $\partial F(w):=\emptyset$~\cite{1998Rockafellar:Wets}.

\emph{Relative} notions measure smoothness, Lipschitz continuity or convexity in the geometry of a reference function $h$ instead of the Euclidean norm, and are useful when $F$ is not Euclidean-Lipschitz or smooth~\cite{2017Bauschke:Bolte:Teboulle,2018Lu:Freund:Nesterov,2019Lu}. In robust mean estimation, $F$ is not directly accessible: we only query subgradients of a surrogate $\hat F$, so we state relative notions for $\hat F$ as a \emph{subgradient oracle} for $F$. In what follows, $\hat F:\re^n\to\re\cup\{+\infty\}$ is proper and closed with $\dom(\hat F)=\dom(F)$, $\Delta^\circ:=\Delta\cap\interior(\dom(h))\cap\dom(F)$, $\Delta':=\Delta\cap\dom(h)\cap\dom(F)$, and $w^*\in\Delta'$ minimizes $F$ over $\Delta$.

\begin{definition}[Relative inexact subgradient oracle]\label{def:inexact:oracle}
We say that $\hat F$ is, relative to $h$ on $\Delta$,
\begin{enumerate}
    \item[\rm (a)] an \emph{$M$-Lipschitz continuous inexact oracle}, with $M\ge0$, if, for all $w\in\Delta^\circ$, $g\in\partial\hat F(w)$ and $w'\in\Delta'$,
    $$
      \langle g,w-w'\rangle\le M\sqrt{2D_h(w'\|w)}, 
    $$
    \item[\rm (b)] a \emph{$(\rho,\zeta,\tau)$-weakly star-convex inexact oracle for $F$ at $w^*$}, with $\rho,\zeta,\tau\ge0$, if, for all $w\in\Delta^\circ$ and $g\in\partial\hat F(w)$,
    \[
      F(w^*)\ge F(w) + \langle g,w^*-w\rangle -\frac\rho2\Big(\sqrt{2D_h(w^*\|w)}+\zeta\Big)^2-\tau.
    \]
\end{enumerate}
\end{definition}
When $\hat F=F$: item~(a) is a directional form of the relative continuity
of~\cite{2019Lu}; item~(b) with $\rho=\tau=0$ is star convexity at
$w^*$~\cite{2020Hinder:Sidford:Sohoni}; and if $\zeta=\tau=0$,
$h=\frac12\|\cdot\|_2^2$, and (b) holds with every point of $\Delta$ in place of the minimizer $w^*$, then (b) is $\rho$-weak convexity, i.e., convexity of $F+\frac\rho2\|\cdot\|_2^2$~\cite{2019Davis:Drusvyatskiy}. For $\zeta,\tau>0$, item~(b) allows an additional slack $\rho\,\zeta\sqrt{2D_h(w^*\|w)}+\frac\rho2\zeta^2+\tau$. 

Finally, we quantify how fast $F$ grows away from $w^*$, up to a radius at which growth may fail.

\begin{definition}[Relative inexact quadratic growth]\label{def:relative:iqg}
Let $m>0$ and $r\ge0$. We say that $F$ has \emph{inexact quadratic growth relative to $h$} on $\Delta$, abbreviated $\IQG_h(m,r)$, with growth constant $m$ and inexactness radius $r$, if
\begin{align}
    F(w) \ge F(w^*) + m\Big(\sqrt{2D_h(w^*\|w)}-r\Big)_+^2
    \qquad\text{for all } w\in\Delta^\circ.
\end{align}
\end{definition}
When $r=0$ and $h=\frac{1}{2}\|\cdot\|_2^2$, this condition is related to the quadratic growth property of \cite{2019Necoara:Nesterov:Glineur}. For $r>0$, growth is required only where $\sqrt{2D_h(w^*\|w)}>r$; thus $r$ is the statistical precision below which the landscape of $F$ carries no information.

\section{Robust averages and counting}\label{s:robust:averages:counting}

This deterministic section relates weighted averages over $\Delta_{n,k}$ to trimmed and truncated empirical means; it is used in Sections~\ref{s:counting:truncation:in:probability} and~\ref{s:stability:spectral:signature}. Throughout, $\bz_{1:n}:=(\bz_1,\ldots,\bz_n)$ are points in $\re^d$, $w\in\Delta_n$, $k\in[n-1]$, $B>0$, and 
\begin{align}
  {\sf clip}_B(t):=\max\{\min\{t,B\},-B\},\quad(t\in\re),
  \label{eq:clip:B}
\end{align}
is the truncation function; see also \eqref{eq:rem:B} in Section \ref{s:counting:truncation:in:probability}. For $\bv\in\mbS^{d-1}$ and $\bfM\in\calM_d$, the \emph{weighted} and \emph{truncated empirical means} of the linear values $\langle\bz_i,\bv\rangle$ and the quadratic values $\langle\bz_i^{\otimes2},\bfM\rangle$ are
\begin{align}
  \sfT_1(\bz_{1:n}\,|\,w,\bv)&:=\sum_{i=1}^nw_i\langle\bz_i,\bv\rangle,&
  \sfT_{1,B}(\bz_{1:n}\,|\,\bv)&:=\frac1n\sum_{i=1}^n{\sf clip}_B\big(\langle\bz_i,\bv\rangle\big),\\
  \sfT_2(\bz_{1:n}\,|\,w,\bfM)&:=\sum_{i=1}^nw_i\langle\bz_i^{\otimes2},\bfM\rangle,&
  \sfT_{2,B}(\bz_{1:n}\,|\,\bfM)&:=\frac1n\sum_{i=1}^n
  {\sf clip}_{B}^2(\|\bfC_{\bfM}\bz_i\|_2),
\end{align}
 and the \emph{trimmed empirical mean} of the quadratic values is
\begin{align}
  \sfT_{2,k}(\bz_{1:n}\,|\,\bfM)
  :=\min_{S\subseteq[n],\,|S|=n-k}\frac1{n-k}\sum_{i\in S}\langle\bz_i^{\otimes2},\bfM\rangle.
\end{align}

The following direct consequence of Lemma~\ref{lemma:capped:extreme} shows that minimizing a weighted average over the capped simplex amounts to trimming the $k$ largest values.
\begin{lemma}[Minimal weighted average = upper-trimmed average]
\label{lemma:min:weight:trimming}
For every $\bfM\in\calM_d$,
$
    \sfT_{2,k}(\bz_{1:n}\,|\,\bfM)
    = \min_{w\in\Delta_{n,k}}\sum_{i=1}^nw_i\langle\bz_i^{\otimes2},\bfM\rangle,
$
and the minimum is attained at the weights $\frac{\unit_S}{n-k}$ with any set $S\subseteq[n]$ of $n-k$ indices carrying the $n-k$ smallest values of $\langle\bz_i^{\otimes2},\bfM\rangle$.
\end{lemma}

We use the counting condition from the analysis of trimmed means~\cite{2019Oliveira:Orenstein,2021Lugosi:Mendelson,2024Oliveira:Rico,2026Abdalla:Zhivotovskiy,2025Oliveira:Resende}, formulated with the norms $\|\bfC_{\bfM}\bh\|_2=\langle\bh^{\otimes2},\bfM\rangle^{1/2}$, $\bfM\in\calM_d$, which simplifies the analysis (see Section~\ref{ss:counting:truncation:preliminaries}).
\begin{definition}[Counting condition]\label{def:counting}
$\bz_{1:n}$ satisfies the \emph{counting condition}
$\Count_2(\bz_{1:n},B,k)$ if $\#\{i\in[n]:\|\bfC_{\bfM}\bz_i\|_2>B\}\le k$ for all $\bfM\in\calM_d$.
\end{definition}

Under this condition, trimmed and truncated averages of the quadratic values are uniformly close; see~\cite[Claim~1]{2024Oliveira:Rico} for a related result. See also \cite{2026Oliveira:Rico:Thompson}. The proof is in Appendix~\ref{app:trim:trunc}.
\begin{lemma}[Counting $\Rightarrow$ trimming $\approx$ truncation]\label{lemma:trim:trunc}
If $\Count_2(\bz_{1:n},B,k)$ holds, then, for every $\bfM\in\calM_d$,
$
  \sfT_{2,B}(\bz_{1:n}\,|\,\bfM)-\frac{B^2k}{n}
  \le\sfT_{2,k}(\bz_{1:n}\,|\,\bfM)
  \le\sfT_{2,B}(\bz_{1:n}\,|\,\bfM).
$
\end{lemma}

\section{Counting and truncation in high probability}\label{s:counting:truncation:in:probability}

This section shows that, with high probability, the clean sample satisfies the counting condition (Lemma~\ref{lemma:counting}) and its truncated empirical processes concentrate (Lemma~\ref{lemma:truncated:processes}); Corollary~\ref{cor:favorable:event} combines both. Throughout, Assumption~\ref{assump:finite:covariance:contamination} holds and we use the notation in Section \ref{sub:p:sample}, in particular, $\bxi:=\bx-\bmmu$, $\bxi_i:=\bx_i-\bmmu$, $\tilde\bxi_i:=\tilde\bx_i-\bmmu$, $B>0$ is fixed, we define the function
\begin{align}
  {\sf rem}_B(t) := t - {\sf clip}_B(t), \quad (t\in\re),
  \label{eq:rem:B}
\end{align}
where ${\sf clip}_B$ is defined in \eqref{eq:clip:B}. $\Gamma$ is the standard Gaussian measure on $\re^d$, and $\Gamma f(\bg,y):=\int_{\re^d}f(\bg,y)\,\Gamma(d\bg)$ when well defined.

\subsection{Preliminaries}\label{ss:counting:truncation:preliminaries}

The following lemma is proved in Appendix~\ref{app:CM:moment}.
\begin{lemma}\label{lemma:CM:moment}
  For all $\bfM\in\calM_d$, $\esp\|\bfC_{\bfM}\bxi\| _2^p\le \nu_p^p$.
\end{lemma}

Proposition~\ref{prop:multvect}, proved in Appendix~\ref{app:multvect}, writes $\|\bfC_{\bfM}\bh\|_2$ as an average of one-dimensional projections. 
It lets us move back and forth expressions involving density matrices and vectors. As such, it is a more efficient alternative to the randomized rounding arguments of~\cite{2019Cheng:Diakonikolas:Ge,2020Hopkins:Li:Zhang,2020Lei:Luh:Venkat:Zhang,2022Depersin:Lecue}.

\begin{proposition}\label{prop:multvect}
For any $\bfM\in\calM_d$, there exists a probability measure
$\bar{\Gamma}_{\bfC_{\bfM}}$ supported on $\mbS^{d-1}$, and depending only on
$\bfC_{\bfM}$, such that, for all $\bh\in\re^d$,
\[
  \|\bfC_{\bfM}\bh\|_2 = \sqrt{\frac{\pi}{2}}\,\big(\Gamma\|\bfC_{\bfM}^\top\bg\|_2\big)
  \int_{\mbS^{d-1}}|\langle\bar\bg,\bh\rangle|\,\bar{\Gamma}_{\bfC_{\bfM}}(d\bar\bg),
\]
where $0<\Gamma\|\bfC_{\bfM}^\top\bg\|_2\le1$.
\end{proposition}

We also need to control the truncation bias; the proof is in
Appendix~\ref{app:bias}.
\begin{lemma}[Truncation bias]\label{lemma:bias}
For all $\bv\in\mbS^{d-1}$ and $p\ge1$,
$
  \big|\esp[{\sf clip}_B(\langle\bv,\bxi\rangle)]\big|
  =\big|\esp[{\sf rem}_B(\langle\bv,\bxi\rangle)]\big|
  \le\frac{\nu_p^p}{B^{p-1}}.
$
Moreover, for all $\bfM\in\calM_d$ and $p\ge2$,
$
  0\le\langle\bfSigma,\bfM\rangle-\esp\big[{\sf clip}_B^2(\|\bfC_{\bfM}\bxi\|_2)\big]
  \le\frac{\nu_p^p}{B^{p-2}}.
$
\end{lemma}

\subsection{Truncated empirical and Rademacher complexities}

We bound the empirical and Rademacher complexities of the truncated linear and quadratic processes, used in Sections~\ref{ss:counting:arguments} and~\ref{ss:truncation:arguments}. Let
\begin{align}
  \calC_n:=\esp\left\|\frac1n\sum_{i=1}^n\bxi_i\right\|_2,
\end{align}
and let $\epsilon_1,\ldots,\epsilon_n$ be i.i.d.\ Rademacher variables independent of $(\bxi_i)_{i\in[n]}$.

\begin{proposition}[Truncated empirical and Rademacher complexities]\label{prop:rademacher}
Define
\begin{align}
  \calE_{1,B}(\mbS^{d-1})&:=\esp\sup_{\bv\in\mbS^{d-1}}\Big|\frac1n\sum_{i=1}^n
  \Big({\sf clip}_B(\langle\bv,\bxi_i\rangle)-\esp[{\sf clip}_B(\langle\bv,\bxi\rangle)]\Big)\Big|,\\
  \calE_{2,B}(\calM_d)&:=\esp\sup_{\bfM\in\calM_d}\Big|\frac1n\sum_{i=1}^n
  \Big({\sf clip}_B^2(\|\bfC_{\bfM}\bxi_i\|_2)-\esp[{\sf clip}_B^2(\|\bfC_{\bfM}\bxi\|_2)]\Big)\Big|,\\
  \calR_{q,B}(\calM_d)&:=\esp\sup_{\bfM\in\calM_d}\Big|\frac1n\sum_{i=1}^n
  \epsilon_i\,{\sf clip}_B^q(\|\bfC_{\bfM}\bxi_i\|_2)\Big|,\qquad q\in\{1,2\}.
\end{align}
Then $\calE_{2,B}(\calM_d)\le2\calR_{2,B}(\calM_d)$,
$\calR_{2,B}(\calM_d)\le4B\,\calR_{1,B}(\calM_d)$,
$\calR_{1,B}(\calM_d)\le4\sqrt{2\pi}\,\calC_n$, and
$\calE_{1,B}(\mbS^{d-1})\le4\calC_n$.
\end{proposition}
\begin{proof}
The inequality $\calE_{2,B}(\calM_d)\le2\calR_{2,B}(\calM_d)$ is the classical
symmetrization inequality; see~\cite[Lemma~11.4]{2013Boucheron:Lugosi:Massart}. Now define the function $h_{B}(t):={\sf clip}_B^2(t)$ for $t\in\re$. This  function is $2B$-Lipschitz, satisfies $h_{B}(0)=0$, and moreover:
$
{\sf clip}_B^2(\|\bfC_{\bfM}\bxi_i\|_2) = h_{B}({\sf clip}_B(\|\bfC_{\bfM}\bxi_i\|_2)).
$

As a result, the Ledoux-Talagrand contraction theorem for Rademacher complexities \cite{1991Ledoux:Talagrand} (see \cite[Theorem 11.6]{2013Boucheron:Lugosi:Massart}) implies:
$$
  \esp\left[\sup_{\bfM\in\calM_d}\left|\frac{1}{n}\sum_{i=1}^n\epsilon_ih_{B}({\sf clip}_B(\|\bfC_{\bfM}\bxi_i\|_2))\right|\right]\leq 4B\esp\left[\sup_{\bfM\in\calM_d}\left|\frac{1}{n}\sum_{i=1}^n\epsilon_i{\sf clip}_B(\|\bfC_{\bfM}\bxi_i\|_2)\right|\right],
$$
that is, $\calR_{2,B}(\calM_d)\leq 4B\calR_{1,B}(\calM_d)$.

To bound $\calR_{1,B}(\calM_d)$, we apply the contraction principle again to remove the function ${\sf clip}_B$ (which is $1$-Lipschitz):
$
\calR_{1,B}(\calM_d)\leq 2\esp\left[\sup_{\bfM\in\calM_d}\left|\frac{1}{n}\sum_{i=1}^n\epsilon_i\|\bfC_{\bfM}\bxi_i\|_2\right|\right].
$

Now recall from Proposition \ref{prop:multvect} that, for given $\bfM\in\calM_d$, 
\[
  \|\bfC_{\bfM}\bxi_i\|_2 = \sqrt{\frac{\pi}{2}}\,\big(\Gamma\|\bfC_{\bfM}^\top\bg\|_2\big)
  \int_{\mbS^{d-1}}|\langle\bar\bg,\bxi_i\rangle|\,\bar{\Gamma}_{\bfC_{\bfM}}(d\bar\bg). 
\]
Therefore, 
\begin{align}
    \left|\frac{1}{n}\sum_{i=1}^n\epsilon_i\|\bfC_{\bfM}\bxi_i\|_2\right| &= \sqrt{\frac{\pi}{2}}\,\Gamma\|\bfC_{\bfM}^\top\bg\|_2\,\left|\int_{\mbS^{d-1}}\frac{1}{n}\sum_{i=1}^n\epsilon_i|\langle\bar{\bg},\bxi_i\rangle|\bar{\Gamma}_{\bfC_{\bfM}}(d\bar{\bg})\right|\\ 
    &\le\sqrt{\frac{\pi}{2}}\,\Gamma\|\bfC_{\bfM}^\top\bg\|_2\,\int_{\mbS^{d-1}}\left|\frac{1}{n}\sum_{i=1}^n\epsilon_i|\langle\bar{\bg},\bxi_i\rangle|\right|\bar{\Gamma}_{\bfC_{\bfM}}(d\bar{\bg})\\ 
    &\leq \sqrt{\frac{\pi}{2}}\sup_{\bv\in\mbS^{d-1}}\left|\frac{1}{n}\sum_{i=1}^n\epsilon_i|\langle\bv,\bxi_i\rangle|\right|.
\end{align}
In the first inequality above we used Jensen's inequality, and in the second, we used $\Gamma\|\bfC_{\bfM}^\top\bg\|_2\le1$ and bounded the integrand by the sup and used that $\bar\Gamma_{\bfC_{\bfM}}$ is a probability measure. 

Since $|\cdot|$ is $1$-Lipschitz and $|0|=0$, we may apply Ledoux-Talagrand contraction (again see \cite[Theorem 11.6]{2013Boucheron:Lugosi:Massart}), followed by desymmetrization \cite[Lemma 11.4]{2013Boucheron:Lugosi:Massart} (with $\esp[\bxi_i]=0$) to obtain:
\begin{align}
  \esp\left[\sup_{\bv\in\mbS^{d-1}}\left|\frac{1}{n}\sum_{i=1}^n\epsilon_i|\langle\bv,\bxi_i\rangle|\right|\right] 
  &\leq 2\esp\left[\sup_{\bv\in\mbS^{d-1}}\left|\frac{1}{n}\sum_{i=1}^n\epsilon_i\,\langle\bv,\bxi_i\rangle\right|\right]\\ 
  &\le 2\esp\left[\sup_{\bv\in\mbS^{d-1}}\left|\left\langle \bv,\frac{1}{n}\sum_{i=1}^n\epsilon_i\,\bxi_i\right\rangle\right|\right] \\ 
  & \leq 4\esp\left[\sup_{\bv\in\mbS^{d-1}}\left|\left\langle \bv,\frac{1}{n}\sum_{i=1}^n\bxi_i\right\rangle\right|\right]\\ 
  & = 4\esp\left[\left\|\frac{1}{n}\sum_{i=1}^n\bxi_i\right\|_2\right] = 4\calC_n.
\end{align}
Combining the preceding displays shows $\calR_{1,B}(\calM_d)\leq 4\sqrt{2\pi}\calC_n$. The inequality $\calE_{1,B}(\mbS^{d-1})\leq 4\calC_n$ follows similarly via symmetrization, contraction and desymmetrization, but using the fact that $\sup_{\bv\in\mbS^{d-1}}|\sum_i\epsilon_i\mathsf{clip}_B(\langle\bv,\bxi_i\rangle)|=\sup_{\bv\in\mbS^{d-1}}\sum_i\epsilon_i\mathsf{clip}_B(\langle\bv,\bxi_i\rangle)$ since $\mathsf{clip}_B$ is odd; we omit the details.
\end{proof}

\subsection{Counting arguments}\label{ss:counting:arguments}

We use the following result of~\cite{2025Oliveira:Resende} (see also~\cite{2021Lugosi:Mendelson}), restated in our notation.
    \begin{lemma}[Counting Lemma; Lemma 8.1 in \cite{2025Oliveira:Resende}] 
    \label{lemma:counting:lemma}
    Suppose $Z_1,\dots,Z_n$ is an i.i.d. sample of random variables taking values in some measurable space $(\bfZ,\calZ)$. Moreover, let $\calF$ be countable family of measurable functions $f:\bfZ\to\re$ such that $\esp[|f(Z_1)|]<+\infty$. For $B>0$, define
    $
    {\rm Rad}_B(\calF):=\esp\left[\sup_{f\in\calF}\left|\frac{1}{n}\sum_{i=1}^n\epsilon_i{\sf clip}_B(f(Z_i))\right|\right].
    $
    Then, if $s\in\mathbb{N}$ and $B>0$ satisfy
    \begin{align}
      \frac{8}{B}\,{\rm Rad}_B(\calF) + \sup_{f\in\calF}\prob\{|f(Z_1)|>B/2\}\leq \frac{s}{8n},\label{eq:conditionforcounting}
    \end{align}
    then 
    $
      \prob\{\forall f\in\calF\,:\, \#\{i\in[n]\,:\, |f(Z_i)|>B\}\leq s\}\geq 1-e^{-s}.
    $
    \end{lemma}

\begin{lemma}[Counting condition in high probability]\label{lemma:counting}
      If
      \begin{align}
        B \ge \max\left\{ 2\nu_p\left(\frac{16\,n}{s}\right)^{1/p}, 512\sqrt{2\pi}\calC_n\left(\frac{n}{s}\right) \right\}, \label{lemma:counting:threshold:condition}
      \end{align}
  then, with probability at least $1-e^{-s}$, $\Count_{2}(\bxi_{1:n},B,s)$ holds, that is, 
  $
    \forall\bfM\in\calM_d,~\#\{i\in[n] : \|\bfC_{\bfM}(\bxi_i)\|_2>B\}\leq s. 
  $
\end{lemma}
\begin{proof}
  We apply Lemma~\ref{lemma:counting:lemma} with $(\bfZ,\calZ):=(\re^d,\calB(\re^d))$, $Z_i:=\bxi_i$ and $\calF:=\{f_{\bfM}:=\|\bfC_{\bfM}(\cdot)\|_2:\bfM\in\calM_d\}$. Since $\bfM\mapsto f_{\bfM}(\bv)$ is continuous for every fixed vector $\bv\in\re^d$, suprema over $\calF$ coincide with suprema over a countable dense subfamily, to which the lemma applies. Moreover, ${\rm Rad}_B(\calF)=\calR_{1,B}(\calM_d)\le4\sqrt{2\pi}\,\calC_n$ by Proposition~\ref{prop:rademacher}.

  By Markov's inequality and Lemma~\ref{lemma:CM:moment}, for $p\ge2$,
  \begin{align}
    \sup_{f\in\calF}\prob\{|f(Z_1)|>B/2\}
    &= \sup_{\bfM\in\calM_d}\prob\{\|\bfC_{\bfM}\bxi_1\|_2>B/2\} \\
    &\le \sup_{\bfM\in\calM_d}\frac{\esp\|\bfC_{\bfM}\bxi_1\|_2^p}{(B/2)^p}
    \le\Big(\frac{2\nu_p}{B}\Big)^p.
  \end{align}

  Hence~\eqref{eq:conditionforcounting} holds whenever $\frac{32\sqrt{2\pi}\calC_n}{B}+\big(\frac{2\nu_p}{B}\big)^p\le\frac{s}{8n}$, which follows from~\eqref{lemma:counting:threshold:condition}, since it gives $\frac{32\sqrt{2\pi}\calC_n}{B}\le\frac{s}{16n}$ and $\big(\frac{2\nu_p}{B}\big)^p\le\frac{s}{16n}$. Lemma~\ref{lemma:counting:lemma} concludes the proof.
\end{proof}

\subsection{Truncation arguments}\label{ss:truncation:arguments}

We use Bousquet's version of Talagrand's inequality~\cite{2002Bousquet}; see also~\cite[Thm.~12.5]{2013Boucheron:Lugosi:Massart}.

\begin{theorem}[Bousquet's inequality~\cite{2002Bousquet}]\label{thm:bousquet}
Let $Z_1,\ldots,Z_n$ be i.i.d.\ random elements of a measurable space
$(\bfZ,\calZ)$, and let $\calF$ be a countable family of measurable functions
$f:\bfZ\to\re$ with $\esp|f(Z_1)|<\infty$ and $|f(Z_1)-\esp f(Z_1)|\le C$ almost surely, for some constant $C>0$. Let
$
  S:=\sup_{f\in\calF}|\frac1n\sum_{i=1}^nf(Z_i)-\esp f(Z_1)|
$
and
$
  \sigma^2(\calF):=\sup_{f\in\calF}{\rm Var}\big(f(Z_1)\big).
$
Then, for every $t>0$,
\[
  \prob\left\{S-\esp S\le\sqrt{\frac{2t}{n}\big(\sigma^2(\calF)+2C\,\esp S\big)}
  +\frac{Ct}{3n}\right\}\ge1-e^{-t}.
\]
\end{theorem}

In particular, for every $t>0$,\footnote{This follows from the elementary facts
$\sqrt{a+b}\le\sqrt a+\sqrt b$ and $2\sqrt{xy}\le x+y$ for $a,b,x,y\ge0$.}
\begin{align}
  \prob\left\{
    S-\esp[S]\le \frac{4}{3}\left(\esp S+ \sigma(\calF)\sqrt{\frac{2t}{n}}
  + \frac{Ct}{n}\right)
\right\}\ge1-e^{-t}.
  \label{eq:bousquet:simplified}
\end{align}

\begin{lemma}[Truncated processes with high probability]
\label{lemma:truncated:processes}
Let $\delta\in(0,1)$ and, for $b\in\{1,2\}$,
$S_b:=\sup_{f\in\calF_b}|\frac1n\sum_{i=1}^nf(\bxi_i)-\esp f(\bxi)|$, where
\[
  \calF_1:=\left\{{\sf clip}_B(\langle\bv,\cdot\rangle):\bv\in\mbS^{d-1}\right\}
  ~~\text{ and }~~
  \calF_2:=\left\{{\sf clip}_B^2(\|\bfC_{\bfM}(\cdot)\|_2):\bfM\in\calM_d\right\}.
\]
Then there exist events ${\rm E}_1(\delta)$ and ${\rm E}_2(\delta)$, each of probability at least $1-\delta$, such that
\begin{align}
  S_1&\le\frac{28}{3}\,\calC_n+\frac{4\sqrt2}{3}\,\nu_2\sqrt{\frac{\log(1/\delta)}{n}}
  +\frac{8B\log(1/\delta)}{3n}\quad\text{on }{\rm E}_1(\delta),\\
    S_2&\le B\left(\frac{224\sqrt{2\pi}}{3}\,\calC_n
  +\frac{4\sqrt2}{3}\,\nu_2\sqrt{\frac{\log(1/\delta)}{n}}
  +\frac{4B\log(1/\delta)}{3n}\right)\quad\text{on }{\rm E}_2(\delta).
\end{align}
\end{lemma}
\begin{proof}
We apply~\eqref{eq:bousquet:simplified} with $Z_i:=\bxi_i$,
$(\bfZ,\calZ):=(\re^d,\calB(\re^d))$ and $t:=\log(1/\delta)$, to $\calF_1$ and to
$\calF_2$. Both families are uncountable, but they depend continuously on
$\bv\in\mbS^{d-1}$ and $\bfM\in\calM_d$, so $S_1$ and $S_2$ coincide with suprema
over countable dense subfamilies, to which Theorem~\ref{thm:bousquet} applies.

The functions in $\calF_1$ take values in $[-B,B]$, so $|f-\esp f|\le2B$; functions in $\calF_2$ take values in $[0,B^2]$, so $|f-\esp f|\le B^2$.

Proposition~\ref{prop:rademacher} gives $\esp S_1\le4\calC_n$ and $\esp S_2\le2\calR_{2,B}(\calM_d)\le8B\,\calR_{1,B}(\calM_d)\le32\sqrt{2\pi}\,B\,\calC_n$.

Since $|{\sf clip}_B(t)|=\min\{B,|t|\}$,
\[
  \sigma^2(\calF_1)\le\sup_{\bv\in\mbS^{d-1}}\esp\big[{\sf clip}_B^2(\langle\bv,\bxi\rangle)\big]
  \le\sup_{\bv\in\mbS^{d-1}}\esp\langle\bv,\bxi\rangle^2=\nu_2^2,
\]
and, by Lemma~\ref{lemma:CM:moment} with $p=2$,
\[
  \sigma^2(\calF_2)\le\sup_{\bfM\in\calM_d}\esp\left[{\sf clip}_B^4(\|\bfC_{\bfM}\bxi\|_2)\right]
  \le B^2\sup_{\bfM\in\calM_d}\esp\|\bfC_{\bfM}\bxi\|_2^2\le B^2\nu_2^2.
\]
Substituting these bounds into~\eqref{eq:bousquet:simplified} gives the two
claims.
\end{proof}

\section{Spectral signature, favorable sample and stability}\label{s:stability:spectral:signature}

Throughout this section, Assumption~\ref{assump:finite:covariance:contamination} holds. We give unified proofs, valid for all $p\ge2$, of \emph{stability} and of the \emph{spectral signature} property, organized around the notion of a \emph{favorable sample}, which underlies our estimation theory. The proofs combine the deterministic results of Section~\ref{s:robust:averages:counting} with the high-probability bounds of Section~\ref{s:counting:truncation:in:probability}.

We first recall, in our notation, the \emph{distributional} definition of stability of~\cite[Def.~2.1]{2023Diakonikolas:Kane}.\footnote{In~\cite[Def.~2.1]{2023Diakonikolas:Kane}, $\bfSigma=\bfI_d$ for convenience; see the comments therein.} Using the notation in Section \ref{sub:p:sample}:
\begin{definition}[Stability with respect to a distribution in $\Delta_n$]
\label{def:stab:mass}
Let $k\in[n]$ with $k<n/2$, $R>0$ and $\omega\in\Delta_n$. The sample $\tilde\bx_{1:n}$ is \emph{$(R,\varepsilon_{n,k})$-stable with respect to $\omega$} if, for every $w\in\Delta_{n,k}(\omega)$,
\begin{itemize}
  \item[\rm (a)] $\|\bmmu_w-\bmmu\|_2\le R$,
  \item[\rm (b)] $\|\bfSigma_w(\bmmu)-\bfSigma\|\le\frac{R^2}{\varepsilon_{n,k}}$.
\end{itemize}
\end{definition}
The \emph{subsample} version of this definition~\cite[Def.~2.1]{2023Diakonikolas:Kane}, used in~\cite{2016Diakonikolas:Kamath:Kane:Li:Moitra:Stewart,2020Diakonikolas:Kane:Pensia}, instead removes $k$ \emph{points} from the sample. Since $\Delta_{n,k}(\unit_{[n]}/n)=\Delta_{n,k}$, and $n\varepsilon_{n,k}=k\in\mathbb{N}$, the two versions coincide~\cite[Lemma~2.2]{2023Diakonikolas:Kane}. 
 
For the algorithmic analysis of our estimators, a different notion suffices, which we call a \emph{favorable sample}.

\begin{definition}[Favorable sample]\label{def:favorable}
Let $k\in[n]$ with $k<n$, and $r_1>0$. The sample $\tilde\bx_{1:n}$ is
\emph{$(r_1,\varepsilon_{n,k})$-favorable} if the following hold:
\begin{enumerate}
  \item[\rm (a)] for all $w\in\Delta_{n,k}$,
  $
    \|\bmmu_w-\bmmu\|_2\le r_1+\sqrt{\frac{\varepsilon_{n,k}}{1-\varepsilon_{n,k}}
    \Big(\lambda_1(\bfSigma_w(\bmmu)-\bfSigma) 
    +\frac{r_1^2}{\varepsilon_{n,k}}\Big)},
  $
  \item[\rm (b)] there exist \emph{good weights} $w^*\in\Delta_{n,k}$ with $\lambda_1\big(\bfSigma_{w^*}(\bmmu)-\bfSigma\big)\le r_1^2/\varepsilon_{n,k}$,
  \item[\rm (c)] $\min_{w\in\Delta_{n,k}}\lambda_d\big(\bfSigma_w(\bmmu)-\bfSigma\big)
  \ge-r_1^2/\varepsilon_{n,k}$.
\end{enumerate}
\end{definition}

Following~\cite{2020Hopkins:Li:Zhang}, we call item~(a) the \emph{spectral signature} property: if the weighted second moment does not exceed $\bfSigma$ by much, the weighted mean is close to $\bmmu$. Items~(b) and~(c) give $\|\bfSigma_{w^*}(\bmmu)-\bfSigma\|\le r_1^2/\varepsilon_{n,k}$, as for the good weights of~\cite[Def.~1.1]{2020Hopkins:Li:Zhang}, but relative to $\bfSigma$. Moreover, item~(c) holds \emph{uniformly} over $\Delta_{n,k}$; to our knowledge, such a uniform lower bound is new in algorithmic robust statistics.

We now state a more fundamental version, whose conditions are those shown to hold with high probability in Sections~\ref{s:robust:averages:counting} and~\ref{s:counting:truncation:in:probability}.

\begin{definition}[Favorable sample via counting--truncation]
\label{def:favorable:counting}
Suppose Assumption \ref{assump:finite:covariance:contamination} holds. Let $k\in[n]$ with $\varepsilon n\le k<n$, and $B,\bar r_1>0$. The sample
$\bx_{1:n}$ is \emph{$(B,\bar r_1,\varepsilon_{n,k})$-favorable via counting--truncation} if
\begin{enumerate}
    \item[\rm (i)] $\Count_2\big(\bxi_{1:n},B,k-\lfloor\varepsilon n\rfloor\big)$ holds,
    \item[\rm (ii)]
    $
      \sup_{\bv\in\mbS^{d-1}}\left|
      \frac{1}{n}\sum_{i=1}^n{\sf clip}_{B}(\langle\bxi_i,\bv\rangle) - \esp[{\sf clip}_{B}(\langle\bxi,\bv\rangle)]
    \right| \le \bar r_1,
    $
    \item[\rm (iii)]
    $
    \sup_{\bfM\in\calM_d}\left|
      \frac{1}{n}\sum_{i=1}^n{\sf clip}_{B}^2(\|\bfC_{\bfM}\bxi_i\|_2) - \esp[{\sf clip}_{B}^2(\|\bfC_{\bfM}\bxi\|_2)]
    \right| \le B\bar r_1.
    $
  \end{enumerate}
\end{definition}

Combining Lemmas~\ref{lemma:counting} and~\ref{lemma:truncated:processes}, the next corollary, proved in Appendix~\ref{app:favorable:event}, shows that Definition~\ref{def:favorable:counting} holds with high probability.
\begin{corollary}[Favorability via counting--truncation with high probability]\label{cor:favorable:event}
Suppose that Assumption~\ref{assump:finite:covariance:contamination} holds. Let $\delta\in(0,1)$ and $k\in[n]$ with $\varepsilon n\le k<n$, and assume that $s:=k-\lfloor\varepsilon n\rfloor$ satisfies $s\ge\log(3/\delta)$. Let 
$
  B \ge\max\{2\nu_p(\frac{16\,n}{s})^{1/p},
    512\sqrt{2\pi}\,\calC_n\,(\frac ns)\}
$
and 
$
  \bar r_1(\delta) :=\frac{224\sqrt{2\pi}}{3}\,\calC_n
    +\frac{4\sqrt2}{3}\,\nu_2\sqrt{\frac{\log(3/\delta)}{n}}
    +\frac{8Bs}{3n}.
$
Let ${\rm E}(\delta)$ be the intersection of the event $\Count_2(\bxi_{1:n},B,s)$ from Lemma~\ref{lemma:counting} and the events ${\rm E}_1(\delta/3)$ and ${\rm E}_2(\delta/3)$ from Lemma~\ref{lemma:truncated:processes}. Then $\prob({\rm E}(\delta))\ge1-\delta$, and, on ${\rm E}(\delta)$, $\bx_{1:n}$ is $(B,\bar r_1(\delta),\varepsilon_{n,k})$-favorable via counting--truncation.
\end{corollary}

In Section~\ref{ss:favorable}, we show that a sample that is favorable via counting--truncation is also favorable with a proper radius. In Section~\ref{ss:stability} we show that a sample that is favorable via counting--truncation is stable. A key step is item~(a) of Definition~\ref{def:favorable}, the spectral signature property, established in Section \ref{ss:spectral:signature}.

\subsection{Spectral signature}\label{ss:spectral:signature}

The next lemma establishes the spectral signature property.
\begin{lemma}[Spectral signature]\label{lemma:spectral:signature}
  Let $k\in[n]$ with $\varepsilon n\le k<n$ and $B,\bar r_1>0$ and suppose the sample $\bx_{1:n}$ is \emph{$(B,\bar r_1,\varepsilon_{n,k})$-favorable via counting-truncation}. Define $r_1':=4B\varepsilon_{n,k} + \bar r_1 + \frac{\nu_p^p}{B^{p-1}}$ and $r'_2 := 2B^2\varepsilon_{n,k} + B\bar r_1 + \frac{\nu_p^p}{B^{p-2}}$.
  Then for all $w\in\Delta_{n,k}$ and $\bv\in\mbS^{d-1}$, 
  \begin{align}
  |\langle \bv,\bmmu_{w}-\bmmu\rangle|
  &\le r_{1}' 
  + \sqrt{\frac{\varepsilon_{n,k}}{1-\varepsilon_{n,k}}\left(
     \langle\bv,[\bfSigma_w(\bmmu)-\bfSigma]\bv\rangle + r_2'
  \right)}.
\end{align}
\end{lemma} 
\begin{proof}
Fix $\bv\in\mbS^{d-1}$. We have:
\begin{align}
  \langle \bv,\bmmu_{w}-\bmmu\rangle = \sum_{i=1}^nw_i\,\langle \tilde\bxi_i,\bv\rangle 
  &= \sum_{i=1}^nw_i\,{\sf clip}_B(\langle \tilde\bxi_i,\bv\rangle) 
  + \sum_{i=1}^nw_i\,{\sf rem}_B(\langle \tilde\bxi_i,\bv\rangle) \\ 
  & \le \left|\sum_{i=1}^nw_i\,{\sf clip}_B(\langle \tilde\bxi_i,\bv\rangle) 
  - \frac{1}{n}\sum_{i=1}^n{\sf clip}_B(\langle \bxi_i,\bv\rangle)\right| \\ 
  & + \left|\frac{1}{n}\sum_{i=1}^n{\sf clip}_B(\langle \bxi_i,\bv\rangle) 
  - \esp[{\sf clip}_B(\langle \bxi,\bv\rangle)]\right| \\ 
  & + |\esp[{\sf clip}_B(\langle \bxi,\bv\rangle)]| \\ 
  & + \sum_{i=1}^nw_i\,|{\sf rem}_B(\langle \tilde\bxi_i,\bv\rangle)|\\
  &=: (I) + (II) + (III) + (IV).
\end{align}

To bound $(I)$, we split it as
\begin{align}
  (I)&\le\left|\sum_{i=1}^n\left(w_i-\frac1n\right){\sf clip}_B(\langle\tilde\bxi_i,\bv\rangle)\right|
  +\Big|\frac1n\sum_{i=1}^n\big({\sf clip}_B(\langle\tilde\bxi_i,\bv\rangle)
  -{\sf clip}_B(\langle\bxi_i,\bv\rangle)\big)\Big|.
\end{align}
Denote by $(I_1)$ and $(I_2)$ the two sums above, respectively. 

Lemma \ref{lemma:weight:dif} applied to $a_i:={\sf clip}_B(\langle\tilde\bxi_i,\bv\rangle)$ gives $(I_1) \le \frac{2Bk}{n}$. 
By Assumption~\ref{assump:finite:covariance:contamination}, $\tilde\bxi_i=\bxi_i$ for all but at most $\varepsilon n\le k$ indices, and each remaining term is at most $2B$ in absolute value, so $(I_2)\le\frac{2Bk}{n}$. By item~(ii) in Definition \ref{def:favorable:counting}, $(II)\le \bar r_1$, and by Lemma~\ref{lemma:bias}, $(III)\le\nu_p^p/B^{p-1}$. We obtain:
\begin{align}
  \langle \bv,\bmmu_{w}-\bmmu\rangle \le \frac{4Bk}{n} + \bar r_1 + \frac{\nu_p^p}{B^{p-1}} + (IV) =: r'_1 + (IV).  
  \label{lemma:spectral:signature:eq1}
\end{align}

We now claim that 
\begin{align}\label{eq:claimproperty}
  \mbox{\underline{Claim:}}\;\;(IV) \le \sqrt{\frac{k}{n-k}\left(
    r_2' + \langle\bv,[\bfSigma_w(\bmmu)-\bfSigma]\bv\rangle
  \right)}.
\end{align}
If we prove the claim for arbitrary $\bv\in\mbS^{d-1}$, by \eqref{lemma:spectral:signature:eq1}, we prove the lemma (by using it for $-\bv$ as well).

To prove this claim, define the set 
$
  S_B:=\{i\in[n]\,:\, {\sf rem}_B(\langle \tilde\bxi_i,\bv\rangle)\neq 0\} = \{i\in[n]\,:\,|\langle \tilde\bxi_i,\bv\rangle|>B\}.
$
By item (i) in Definition \ref{def:favorable:counting} and Assumption \ref{assump:finite:covariance:contamination}, $\Count_2(\tilde\bxi_{1:n},B,k)$ holds. By this fact and the fact that $|\langle \tilde\bxi_i,\bv\rangle|=\|\bfC_{\bv^{\otimes2}}\tilde\bxi_i\|_2$ for any $\bv\in\mbS^{d-1}$, we have that $\# S_{B}\le k$. We obtain:
\begin{align}
  (IV) &= \sum_{i=1}^nw_i\,|{\sf rem}_B(\langle \tilde\bxi_i,\bv\rangle)|
    = \sum_{i\in S_B}w_i\,|{\sf rem}_B(\langle \tilde\bxi_i,\bv\rangle)|\\ 
    &\leq  \sqrt{\sum_{i\in S_B}w_i}\sqrt{\sum_{i=1}^nw_i\,{\sf rem}_B^2(\langle \tilde\bxi_i,\bv\rangle)}
    \leq  \sqrt{\frac{k}{n-k}}\sqrt{\sum_{i=1}^nw_i\,{\sf rem}_B^2(\langle \tilde\bxi_i,\bv\rangle)}.
    \label{lemma:spectral:signature:eq2}
\end{align}
Above, we used Cauchy-Schwarz in the first inequality and the fact that $\sum_{i\in S_B}w_i\le|S_B|\cdot\frac1{n-k}\le\frac{k}{n-k}$ since $w\in\Delta_{n,k}$ in the second inequality. 

Now, for any $t\in\re$, ${\sf rem}_B^2(t) = (|t|-|{\sf clip}_B(t)|)^2\leq t^2 -{\sf clip}_B^2(t)$. Therefore, 
\begin{align}
  \sum_{i=1}^nw_i\,{\sf rem}_B^2(\langle \tilde\bxi_i,\bv\rangle) \leq  
  \sum_{i=1}^nw_i \langle \tilde\bxi_i,\bv\rangle^2 -\sum_{i=1}^nw_i\,{\sf clip}_B^2(\langle \tilde\bxi_i,\bv\rangle).
  \label{lemma:spectral:signature:eq3}
\end{align}

Note that 
\begin{align}
  \sum_{i=1}^nw_i \langle \tilde\bxi_i,\bv\rangle^2  &= \langle \bv,\bfSigma_{w}(\bmmu)\bv\rangle
  =  \langle \bv,[\bfSigma_{w}(\bmmu)-\bfSigma]\bv\rangle 
  + \langle \bv,\bfSigma\bv\rangle.
  \label{lemma:spectral:signature:eq4}
\end{align}

From Lemma \ref{lemma:weight:dif} applied to $a_i:={\sf clip}_B^2(\langle\tilde\bxi_i,\bv\rangle)$, 
\begin{align}
  -\sum_{i=1}^nw_i\,{\sf clip}_B^2(\langle \tilde\bxi_i,\bv\rangle)  
  &\le -\frac{1}{n}\sum_{i=1}^n{\sf clip}_B^2(\langle \tilde\bxi_i,\bv\rangle) 
  + \frac{B^2k}{n}\\ 
  &\le -\frac{1}{n}\sum_{i=1}^n{\sf clip}_B^2(\langle \bxi_i,\bv\rangle) 
  + \frac{2B^2k}{n} \\
  &\le -\esp[{\sf clip}_B^2(\langle \bxi_1,\bv\rangle)] + B\bar r_1 + \frac{2B^2k}{n}\\ 
  &\le -\langle \bv,\bfSigma\bv\rangle 
  + \frac{\nu_p^p}{B^{p-2}} + B\bar r_1 + \frac{2B^2k}{n}.
  \label{lemma:spectral:signature:eq5}
\end{align}
In the second inequality above we used that at most $\varepsilon n\le k$ terms change, and each change lies in $[-B^2,B^2]$, so this costs another $\frac{B^2k}{n}$. The third inequality follows from item (iii) in Definition \ref{def:favorable:counting} and the fourth inequality follows from Lemma \ref{lemma:bias}. 

Combining \eqref{lemma:spectral:signature:eq3}, \eqref{lemma:spectral:signature:eq4} and \eqref{lemma:spectral:signature:eq5}, we obtain
$
  \sum_{i=1}^nw_i\,{\sf rem}_B^2(\langle \tilde\bxi_i,\bv\rangle) 
  \leq \langle \bv,[\bfSigma_{w}(\bmmu)-\bfSigma]\bv\rangle 
  + r_2',
$
which we can plug back into \eqref{lemma:spectral:signature:eq2} to finish the proof of the claim.
\end{proof}

\subsection{Favorability}\label{ss:favorable}

Lemma~\ref{lemma:spectral:signature} yields the following.

\begin{proposition}[Favorable via counting--truncation $\Rightarrow$ favorable]\label{prop:counting:favorable}
Suppose Assumption \ref{assump:finite:covariance:contamination} holds. Let $k\in[n]$ with $\varepsilon n\le k<n$, and $B,\bar r_1>0$ and suppose that $\bx_{1:n}$ is \emph{$(B,\bar r_1,\varepsilon_{n,k})$-favorable via counting--truncation}. Then $\tilde\bx_{1:n}$ is \emph{$(r_1,\varepsilon_{n,k})$-favorable} with 
\begin{align}
  r_1 := \left(
    4B\varepsilon_{n,k} + \bar r_1 + \frac{\nu_p^p}{B^{p-1}}
  \right)
  \bigvee 
  \left(
    \varepsilon_{n,k}\left(
      2B^2\varepsilon_{n,k} + B\bar r_1 + \frac{\nu_p^p}{B^{p-2}}
    \right)
  \right)^{\frac{1}{2}}. 
  \label{prop:counting:favorable:r:1}
\end{align}
\end{proposition}
\begin{proof}
Items (i)-(iii) in Definition \ref{def:favorable:counting} and Lemma \ref{lemma:spectral:signature} imply immediately that, for any $w\in\Delta_{n,k}$,
\begin{align}
  \|\bmmu_{w}-\bmmu\|_2 &\le r_{1}' 
  + \sqrt{\frac{\varepsilon_{n,k}}{1-\varepsilon_{n,k}}
    \left(
      \lambda_1(\bfSigma_w(\bmmu)-\bfSigma) + r_2'
    \right)},\label{lemma:favorable:counting:truncation:to:favorable:eq1} 
\end{align}
where $r_1':=\frac{4Bk}{n} + \bar r_1 + \frac{\nu_p^p}{B^{p-1}}$ and $r'_2 := \frac{2B^2k}{n} + B\bar r_1 + \frac{\nu_p^p}{B^{p-2}}$. 

Let $w^*\in\argmin_{w\in\Delta_{n,k}}\lambda_1(\bfSigma_w(\bmmu)-\bfSigma)$. Since $\Delta_{n,k}$ and $\calM_d$ are compact and convex, and $(w,\bfM)\mapsto\langle\bfSigma_w(\bmmu)-\bfSigma,\bfM\rangle$ is bilinear, von Neumann’s minimax theorem allows us to exchange the order of $\min$ and $\sup$: 
\begin{align}
  \lambda_1(\bfSigma_{w^*}(\bmmu)-\bfSigma)
  = \sup_{\bfM\in\calM_d}\min_{w\in\Delta_{n,k}}
  \langle\bfSigma_w(\bmmu)-\bfSigma,\bfM\rangle.
\end{align}
From Lemma \ref{lemma:min:weight:trimming}, 
\begin{align}
  \sup_{\bfM\in\calM_d}\min_{w\in\Delta_{n,k}}
  \langle\bfSigma_w(\bmmu)-\bfSigma,\bfM\rangle = 
  \sup_{\bfM\in\calM_d} \min_{S\subseteq[n],\,|S|=n-k}\frac1{n-k}\sum_{i\in S}\langle\tilde\bxi_i^{\otimes2}-\bfSigma,\bfM\rangle.
\end{align}
By item (i) in Definition \ref{def:favorable:counting} and Assumption \ref{assump:finite:covariance:contamination}, $\Count_2(\tilde\bxi_{1:n},B,k)$ holds. By this fact and Lemma \ref{lemma:trim:trunc},
\begin{align}
  &\sup_{\bfM\in\calM_d} \min_{S\subseteq[n],\,|S|=n-k}\frac1{n-k}\sum_{i\in S}\langle\tilde\bxi_i^{\otimes2}-\bfSigma,\bfM\rangle\\
  &\le \sup_{\bfM\in\calM_d}\frac1n\sum_{i=1}^n
  \left(
    {\sf clip}_{B}^2(\|\bfC_{\bfM}\tilde\bxi_i\|_2) - \langle\bfSigma,\bfM\rangle
  \right) =: (I). 
\end{align}
We can write 
\begin{align}
  (I) &\le \sup_{\bfM\in\calM_d}\frac1n\sum_{i=1}^n
  \left(
    {\sf clip}_{B}^2(\|\bfC_{\bfM}\tilde\bxi_i\|_2) 
    - {\sf clip}_{B}^2(\|\bfC_{\bfM}\bxi_i\|_2)
  \right) \\
  &+ \sup_{\bfM\in\calM_d}\frac1n\sum_{i=1}^n
  \left(
    {\sf clip}_{B}^2(\|\bfC_{\bfM}\bxi_i\|_2) 
    - \esp[{\sf clip}_{B}^2(\|\bfC_{\bfM}\bxi\|_2)]
  \right) =: (I_1) + (I_2),\label{lemma:favorable:counting:truncation:to:favorable:eq2} 
\end{align}
where we used that 
$
  \esp\big[{\sf clip}_B^2(\|\bfC_{\bfM}\bxi\|_2)\big]\le\langle\bfSigma,\bfM\rangle
$
by Lemma \ref{lemma:bias}. By Assumption~\ref{assump:finite:covariance:contamination},
$\tilde\bxi_i=\bxi_i$ for all but at most $\varepsilon n\le k$ indices, and each remaining term of $(I_1)$ lies in $[-B^2,B^2]$, since both truncated values lie in $[0,B^2]$; hence $(I_1)\le\frac{B^2k}{n}$. By item~(iii) of Definition~\ref{def:favorable:counting}, $(I_2)\le B\bar r_1$. We conclude that 
\begin{align}
   \lambda_1\big(\bfSigma_{w^*}(\bmmu)-\bfSigma\big)\le\frac{B^2k}{n}+B\bar r_1.
   \label{lemma:favorable:counting:truncation:to:favorable:eq3}
\end{align}

It remains to show the uniform lower bound over weights in $\Delta_{n,k}$. We have that 
\begin{align}
  \min_{w\in\Delta_{n,k}}\lambda_d(\bfSigma_{w}(\bmmu)-\bfSigma)
  = \min_{\bfM\in\calM_d}\min_{w\in\Delta_{n,k}}
  \langle\bfSigma_w(\bmmu)-\bfSigma,\bfM\rangle.
\end{align}
From Lemma \ref{lemma:min:weight:trimming}, 
\begin{align}
  \min_{\bfM\in\calM_d}\min_{w\in\Delta_{n,k}}
  \langle\bfSigma_w(\bmmu)-\bfSigma,\bfM\rangle = 
  \min_{\bfM\in\calM_d} \min_{S\subseteq[n],\,|S|=n-k}\frac1{n-k}\sum_{i\in S}\langle\tilde\bxi_i^{\otimes2}-\bfSigma,\bfM\rangle.
\end{align}
The validity of $\Count_2(\tilde\bxi_{1:n},B,k)$ and Lemma \ref{lemma:trim:trunc} imply
\begin{align}
  &\min_{\bfM\in\calM_d} \min_{S\subseteq[n],\,|S|=n-k}\frac1{n-k}\sum_{i\in S}\langle\tilde\bxi_i^{\otimes2}-\bfSigma,\bfM\rangle\\
  &\ge \min_{\bfM\in\calM_d}\frac1n\sum_{i=1}^n
  \left(
    {\sf clip}_{B}^2(\|\bfC_{\bfM}\tilde\bxi_i\|_2) - \langle\bfSigma,\bfM\rangle
  \right) 
  - \frac{B^2k}{n} =: (II). 
\end{align}

Proceeding similarly with the arguments in \eqref{lemma:favorable:counting:truncation:to:favorable:eq2} 
but using the fact that 
\[
  \esp\big[{\sf clip}_B^2(\|\bfC_{\bfM}\bxi\|_2)\big]\ge\langle\bfSigma,\bfM\rangle
  -\frac{\nu_p^p}{B^{p-2}}
\]
by Lemma \ref{lemma:bias} and using the chain of equalities and inequalities above, we obtain that 
\begin{align}
  \min_{w\in\Delta_{n,k}}\lambda_d(\bfSigma_{w}(\bmmu)-\bfSigma)
  \ge -\frac{2B^2k}{n} -B\bar r_1 -\frac{\nu_p^p}{B^{p-2}}.
  \label{lemma:favorable:counting:truncation:to:favorable:eq4}
\end{align}

We conclude from \eqref{lemma:favorable:counting:truncation:to:favorable:eq1}, \eqref{lemma:favorable:counting:truncation:to:favorable:eq3} and \eqref{lemma:favorable:counting:truncation:to:favorable:eq4} that the sample $\tilde\bx_{1:n}$ is $(r_1,\varepsilon_{n,k})$-favorable with 
$
  r_1 = r'_1\vee \left( \varepsilon_{n,k}r_2' \right)^{\frac{1}{2}}.
$
This finishes the proof. 
\end{proof}

Next, we prove:
\begin{theorem}[Favorability with the optimal sub-Gaussian rate]\label{thm:favorability:subgaussian:rate}
  Suppose that Assumption~\ref{assump:finite:covariance:contamination} holds. Let $\delta\in(0,1)$ and $\eta\in(0,1]$ and suppose that 
\begin{align}
  k:=\lceil(1+\eta)\varepsilon n\rceil+\lceil\log(3/\delta)\rceil < n.
  \label{eq:choice:k}
\end{align}  
Then, on the event ${\rm E}(\delta)$ of Corollary~\ref{cor:favorable:event}, the sample $\tilde\bx_{1:n}$ is $(O(r),\varepsilon_{n,k})$-favorable with $r=r(n,\delta,\varepsilon)$ given by 
\begin{align}
  r(n,\delta,\varepsilon) := \eta^{-1}\sqrt{\frac{\tr(\bfSigma)}{n}}
  +\eta^{-\frac12}\nu_2\sqrt{\frac{\log(3/\delta)}{n}}
  +\eta^{-\frac12}\inf_{p\ge2}\nu_p\,\varepsilon^{1-\frac1p}.
  \label{eq:subgaussian:rate:alpha:v2}
\end{align}
\end{theorem}
\begin{proof}
  The integer $s:=k-\lfloor\varepsilon n\rfloor$ satisfies $\eta\varepsilon n+\log(3/\delta)\le s\le\eta\varepsilon n+\log(3/\delta)+3$; in particular, $\varepsilon n\le k<n$ and $s\ge\log(3/\delta)$, as required in Corollary~\ref{cor:favorable:event}, and $k/s\le4/\eta$. Define
\begin{align}
  B &:= \max\Big\{2\nu_p\Big(\frac{16\,n}{s}\Big)^{1/p},\
  512\sqrt{2\pi}\,\calC_n\,\frac ns\Big\},
  \label{eq:main:B}\\
  \bar r_1 := \bar r_1(\delta)&:=\frac{224\sqrt{2\pi}}{3}\,\calC_n
  +\frac{4\sqrt2}{3}\,\nu_2\sqrt{\frac{\log(3/\delta)}{n}}
  +\frac{8Bs}{3n},
  \label{eq:main:r:1:delta}
\end{align}
and let ${\rm E}(\delta)$ be the event of Corollary~\ref{cor:favorable:event} for these $B$ and $\bar r_1$. On ${\rm E}(\delta)$, the $\bx_{1:n}$ is $(B,\bar r_1,\varepsilon_{n,k})$-favorable via counting--truncation, so, by Proposition~\ref{prop:counting:favorable}, $\tilde\bx_{1:n}$ is $(r_1,\varepsilon_{n,k})$-favorable, with $r_1$ given by~\eqref{prop:counting:favorable:r:1}.

Since $B\ge2\nu_p(16n/s)^{1/p}$ and $s\le k$, the truncation bias satisfies $\nu_p^p/B^{p-1}\le B\varepsilon_{n,k}/64$. Hence the second term in the maximum defining $r_1$ is at most $2.05\,B\varepsilon_{n,k}+\frac12\bar r_1$, which is smaller than the first term, so
\begin{align}
  r_1=4B\varepsilon_{n,k}+\bar r_1+\frac{\nu_p^p}{B^{p-1}}
  \le \bar r_1+\frac{257}{64}\,B\varepsilon_{n,k}.
  \label{eq:favorable:radius}
\end{align}
By~\eqref{eq:main:B}--\eqref{eq:favorable:radius}, $k/s\le4/\eta$, $\eta^{-1/p}\le\eta^{-1/2}$, $\varepsilon_{n,k}\le\varepsilon+s/n$ and the well known bound $\calC_n\le\sqrt{\tr(\bfSigma)/n}$,
\begin{align}
    r_1 \lesssim \eta^{-1}\sqrt{\frac{\tr(\bfSigma)}{n}}
    +\nu_2\sqrt{\frac{\log(3/\delta)}{n}}
    +\eta^{-1/2}\nu_p\Big(\varepsilon^{1-1/p}
    +\Big(\frac{s}{n}\Big)^{1-1/p}\Big).
\end{align}

Since $s\le k$, it suffices to bound $k$. The exponent $p$, and hence $B$ and ${\rm E}(\delta)$, is fixed in advance: if $\varepsilon n\ge\log(3/\delta)$, we choose $p\ge2$ such that $\nu_p\varepsilon^{1-1/p}\le2\inf_{q\ge2}\nu_q\varepsilon^{1-1/q}$; otherwise, we choose $p=2$. In the first case, $k\lesssim\varepsilon n$, so the last term is at most of order $\eta^{-1/2}\nu_p\varepsilon^{1-1/p}\lesssim\eta^{-1/2}\inf_{q\ge2}\nu_q\varepsilon^{1-1/q}$. In the second case, $k\lesssim\log(3/\delta)$, so the last term is at most of order $\eta^{-1/2}\nu_2\sqrt{\log(3/\delta)/n}$. In both cases, $r_1\lesssim r$ on ${\rm E}(\delta)$, which is the optimal sub-Gaussian rate up to factors depending on $\eta$.

\end{proof}

\subsection{Stability}\label{ss:stability}

In this subsection, we show that a sample that is favorable via counting--truncation is stable with respect to its good weights; the proof is in Appendix~\ref{app:favorable:stable}.
\begin{proposition}[Favorable via counting--truncation $\Rightarrow$ stable]\label{prop:favorable:stable}
Suppose Assumption \ref{assump:finite:covariance:contamination} holds. Let $k\in[n]$ with $\varepsilon n\le k<n/2$ and $B,\bar r_1>0$, and suppose that $\bx_{1:n}$ is $(B,\bar r_1,\varepsilon_{n,k})$-favorable via counting--truncation. Let $r_1$ given by~\eqref{prop:counting:favorable:r:1} and 
\begin{align}
  R_1 := \left(
    4B\varepsilon_{n,2k} + \bar r_1 + \frac{\nu_p^p}{B^{p-1}}
  \right)
  \bigvee 
  \left(
    \varepsilon_{n,2k}\left(
      2B^2\varepsilon_{n,2k} + B\bar r_1 + \frac{\nu_p^p}{B^{p-2}}
    \right)
  \right)^{\frac{1}{2}}. 
  \label{prop:favorable:stable:R:1}
\end{align}
Then, for every $w^*\in\argmin_{w\in\Delta_{n,k}}\lambda_1(\bfSigma_w(\bmmu)-\bfSigma)$, the sample $\tilde\bx_{1:n}$ is $(R,\varepsilon_{n,k})$-stable with respect to $w^*$, with
\begin{align}
  R:=R_1+\bigg(\frac{1}{1-\varepsilon_{n,2k}}\Big(R_1^2
  +\frac{2(r_1^2+\varepsilon_{n,k}^2\|\bfSigma\|)}{1-\varepsilon_{n,k}}\Big)\bigg)^{1/2}.
  \label{prop:favorable:stable:r}
\end{align}
\end{proposition}

Next, we prove:
\begin{theorem}[Stability with the optimal sub-Gaussian rate]\label{thm:stability:subgaussian:rate}
  Suppose the assumptions of Theorem \ref{thm:favorability:subgaussian:rate} hold and $k$, as in \eqref{eq:choice:k}, satisfies $k<n/2$. Then, on the event ${\rm E}(\delta)$ of Corollary~\ref{cor:favorable:event}, the sample $\tilde\bx_{1:n}$ is $(R,\varepsilon_{n,k})$-stable with respect to every $w^*\in\argmin_{w\in\Delta_{n,k}}\lambda_1(\bfSigma_w(\bmmu)-\bfSigma)$ with $R=O(r)(1-2\varepsilon_{n,k})^{-1/2}$ with $r=r(n,\delta,\varepsilon)$ given by \eqref{eq:subgaussian:rate:alpha:v2}.
\end{theorem}
\begin{proof}
  We refer to the proof of Theorem \ref{thm:favorability:subgaussian:rate} and definitions therein. If, moreover, $k<n/2$, then, since on the event ${\rm E}(\delta)$ the sample $\bx_{1:n}$ is $(B,\bar r_1,\varepsilon_{n,k})$-favorable via counting--truncation, Proposition~\ref{prop:favorable:stable} shows that the sample $\tilde\bx_{1:n}$ is $(R,\varepsilon_{n,k})$-stable with respect to every $w^*\in\argmin_{w\in\Delta_{n,k}}\lambda_1(\bfSigma_w(\bmmu)-\bfSigma)$, with $R$ given by~\eqref{prop:favorable:stable:r} and $R_1$ by~\eqref{prop:favorable:stable:R:1}. Since $1-\varepsilon_{n,k}>\frac12$, $R_1\le2\sqrt2\,r_1$, and $\varepsilon_{n,k}\sqrt{\|\bfSigma\|}\le\frac18 r_1$ (because $\nu_p\ge\nu_2=\sqrt{\|\bfSigma\|}$ and $r_1\ge4B\varepsilon_{n,k}$), we have
  $
    \sqrt2\,r_1\le R\lesssim\frac{r_1}{\sqrt{1-2\varepsilon_{n,k}}}. 
  $
  Hence, for $\varepsilon_{n,k}$ bounded away from $\frac12$, \eqref{eq:subgaussian:rate:alpha:v2} also bounds the stability radius $R$, up to an absolute constant.  
\end{proof}

\section{Optimization landscape of robust mean estimation}\label{s:landscape}

Recall definitions \eqref{eq:landscape:problem}-\eqref{eq:surrogate:spectral:function} of the functions $F$ and $F^*$ in Section \ref{s:main:results}. Throughout this section, $k\in[n]$ with $k<n$ and $r_1>0$ are fixed, and $\tilde\bx_{1:n}$ is $(r_1,\varepsilon_{n,k})$-favorable. Without loss of generality, we take its good weights to be $w^*\in\argmin_{w\in\Delta_{n,k}}\lambda_1(\bfSigma_w(\bmmu)-\bfSigma)$, since the minimum is attained and any minimizer satisfies item~(b). Items~(b) and~(c) then give
\begin{align}
  -\frac{r_1^2}{\varepsilon_{n,k}}\le F^* = \lambda_1(\bfSigma_{w^*}(\bmmu)-\bfSigma)
  \le\frac{r_1^2}{\varepsilon_{n,k}},
  \label{def:w:*:r:1}
\end{align}
and item~(a) at $w^*$, together with~\eqref{def:w:*:r:1}, gives
\begin{align}
  \|\bmmu_{w^*}-\bmmu\|_2\le\left(1+\sqrt{\frac{2}{1-\varepsilon_{n,k}}}\right)r_1
  =:r_{1,*}.\label{def:r:1:*}
\end{align}
Fix $\alpha\ge0$ and define
\begin{align}
  R_{\alpha} := \max_{\omega\in\Delta_{n,k}}\left\{
    \sqrt{\alpha}\|\omega-w^*\|_2
  \right\}
  ~~\mbox{ and }~~
  r_{1,\alpha} := 2r_{1,*} + R_{\alpha}.\label{def:R:alpha:r:1:alpha}
\end{align}
The error of centralization to $\bfSigma$ in the surrogate enters only through
\begin{align}
    \tau:=\lambda_1(\hat\bfSigma-\bfSigma)+\lambda_1(\bfSigma-\hat\bfSigma) \ge0,
  \label{eq:tau}
\end{align}
the difference between the largest and smallest eigenvalues of $\hat\bfSigma-\bfSigma$. Note that $\tau\le2\|\hat\bfSigma-\bfSigma\|$, and $\tau=0$ whenever $\hat\bfSigma-\bfSigma$ is a multiple of $\bfI_d$. 

\subsection{Relative inexact quadratic growth, via spectral signature}

Next, we show the objective $F$ satisfies tht inexact quadratic growth property relative to the mean-induced divergence with quadrature constant $\frac{n-k}{k}=\frac{1-\varepsilon_{n,k}}{\varepsilon_{n,k}}$ and radius $r_{1,\alpha}$, for any $\alpha\ge0$. See  Definition \ref{def:relative:iqg} in Section \ref{ss:relative:regularity}.

\begin{proposition}[Spectral signature $\Rightarrow$ mean-induced $\IQG$]\label{prop:IQG}
For any $w\in\Delta_{n,k}$, 
\begin{align}
  \|\bmmu_{w}-\bmmu\|_2 \le r_{1,*}+\sqrt{\frac{\varepsilon_{n,k}}{1-\varepsilon_{n,k}}(F(w)-F^*)}.\label{eq:rates:signature}
\end{align}
In addition, for any $\alpha\ge0$, $F$ satisfies $\IQG_{h_{\mbX,\alpha}}(\frac{n-k}{k},r_{1,\alpha})$ on $\Delta_{n,k}$, i.e., 
\begin{align}
    F(w) \ge F(w^*) + \frac{n-k}{k}\Big(\sqrt{2D_{h_{\mbX,\alpha}}(w^*\|w)}-r_{1,\alpha}\Big)_+^2
    \qquad\text{for all } w\in\Delta_{n,k}.
\end{align}
\end{proposition}
\begin{proof}
  Fix $w\in\Delta_{n,k}$. By item~(a) of Definition~\ref{def:favorable}, $F^*\le r_1^2/\varepsilon_{n,k}$ and $\sqrt{x+y}\le\sqrt x+\sqrt y$,
\[
  \|\bmmu_{w}-\bmmu\|_2
  \le r_1+\sqrt{\frac{\varepsilon_{n,k}}{1-\varepsilon_{n,k}}\Big(F(w)-F^*+\frac{2r_1^2}{\varepsilon_{n,k}}\Big)}
  \le r_{1,*}+\sqrt{\frac{\varepsilon_{n,k}}{1-\varepsilon_{n,k}}(F(w)-F^*)}.
\]
This proves \eqref{eq:rates:signature}. 

From this fact and the bound \eqref{def:r:1:*} on $\|\bmmu_{w^*}-\bmmu\|_2$,
\begin{align}
  \|\bmmu_w-\bmmu_{w^*}\|_2 &\le\|\bmmu_w-\bmmu\|_2+\|\bmmu_{w^*}-\bmmu\|_2\\
  &\le 2r_{1,*}+\sqrt{\frac{\varepsilon_{n,k}}{1-\varepsilon_{n,k}}
    \big(F(w)-F^*\big)}.\label{eq:spectral:centrality:at:bmmu:w:*}
\end{align}

Adding $d_{\alpha}(w^*,w) := \sqrt{\alpha}\|w-w^*\|_2$ to both sides of the above inequality, using that $\sqrt{x+y}\le \sqrt{x}+\sqrt{y}$ for any $x,y\ge0$ and the fact that $d_{\alpha}(w^*,w) \le R_{\alpha}$,
\begin{align}
    \sqrt{\|\bmmu_w-\bmmu_{w^*}\|_2^2
    + d^2_{\alpha}(w^*,w)}
    &\le \|\bmmu_w-\bmmu_{w^*}\|_2
    + d_{\alpha}(w^*,w)\\
    &\le 2r_{1,*} + R_{\alpha} + \sqrt{\frac{\varepsilon_{n,k}}{1-\varepsilon_{n,k}}\left(F(w)-F^*\right)}.
\end{align}
This inequality, together with the facts that $\|\bmmu_{w}-\bmmu_{w^*}\|_2^2 = 2D_{\mbX}(w^*\Vert w)$ and $d_{\alpha}^2(w^*,w) = 2\alpha D_{\E}(w^*\Vert w)$, yields:
\begin{align}
    \Big(
    \sqrt{
      2D_{\mbX}(w^*\Vert w)
      + 2\alpha D_{\E}(w^*\Vert w)
    }
    - 2r_{1,*} - R_\alpha
    \Big)_+
    \le \sqrt{\frac{\varepsilon_{n,k}}{1-\varepsilon_{n,k}}
    \left(F(w)-F^*\right)}.
\end{align}
Squaring the inequality above finishes the proof.
\end{proof}

\subsection{Relative inexact subgradient oracle}

Each $\hat f(\cdot,\bfM)$ is a concave quadratic, so $\hat F$ is weakly convex, and its Fréchet subdifferential $\partial\hat F(w)$ is nonempty and coincides with the Clarke subdifferential for every $w\in\re^n$~\cite{1998Rockafellar:Wets,2019Davis:Drusvyatskiy}. Its Euclidean weak-convexity modulus, however, grows with $n$, $d$ and the size of the outliers \cite{2020Cheng:Diakonikolas:Ge:Soltanolkotabi}. The next proposition shows that the mean-induced geometry removes this dependence: relative to $h_{\mbX}$, $\partial\hat F$ is a $(4,r_{1,*},\tau)$-weakly star-convex inexact oracle for $F$ at $w^*$. See Definition \ref{def:inexact:oracle}(b) in Section \ref{ss:relative:regularity}.

\begin{proposition}[Relative inexact weakly star-convex oracle]\label{prop:weakly:convex:oracle}
For all $w\in\Delta_n$ and $g\in\partial\hat F(w)$,
\begin{align}
  F(w^*)\ge F(w)+\langle g,w^*-w\rangle-2D_{\mbX}(w^*\|w)
  -\Big(\sqrt{2D_{\mbX}(w^*\|w)}+r_{1,*}\Big)^2 -\tau.
\end{align}
In particular, $\partial\hat F$ is a $(4,r_{1,*},\tau)$-weakly star-convex inexact oracle for $F$ at $w^*$ relative to $h_{\mbX}$.
\end{proposition}

To prove Proposition \ref{prop:weakly:convex:oracle}, we need the following lemma.
\begin{lemma}[Directional derivatives]
\label{lemma:marginal:directional:derivatives}
Extend $\hat f(\cdot,\bfM)$ to $\re^n$ by
$\hat f(w,\bfM):=\sum_iw_i\langle\tilde\bx_i^{\otimes2},\bfM\rangle
-\langle\bmmu_w^{\otimes2},\bfM\rangle-\langle\hat\bfSigma,\bfM\rangle$, with $\bmmu_w:=\sum_iw_i\tilde\bx_i$. Then, for every $\bfM\in\calM_d$,
$w\in\Delta_n$ and $w',w''\in\Delta_n$,
\begin{align}
    \big\langle\nabla_1\hat f(w,\bfM),w'-w''\big\rangle
    =\big\langle\bfSigma_{w'}(\bmmu_w)-\bfSigma_{w''}(\bmmu_w),\bfM\big\rangle.
\end{align}
\end{lemma}
\begin{proof}
Since $\hat f(\cdot,\bfM)$ is a polynomial on $\re^n$,
$[\nabla_1\hat f(w,\bfM)]_i=\langle\tilde\bx_i^{\otimes2},\bfM\rangle
-2\langle\bfM\bmmu_w,\tilde\bx_i\rangle$. On the other hand, with
$\delta:=w'-w''$, so that $\sum_i\delta_i=0$,
\begin{align}
    \big\langle\bfSigma_{w'}(\bmmu_w)-\bfSigma_{w''}(\bmmu_w),\bfM\big\rangle
    &=\sum_i\delta_i\big\langle(\tilde\bx_i-\bmmu_w)^{\otimes2},\bfM\big\rangle\\
    &=\sum_i\delta_i\Big(\langle\tilde\bx_i^{\otimes2},\bfM\rangle
    -2\langle\bfM\bmmu_w,\tilde\bx_i\rangle\Big)
    +\langle\bmmu_w^{\otimes2},\bfM\rangle\sum_i\delta_i\\
    &=\big\langle\nabla_1\hat f(w,\bfM),\delta\big\rangle.
\end{align}

\end{proof}

\begin{proof}[Proof of Proposition~\ref{prop:weakly:convex:oracle}]
Fix $w\in\Delta_n$. We first relate $F(w)$ and $\hat F(w)$. By the
bias--variance decomposition,
$\bfSigma_w(\bmmu)=\bfSigma_w+(\bmmu_w-\bmmu)^{\otimes2}$. Hence
\begin{align}
    F(w)
    &=\max_{\bv\in\mbS^{d-1}}\big\langle\bfSigma_w(\bmmu)-\bfSigma,\bv^{\otimes2}\big\rangle\\
    &=\max_{\bv\in\mbS^{d-1}}\Big\{\big\langle\bfSigma_w-\hat\bfSigma,\bv^{\otimes2}\big\rangle
    +\big\langle\hat\bfSigma-\bfSigma,\bv^{\otimes2}\big\rangle
    +\langle\bmmu_w-\bmmu,\bv\rangle^2\Big\}\\
    &\le\hat F(w)+\lambda_1(\hat\bfSigma-\bfSigma)+\|\bmmu_w-\bmmu\|_2^2 .
    \label{prop:weakly:convex:oracle:eq1}
\end{align}

Next, we bound $\langle g,\omega-w\rangle$ for $g\in\partial\hat F(w)$ and
$\omega\in\Delta_n$; later we set $\omega:=w^*$. By Danskin's theorem, every
$g\in\partial\hat F(w)$ is a convex combination of vectors
$\nabla_1\hat f(w,\bv^{\otimes2})$ with $\bv\in\mbS^{d-1}$ such that
$\hat F(w)=\langle\bfSigma_w-\hat\bfSigma,\bv^{\otimes2}\rangle$. Let $\bv$ be such a
vector. By Lemma~\ref{lemma:marginal:directional:derivatives},
\begin{align}
    \big\langle\nabla_1\hat f(w,\bv^{\otimes2}),\omega-w\big\rangle
    =\big\langle\bfSigma_\omega(\bmmu_w)-\bfSigma_w,\bv^{\otimes2}\big\rangle .
\end{align}
By the bias--variance decomposition,
\begin{align}
    \bfSigma_\omega(\bmmu_w)
    =\bfSigma_\omega+(\bmmu_w-\bmmu_\omega)^{\otimes2}
    =\bfSigma_\omega(\bmmu)-(\bmmu_\omega-\bmmu)^{\otimes2}
    +(\bmmu_w-\bmmu_\omega)^{\otimes2} .
\end{align}
From the previous two identities, adding and subtracting $\bfSigma$ and $\hat\bfSigma$, and using
$\langle\bfSigma_\omega(\bmmu)-\bfSigma,\bv^{\otimes2}\rangle\le F(\omega)$,
$\langle\bfSigma_w-\hat\bfSigma,\bv^{\otimes2}\rangle=\hat F(w)$,
$\langle\bfSigma-\hat\bfSigma,\bv^{\otimes2}\rangle\le\lambda_1(\bfSigma-\hat\bfSigma)$ and
$-\langle\bmmu_\omega-\bmmu,\bv\rangle^2\le0$,
\begin{align}
    \big\langle\nabla_1\hat f(w,\bv^{\otimes2}),\omega-w\big\rangle
    &=\big\langle\bfSigma_\omega(\bmmu)-\bfSigma,\bv^{\otimes2}\big\rangle
    -\big\langle\bfSigma_w-\hat\bfSigma,\bv^{\otimes2}\big\rangle
    +\big\langle\bfSigma-\hat\bfSigma,\bv^{\otimes2}\big\rangle\\
    &~~~-\langle\bmmu_\omega-\bmmu,\bv\rangle^2
    +\langle\bmmu_w-\bmmu_\omega,\bv\rangle^2\\
    &\le F(\omega)-\hat F(w)+\lambda_1(\bfSigma-\hat\bfSigma)+\|\bmmu_w-\bmmu_\omega\|_2^2.
    \label{prop:weakly:convex:oracle:eq2}
\end{align}
Combining~\eqref{prop:weakly:convex:oracle:eq1}
and~\eqref{prop:weakly:convex:oracle:eq2},
\begin{align}
    \big\langle\nabla_1\hat f(w,\bv^{\otimes2}),\omega-w\big\rangle
    \le F(\omega)-F(w)+\|\bmmu_w-\bmmu_\omega\|_2^2+\|\bmmu_w-\bmmu\|_2^2+\tau .
\end{align}
Since the right-hand side does not depend on the maximizer $\bv$, the same bound
holds for every convex combination of the vectors
$\nabla_1\hat f(w,\bv^{\otimes2})$, hence for every $g\in\partial\hat F(w)$.

Now set $\omega:=w^*$. Then $\|\bmmu_w-\bmmu_{w^*}\|_2^2=2D_{\mbX}(w^*\|w)$ and, by the triangle inequality and $\|\bmmu_{w^*}-\bmmu\|_2\le r_{1,*}$, $\|\bmmu_w-\bmmu\|_2\le\sqrt{2D_{\mbX}(w^*\|w)}+r_{1,*}$. Hence
\begin{align}
    F(w^*)\ge F(w)+\langle g,w^*-w\rangle-2D_{\mbX}(w^*\|w)
    -\Big(\sqrt{2D_{\mbX}(w^*\|w)}+r_{1,*}\Big)^2-\tau.
\end{align}
Finally, $2D_{\mbX}(w^*\|w)+\big(\sqrt{2D_{\mbX}(w^*\|w)}+r_{1,*}\big)^2
\le2\big(\sqrt{2D_{\mbX}(w^*\|w)}+r_{1,*}\big)^2$, so $\partial\hat F$ is a $(4,r_{1,*},\tau)$-weakly star-convex inexact oracle for $F$ at $w^*$  relative to $h_{\mbX}$.
\end{proof}

Next, we prove that $\partial\hat F$ is a  Lipschitz continuous oracle relatively to all three divergences considered in this paper; see Definition \ref{def:inexact:oracle}(a) in Section \ref{ss:relative:regularity}.
\begin{lemma}[Relative Lipschitz continuous oracle]\label{lemma:lipschitz}
Recall that $M_{\mbX}:=\max_{i\in[n]}\|\tilde\bx_i\|_2$. For all $w,w'\in\Delta_n$ and $g\in\partial\hat F(w)$,
\begin{align}
  \langle g,w-w'\rangle
  \le\frac{M_{\mbX}^2}{2}\|w-w'\|_1+2M_{\mbX}\|\bmmu_w-\bmmu_{w'}\|_2
  \le\frac52M_{\mbX}^2\|w-w'\|_1.
  \label{lemma:lipschitz:eq1}
\end{align}
Consequently, $\partial\hat F$ is a Lipschitz continuous inexact oracle on $\Delta_n$ relative to $h$, that is,
\begin{align}
  \langle g,w-w'\rangle\le M_h\sqrt{2D_h(w'\|w)},
\end{align}
with constant
$M_{h_{\KL}}=\frac52M_{\mbX}^2$, $M_{h_{\E}}=\frac52\sqrt n\,M_{\mbX}^2$, $M_{h_{\mbX,\alpha}}=M_{\mbX}\sqrt{4+\frac{nM_{\mbX}^2}{4\alpha}}$ for $\alpha>0$, where $w\in\ri(\Delta_n)$ for $h=h_{\KL}$.
\end{lemma}
\begin{proof}
The second and third inequalities are consequences of the first, Pinsker's inequality and the fact that $\|\bmmu_w-\bmmu_{w'}\|_2\le M_{\mbX}
\|w-w'\|_1$ for any $w,w'\in\Delta_n$. 

Next, we prove the first. Fix $w,w'\in\Delta_n$ and let $\delta:=w-w'$, so that $\sum_i\delta_i=0$ and  $\sum_i\delta_i\tilde\bx_i=\mbX^\top\delta=\bmmu_w-\bmmu_{w'}$. Let
$\bv\in\mbS^{d-1}$. By the gradient formula in the proof of
Lemma~\ref{lemma:marginal:directional:derivatives} with $\bfM=\bv^{\otimes2}$,
$[\nabla_1\hat f(w,\bv^{\otimes2})]_i=\langle\bv,\tilde\bx_i\rangle^2
-2\langle\bv,\bmmu_w\rangle\langle\bv,\tilde\bx_i\rangle$. Hence
\begin{align}
  \big\langle\nabla_1\hat f(w,\bv^{\otimes2}),\delta\big\rangle
  =\sum_{i=1}^n\delta_i\langle\bv,\tilde\bx_i\rangle^2
  -2\langle\bv,\bmmu_w\rangle\langle\bv,\bmmu_w-\bmmu_{w'}\rangle .
  \label{lemma:lipschitz:eq2'}
\end{align}

Write $\delta=\delta^+-\delta^-$ with $\delta^\pm\ge0$ coordinatewise. Since
$\sum_i\delta_i=0$, we have $\sum_i\delta_i^+=\sum_i\delta_i^-=\frac12\|\delta\|_1$.
The values $a_i:=\langle\bv,\tilde\bx_i\rangle^2$ lie in $[0,M_{\mbX}^2]$, so
\begin{align}
  \sum_{i=1}^n\delta_ia_i
  =\sum_{i=1}^n\delta_i^+a_i-\sum_{i=1}^n\delta_i^-a_i
  \le\frac{\|\delta\|_1}{2}\Big(\max_ia_i-\min_ia_i\Big)
  \le\frac{M_{\mbX}^2}{2}\|\delta\|_1 .
  \label{lemma:lipschitz:eq3}
\end{align}

Since $\bmmu_w$ is a convex combination of the $\tilde\bx_i$,
$|\langle\bv,\bmmu_w\rangle|\le\|\bmmu_w\|_2\le M_{\mbX}$, and by Cauchy--Schwarz,
\begin{align}
  -2\langle\bv,\bmmu_w\rangle\langle\bv,\bmmu_w-\bmmu_{w'}\rangle
  \le2M_{\mbX}\|\bmmu_w-\bmmu_{w'}\|_2 .
  \label{lemma:lipschitz:eq4}
\end{align}

Combining~\eqref{lemma:lipschitz:eq2'}, \eqref{lemma:lipschitz:eq3}
and~\eqref{lemma:lipschitz:eq4}, we obtain that~\eqref{lemma:lipschitz:eq1} holds for $g=\nabla_1\hat f(w,\bv^{\otimes2})$, for every $\bv\in\mbS^{d-1}$. By Danskin's theorem, every $g\in\partial\hat F(w)$ is a convex combination of such vectors, with $\bv$ a top unit eigenvector of $\bfSigma_w-\hat\bfSigma$. Since the right-hand side of~\eqref{lemma:lipschitz:eq1} does not depend on $\bv$, the bound is preserved under convex combinations.
\end{proof}

\section{Mirror descent on the capped simplex}\label{s:mirror:descent}

We keep the setting and notation of Section~\ref{s:landscape}, recalling the notation in its introduction, \eqref{def:w:*:r:1}, \eqref{def:r:1:*}, \eqref{def:R:alpha:r:1:alpha} and \eqref{eq:tau}.

Throughout this section, $h$ is a reference function on $\Delta_n$ (among $h_{\KL}$, $h_{\E}$ or $h_{\mbX,\alpha}$ for $\alpha>0$) (Section~\ref{s:bregman:geometry}), and we analyze the mirror descent iterates~\eqref{eq:md:update} under the properties established in Section~\ref{s:landscape}: inexact quadratic growth of $F$ (Proposition~\ref{prop:IQG}) and an inexact weakly star-convex oracle $\hat F$ (Proposition~\ref{prop:weakly:convex:oracle}), both relative to the mean-induced divergence, and Lipschitz continuity of $\hat F$ relative to all three divergences (Lemma~\ref{lemma:lipschitz}). Unless otherwise stated, the initialization is arbitrary in $\Delta_{n,k}\cap\interior(\dom h)$, and the iterates remain in this set.

\subsection{A Lyapunov dissipation function}

Propositions~\ref{prop:IQG} and~\ref{prop:weakly:convex:oracle} combine into a single lower bound on the progress of the oracle. For $w\in\Delta_{n,k}$ and $g\in\partial\hat F(w)$, write $D:=D_{\mbX,\alpha}(w^*\|w)$. Since $D_{\mbX}\le D_{\mbX,\alpha}$,
\begin{align}
  \langle g,w-w^*\rangle
  &\ge F(w)-F^*-2(\sqrt{2D_{\mbX}(w^*\|w)}+r_{1,*})^2-\tau\\
  &\ge\frac{n-k}{k}(\sqrt{2D}-r_{1,\alpha})_+^2
  -2(\sqrt{2D}+r_{1,*})^2-\tau.
\end{align}
This motivates the following definition.

\begin{definition}[Lyapunov dissipation function]\label{def:lyapunov:Phi:m:r}
Let $m>2$ and $r,\zeta\ge0$ with $r+\zeta>0$. Let 
$
    \Phi_{m,r,\zeta}(D):=m(\sqrt{2D}-r)_+^2-2(\sqrt{2D}+\zeta)^2,
$
$D\ge0$.
\end{definition}

Let $m>2$, $r,\zeta,\alpha\ge0$ with $r+\zeta>0$ and define
\begin{align}
  c_m &:=2\sqrt2\,(\sqrt m-\sqrt2),
  &&{b}_m: =\sqrt m\,(\sqrt m-\sqrt2),\label{eq:md:constants}\\
  D_{*,m}(r,\zeta) &:=\frac12\Big(\frac{\sqrt m\,r+\sqrt2\,\zeta}{\sqrt m-\sqrt2}\Big)^2, 
  &&D_*:=D_{*,\frac{n-k}{k}}(r_{1,\alpha},r_{1,*}).
  \label{eq:D:*}
\end{align}

The following lemma, proved in Appendix~\ref{app:Phi:G}, gives useful lower bounds on $\Phi_{m,r,\zeta}$ around its positive zero.
\begin{lemma}[Lower bounds on the dissipation]\label{lemma:Phi:G}
$\Phi_{m,r,\zeta}$ is convex on $[0,\infty)$, $\Phi_{m,r,\zeta}(D_{*,m}(r,\zeta))=0$, and, for all $D\ge0$,
\begin{enumerate}
  \item[\rm (i)] $\Phi_{m,r,\zeta}(D)\ge c_m
  \left(1+\frac{\zeta}{\sqrt{2D_{*,m}(r,\zeta)}}\right)\left(D-D_{*,m}(r,\zeta)\right)$;
  \item[\rm (ii)] $\Phi_{m,r,\zeta}(D)\ge b_m\left(\sqrt{2D}-\sqrt{2D_{*,m}(r,\zeta)}\right)_+^2 -2\left(\sqrt{2D_{*,m}(r,\zeta)}+\zeta\right)^2$. 
\end{enumerate}
\end{lemma}
The tangent-line bound~(i) drives the analysis of mirror descent in this section, and the  quadratic-growth bound~(ii) that of the multiplicative weights update in Section~\ref{s:mwu}.

\begin{remark}[Trimming parameter $k$ and the threshold $1/3$]
Since $\Phi_{m,r,\zeta}(D)=2(m-2)D+O(\sqrt D)$ as $D\to\infty$, the dissipation is eventually positive if and only if $m>2$, that is, $\varepsilon_{n,k}<\frac13$ for $m=\frac{n-k}{k}$. This matches the threshold $1/3$ of Zhu, Jiao and Steinhardt~\cite{2022Zhu:Jiao:Steinhardt}.
\end{remark}

\subsection{Convergence theory for any Bregman divergence}

In this subsection we analyze mirror descent for all three reference functions on $\Delta_n$ (among $h_{\KL}$, $h_{\E}$ or $h_{\mbX,\alpha}$ for $\alpha>0$). We start with two auxiliary lemmas. 
\begin{lemma}[General recursion]\label{lemma:recursion:nonsmooth}
Let $M>0$ be the Lipschitz constant of Lemma~\ref{lemma:lipschitz} for the chosen reference function $h$ (among $h_{\KL}$, $h_{\E}$ or $h_{\mbX,\alpha}$ for $\alpha>0$). For every $t\ge1$,
\begin{align}
    D_h(w^*\|w^{t+1}) \le D_h(w^*\|w^t)
    -\gamma_t\,\Phi_{\frac{n-k}{k},r_{1,\alpha},r_{1,*}}\left(D_{\mbX,\alpha}(w^*\|w^t)\right)
    +\gamma_t\tau+\frac{\gamma_t^2M^2}{2}.
    \label{lemma:recursion:nonsmooth:eq}
\end{align}
\end{lemma}

\begin{proof}
By the three-point inequality for the Bregman proximal
step~\eqref{eq:md:update}~\cite[Lemma~4.1]{2003Beck:Teboulle}, applied with the
linear function $w\mapsto\langle g^t,w\rangle$ and the point $w^*$,
\begin{align}
  D_h(w^*\|w^{t+1})
  \le D_h(w^*\|w^t)+\gamma_t\langle g^t,w^*-w^{t+1}\rangle-D_h(w^{t+1}\|w^t) .
  \label{lemma:recursion:nonsmooth:eq1}
\end{align}
We split
$\langle g^t,w^*-w^{t+1}\rangle
=\langle g^t,w^*-w^t\rangle+\langle g^t,w^t-w^{t+1}\rangle$.
By the weak star-convexity of the oracle,
\begin{align}
  \langle g^t,w^*-w^t\rangle
  \le F(w^*)-F(w^t)+2\Big(\sqrt{2D_{\mbX}(w^*\|w^t)}+r_{1,*}\Big)^2 + \tau.
  \label{lemma:recursion:nonsmooth:eq2}
\end{align}
By the $\IQG$ property, $D_{\mbX}\le D_{\mbX,\alpha}$ and the monotonicity of
$D\mapsto(\sqrt{2D}+r_{1,*})^2$,
\begin{align}
    \langle g^t,w^*-w^t\rangle
  &\le-m\Big(\sqrt{2D_{\mbX,\alpha}(w^*\|w^t)}-r_{1,\alpha}\Big)_+^2
  +2\Big(\sqrt{2D_{\mbX,\alpha}(w^*\|w^t)}+r_{1,*}\Big)^2+\tau\\
  &=-\Phi_{m,r_{1,\alpha},r_{1,*}}\big(D_{\mbX,\alpha}(w^*\|w^t)\big)+\tau.
  \label{lemma:recursion:nonsmooth:eq3}
\end{align}
By the Lipschitz continuity of the oracle relative to $h$ and the inequality
$\sqrt2\,ay-y^2\le\frac{a^2}{2}$ for $a,y\ge0$,
\begin{align}
  \gamma_t\langle g^t,w^t-w^{t+1}\rangle-D_h(w^{t+1}\|w^t)
  &\le\gamma_tM\sqrt{2D_h(w^{t+1}\|w^t)}-D_h(w^{t+1}\|w^t) \\
  &\le\frac{\gamma_t^2M^2}{2}.
  \label{lemma:recursion:nonsmooth:eq4}
\end{align}
Inserting~\eqref{lemma:recursion:nonsmooth:eq3}
and~\eqref{lemma:recursion:nonsmooth:eq4}
into~\eqref{lemma:recursion:nonsmooth:eq1}
gives~\eqref{lemma:recursion:nonsmooth:eq}.
\end{proof}

The next lemma, proved in Appendix~\ref{app:recursion:slow}, analyzes this recursion with a constant step small size.

\begin{lemma}[Constant small step size]\label{lemma:recursion:slow}
Let $m>2$ and $r,\zeta\ge0$ with $r+\zeta>0$. Let $(D_t)_{t\ge1}$ and $(\bar D_t)_{t\ge1}$ be nonnegative sequences and $\gamma_t>0$
such that, for all $t\ge1$,
\begin{align}
    D_{t+1}\le D_t-\gamma_t\,\Phi_{m,r,\zeta}(\bar D_t)+\gamma_t\tau+\frac{\gamma_t^2M^2}{2}.
    \label{lemma:recursion:slow:eq0}
\end{align}
Let $T\ge1$ and $\Omega\ge D_1$, and take $\gamma_t:=\frac1M\sqrt{2\Omega/T}$ for
$t\in[T]$. Then
\begin{align}
  \frac1T\sum_{t=1}^T\bar D_t
  \le D_*+\frac{\tau}{c_m}+\frac{M}{2(\sqrt m-\sqrt2)}\sqrt{\frac{\Omega}{T}}.
    \label{lemma:recursion:slow:eq}
\end{align}
\end{lemma}

Combining Lemmas~\ref{lemma:recursion:nonsmooth} and~\ref{lemma:recursion:slow} gives ergodic rates for the distance to the good weights and for the optimality gap.
\begin{proposition}[General Bregman divergence: ergodic rates]\label{prop:ergodic:rates}
Assume that $k<n/3$, so that $m:=\frac{n-k}{k}>2$, and the conditions of Lemma~\ref{lemma:recursion:nonsmooth} hold. Let $\Omega\ge D_h(w^*\|w^1)$, and run~\eqref{eq:md:update} with step size
$$
  \gamma_t:=\frac1M\sqrt{2\Omega/T},~~(t\in[T]).
$$
Define the \emph{statistical} and \emph{optimization errors}, 
\begin{align}
  \mcE_{\rm stat}(\alpha) &:= \frac{\sqrt m\,r_{1,\alpha}+\sqrt2\,r_{1,*}}{\sqrt m-\sqrt2}
  +\sqrt{\frac{2\tau}{c_m}}+r_{1,*},\label{prop:ergodic:rates:stat:error}\\
  \mcE_{\rm opt}(T) &:= \Big(\frac{M}{\sqrt m-\sqrt2}\Big)^{1/2}
  \Big(\frac{\Omega}{T}\Big)^{1/4}.
  \label{prop:ergodic:rates:opt:error}
\end{align}
Then the following hold.
\begin{enumerate}
\item[\rm (i)] \textbf{Ergodic average:} Let $\bar w^T:=\frac1T\sum_{t=1}^Tw^t\in\Delta_{n,k}$. Then
\begin{align}
    F(\bar w^T)-F^* \le\frac1T\sum_{t=1}^T\big(F(w^t)-F^*\big)
    &\le\underbrace{4\,\mcE_{\rm stat}(\alpha)^2+\tau}_{\text{statistical}}
    +\underbrace{\Big(\sqrt2+\frac{4}{\sqrt m-\sqrt2}\Big)M\sqrt{\frac{\Omega}{T}}}_{\text{optimization}}.
    \label{prop:ergodic:rates:gap}
\end{align}
\item[\rm (ii)] \textbf{Distance to mean:}
\begin{align}
  \|\bmmu_{\bar w^T}-\bmmu\|_2 &\le \underbrace{r_{1,*}+2\sqrt{\frac{\varepsilon_{n,k}}{1-\varepsilon_{n,k}}}\,\mcE_{\rm stat}(\alpha)
  +\sqrt{\frac{\varepsilon_{n,k}\,\tau}{1-\varepsilon_{n,k}}}}_{\text{statistical}}\\
  &+\underbrace{\sqrt{\frac{\varepsilon_{n,k}}{1-\varepsilon_{n,k}}}
  \Big(\sqrt2+\frac{4}{\sqrt m-\sqrt2}\Big)^{\frac12}M^{\frac12}\Big(\frac{\Omega}{T}\Big)^{\frac14}}_{\text{optimization}}.\label{prop:ergodic:rates:distance:mean}
\end{align}
\end{enumerate}
In particular, with $w^1=\unit_{[n]}/n$, one can take:
\begin{enumerate}
  \item[\rm (a)] for $h=h_{\KL}$ (entropic mirror descent),
  $M=\frac52M_{\mbX}^2$ and $\Omega=\log\frac{1}{1-\varepsilon_{n,k}}$;
  \item[\rm (b)] for $h=h_{\E}$ (projected subgradient method),
  $M=\frac52M_{\mbX}^2\sqrt n$ and $\Omega=\frac{k}{2n(n-k)}$.
\end{enumerate}
\end{proposition}
\begin{proof}
Let $D_t:=D_h(w^*\|w^t)$ and $\bar D_t:=D_{\mbX,\alpha}(w^*\|w^t)$. By
Lemma~\ref{lemma:recursion:nonsmooth}, these sequences
satisfy~\eqref{lemma:recursion:slow:eq0} with $r:=r_{1,\alpha}$ and $\zeta:=r_{1,*}$, and Lemma~\ref{lemma:recursion:slow} gives
$\frac1T\sum_{t=1}^T\bar D_t\le D_*+\frac{\tau}{c_m}+\frac12\mcE_{\rm opt}(T)^2$. Since $\sqrt{2D_*}=\frac{\sqrt m\,r_{1,\alpha}+\sqrt2\,r_{1,*}}{\sqrt m-\sqrt2}$, we have $\sqrt{2D_*}+\sqrt{2\tau/c_m}=\mcE_{\rm stat}(\alpha)-r_{1,*}$, and expanding the square gives $D_*+\frac{\tau}{c_m}+\frac12\mcE_{\rm opt}(T)^2\le\frac12\big(\mcE_{\rm stat}(\alpha)-r_{1,*}+\mcE_{\rm opt}(T)\big)^2$. Hence, 
\begin{align}
  D_{\mbX,\alpha}(w^*\|\bar w^T)
    \le\frac1T\sum_{t=1}^TD_{\mbX,\alpha}(w^*\|w^t)
    &\le\frac12\big(\mcE_{\rm stat}(\alpha)+\mcE_{\rm opt}(T)\big)^2, \label{prop:ergodic:rates:dist}
\end{align}
where we used Jensen's inequality. 

\emph{Proof of (i).} Combining~\eqref{lemma:recursion:nonsmooth:eq1},
\eqref{lemma:recursion:nonsmooth:eq2} and~\eqref{lemma:recursion:nonsmooth:eq4},
and using $D_{\mbX}\le D_{\mbX,\alpha}$,
\begin{align}
    \gamma_t\big(F(w^t)-F^*\big)
    \le D_t-D_{t+1}+2\gamma_t\Big(\sqrt{2\bar D_t}+r_{1,*}\Big)^2
    +\gamma_t\tau+\frac{\gamma_t^2M^2}{2} .
\end{align}
Summing over $t\in[T]$, using $D_{T+1}\ge0$ and $D_1\le\Omega$, and dividing by
$\gamma T$ with $\gamma:=\frac1M\sqrt{2\Omega/T}$,
\begin{align}
    \frac1T\sum_{t=1}^T\big(F(w^t)-F^*\big)
    & \le\frac{\Omega}{\gamma T}+\frac{\gamma M^2}{2}+\tau
    +\frac2T\sum_{t=1}^T\Big(\sqrt{2\bar D_t}+r_{1,*}\Big)^2 \\
    &= M\sqrt{\frac{2\Omega}{T}}+\tau
    +\frac2T\sum_{t=1}^T\Big(\sqrt{2\bar D_t}+r_{1,*}\Big)^2.
\end{align}
The function $D\mapsto(\sqrt{2D}+r_{1,*})^2$ is concave and nondecreasing on $[0,\infty)$, so, by Jensen's inequality and~\eqref{prop:ergodic:rates:dist},
$\frac1T\sum_t(\sqrt{2\bar D_t}+r_{1,*})^2\le\big(\mcE_{\rm stat}(\alpha)+\mcE_{\rm opt}(T)\big)^2\le2\mcE_{\rm stat}(\alpha)^2+2\mcE_{\rm opt}(T)^2$. Since $M\sqrt{2\Omega/T}+4\mcE_{\rm opt}(T)^2=\big(\sqrt2+\frac{4}{\sqrt m-\sqrt2}\big)M\sqrt{\Omega/T}$, this gives the second inequality in~\eqref{prop:ergodic:rates:gap}. The first follows from the convexity of $F$.

\emph{Proof of (ii).} From \eqref{eq:rates:signature}, 
\[
  \|\bmmu_{\bar w^T}-\bmmu\|_2
  \le r_{1,*}+\sqrt{\frac{\varepsilon_{n,k}}{1-\varepsilon_{n,k}}(F(\bar w^T)-F^*)},
\]
and~\eqref{prop:ergodic:rates:gap} bounds $F(\bar w^T)-F^*$. This proves item (ii).

We now prove items (a) and (b). For~(a), the iterates remain in $\ri(\Delta_n)$, Lemma~\ref{lemma:lipschitz} gives $M=\frac52M_{\mbX}^2$, and Lemma~\ref{lemma:capped:kl} gives $\sup_{w\in\Delta_{n,k}}D_{\KL}(w\|\unit_{[n]}/n)=\log\frac{1}{1-\varepsilon_{n,k}}$. For~(b), $\dom(h_{\E})=\re^n$, so the iterates are the Euclidean projections onto $\Delta_{n,k}$ and are well defined. The inequality in~\eqref{lemma:lipschitz:eq1}, which holds for all $w,w'\in\Delta_n$, and $\|w-w'\|_1\le\sqrt n\,\|w-w'\|_2=\sqrt n\sqrt{2D_{\E}(w'\|w)}$ give $M=\frac52\sqrt n\,M_{\mbX}^2$. Finally, Lemma~\ref{lemma:capped:diameter} gives $\sup_{w\in\Delta_{n,k}}D_{\E}(w\|\unit_{[n]}/n) =\frac12\cdot\frac{k}{n(n-k)}$.
\end{proof}

\subsection{Convergence theory for the mean-induced divergence}

We now specialize to $h=h_{\mbX,\alpha}$, so that mirror descent runs in the same divergence as the quadratic growth of $F$: the mean-induced divergence. Then $D_h(w^*\|w^t)=D_{\mbX,\alpha}(w^*\|w^t)$, and the recursion of Lemma~\ref{lemma:recursion:nonsmooth} involves a single sequence. We will need the following two auxiliary lemmas, proved in Appendices~\ref{app:recursion:G:plus} and~\ref{app:series:bounds}.

\begin{lemma}[Decreasing step sizes]\label{lemma:recursion:G:plus}
Let $m>2$ and $r,\zeta\ge0$ with $r+\zeta>0$. Let $K\ge\sqrt{2m}-1$ and step size
$
  \gamma_t:=\frac{2}{c_m(t+K)}
$
for $t\in\mathbb{N}$. Let $(D_t)_{t\ge1}$ be a nonnegative sequence such that, for all $t\in\mbN$,
\begin{align}
    D_{t+1}\le D_t-\gamma_t\,\Phi_{m,r,\zeta}(D_t)+\gamma_t\tau+\frac{\gamma_t^2M^2}{2}.
    \label{lemma:recursion:G:plus:eq0}
\end{align}
Then, for all $t\in\mbN$,
\begin{align}
  \left(D_t - D_{*,m}^\tau(r,\zeta)\right)_+
    \le\Big(\frac{1+K}{t+K}\Big)^2(D_1-D_{*,m}^\tau(r,\zeta))_+ +\frac{\frac{9M^2}{2c_m^2}}{t+K}, 
\end{align}
where $D_{*,m}^\tau(r,\zeta):=D_{*,m}(r,\zeta)+\frac{\tau}{c_m}$.
\end{lemma}

\begin{lemma}\label{lemma:series:bounds}
Let $1\le S<T$ be integers and $K\ge0$. Then
$\sum_{t=S+1}^T\frac{1}{(t+K)^2}\le\frac{1}{S+K}$ and
$\sum_{t=S+1}^T\frac{1}{t+K}\le\log\frac{T+K}{S+K}$. If, moreover,
$S=(1-\theta)T$ with $\theta\in(0,1)$, then
$\log\frac{T+K}{S+K}\le\log\frac{1}{1-\theta}$.
\end{lemma}

Combining Lemma~\ref{lemma:recursion:nonsmooth} with Lemmas~\ref{lemma:recursion:G:plus} and~\ref{lemma:series:bounds} gives, for decreasing step sizes, last-iterate and windowed rates for the distance to the good weights and for the optimality gap.

\begin{proposition}[Mean-induced divergence: last-iterate and windowed ergodic rates]\label{prop:window:rates}
Assume that $k<n/3$, so that $m:=\frac{n-k}{k}>2$, and the conditions of Lemma~\ref{lemma:recursion:nonsmooth} hold with $h:=h_{\mbX,\alpha}$ for some $\alpha>0$, and $M:=M_{\mbX}\sqrt{4+\nicefrac{nM_{\mbX}^2}{4\alpha}}$. Let $K\ge\sqrt{2m}-1$ and $\Omega\ge D_{\mbX,\alpha}(w^*\|w^1)$, and run~\eqref{eq:md:update} with step size
\begin{align}
  \gamma_t:=\frac{2}{c_m(t+K)},~~~(t\in\mathbb{N}).
\end{align}
Let $\theta\in(0,1)$ and $T\in\mbN$ be such that $S:=(1-\theta)T\in\mbN$, and let
$\bar w^T:=\frac{1}{\theta T}\sum_{t=S+1}^Tw^t\in\Delta_{n,k}$. With $\mcE_{\rm stat}(\alpha)$ as in~\eqref{prop:ergodic:rates:stat:error}, define the \emph{optimization errors}
\begin{align}
  \mcE_{\rm opt}^{\rm last}(t)
  &:=\bigg[2\Big(\frac{1+K}{t+K}\Big)^2\Omega+\frac{\frac{9M^2}{c_m^2}}{t+K}\bigg]^{1/2},
  \label{prop:window:rates:errors:last:opt}\\
  \mcE_{\rm opt}^{\rm win}(T)
  &:=\bigg[\frac{2}{\theta T}\Big(\frac{(1+K)^2\,\Omega}{S+K}
  +\frac{9M^2}{2c_m^2}\log\frac{1}{1-\theta}\Big)\bigg]^{1/2},
  \label{prop:window:rates:errors:win:opt}\\
  \mcE_{\rm opt}^{\Sigma}(T)
  &:=\frac{c_m}{4}\Big(\frac{(1+K)\,\mcE_{\rm stat}(\alpha)^2
  +(S+1+K)\,\mcE_{\rm opt}^{\rm last}(S+1)^2}{\theta T}
  +\mcE_{\rm opt}^{\rm win}(T)^2\Big)\\
  &+4\,\mcE_{\rm opt}^{\rm win}(T)^2
  +\frac{M^2}{c_m\theta T}\log\frac{1}{1-\theta}.
  \label{prop:window:rates:errors:sigma:opt}
\end{align}
Then the following hold.
\begin{enumerate}
\item[\rm (i)] \textbf{Last iterate:} For all $t\in\mbN$,
$D_{\mbX,\alpha}(w^*\|w^t)\le\frac12\big(\mcE_{\rm stat}(\alpha)-r_{1,*}+\mcE_{\rm opt}^{\rm last}(t)\big)^2$ and 
$\|\bmmu_{w^t}-\bmmu\|_2\le\mcE_{\rm stat}(\alpha)+\mcE_{\rm opt}^{\rm last}(t)$. 
\item[\rm (ii)] \textbf{Windowed average:}
\begin{align}
    F(\bar w^T)-F^*\le\underbrace{\Big(4+\frac{c_m}{4\theta}\Big)\mcE_{\rm stat}(\alpha)^2+\tau}_{\text{statistical}}
    +\mcE_{\rm opt}^{\Sigma}(T).
    \label{prop:window:rates:gap}
\end{align}
\item[\rm (iii)] \textbf{Distance to mean:} 
\begin{align}
  \|\bmmu_{\bar w^T}-\bmmu\|_2 \le \underbrace{r_{1,*}+\Big(\frac{4}{m}+\frac{c_m}{4\theta m}\Big)^{1/2}\mcE_{\rm stat}(\alpha)+\sqrt{\frac{\varepsilon_{n,k}}{1-\varepsilon_{n,k}}\tau}}_{\text{statistical}}
  +\underbrace{\sqrt{\frac{\varepsilon_{n,k}}{1-\varepsilon_{n,k}}\mcE_{\rm opt}^{\Sigma}(T)}}_{\text{optimization}}.
  \label{prop:window:rates:distance:mean}
\end{align}
\end{enumerate}
In particular, for $w^1:=\unit_{[n]}/n$, we can take 
$
  \Omega:=2M_{\mbX}^2\varepsilon_{n,k}^2+\frac{\alpha k}{2n(n-k)}.
$
\end{proposition}
\begin{proof}
Let $D_t:=D_{\mbX,\alpha}(w^*\|w^t)$, $D_*^\tau:=D_*+\frac{\tau}{c_m}$ and $A:=\frac{9M^2}{c_m^2}$. Since $D_h(w^*\|w^t)=D_t$,
Lemma~\ref{lemma:recursion:nonsmooth} gives
\begin{align}
  D_{t+1}\le D_t-\gamma_t\Phi_{m,r_{1,\alpha},r_{1,*}}(D_t)+\gamma_t\tau+\frac{\gamma_t^2M^2}{2}.\label{prop:window:rates:distance:eq1}
\end{align}
Since $\sqrt{2D_*}=\frac{\sqrt m\,r_{1,\alpha}+\sqrt2\,r_{1,*}}{\sqrt m-\sqrt2}$, we have $\sqrt{2D_*^\tau}\le\mcE_{\rm stat}(\alpha)-r_{1,*}$, so that, expanding the square, $D_*^\tau+\frac12c^2\le\frac12\big(\mcE_{\rm stat}(\alpha)-r_{1,*}+c\big)^2$ for every $c\ge0$.

\emph{Proof of (i).} Lemma~\ref{lemma:recursion:G:plus} applies with
$r:=r_{1,\alpha}$ and $\zeta:=r_{1,*}$, and, with $(D_1-D_*^\tau)_+\le\Omega$, gives 
\begin{align}
  D_t\le D_*^\tau+\left(\frac{1+K}{t+K}\right)^2\Omega+\frac{A}{2(t+K)}=D_*^\tau+\frac12\mcE_{\rm opt}^{\rm last}(t)^2,
\end{align}
which gives~(i) by the previous paragraph with $c=\mcE_{\rm opt}^{\rm last}(t)$.

We now give the bound on the second inequality in (i). From the bound \eqref{def:r:1:*} on $\|\bmmu_{w^*}-\bmmu\|_2$ and the triangle inequality,
\[
  \|\bmmu_w-\bmmu\|_2\le\|\bmmu_w-\bmmu_{w^*}\|_2+\|\bmmu_{w^*}-\bmmu\|_2
  \le r_{1,*}+\sqrt{2D_{\mbX,\alpha}(w^*\|w)}.
\]
Applying it to $w:=w^t$ and using the already proof on the bound for $D_{\mbX,\alpha}(w^*\|w^t)$ we prove the bound on the second inequality in (i). 

\emph{Proof of (ii), optimality gap.} Since $K\ge\sqrt{2m}-1>0$,
Lemma~\ref{lemma:series:bounds} gives
$\sum_{t=S+1}^T(t+K)^{-2}\le\frac{1}{S+K}$ and
$\sum_{t=S+1}^T(t+K)^{-1}\le\log\frac{1}{1-\theta}$, so the bound in the proof of~(i) gives
\begin{align}
  \frac{1}{\theta T}\sum_{t=S+1}^TD_t\le D_*^\tau+\frac12\mcE_{\rm opt}^{\rm win}(T)^2\le\frac12\big(\mcE_{\rm stat}(\alpha)-r_{1,*}+\mcE_{\rm opt}^{\rm win}(T)\big)^2.\label{prop:window:rates:distance:eq2}
\end{align}

Combining~\eqref{lemma:recursion:nonsmooth:eq1},
\eqref{lemma:recursion:nonsmooth:eq2} and~\eqref{lemma:recursion:nonsmooth:eq4}
and dividing by $\gamma_t$,
\begin{align}
    F(w^t)-F^*\le\frac{D_t-D_{t+1}}{\gamma_t}
    +2\big(\sqrt{2D_t}+r_{1,*}\big)^2 + \tau + \frac{\gamma_tM^2}{2}.
    \label{prop:window:rates:distance:eq3}
\end{align}
Sum over $t=S+1,\ldots,T$ and divide by $\theta T$. Since
$\gamma_t^{-1}-\gamma_{t-1}^{-1}=\frac{c_m}{2}$, summation by parts gives
$$
  \sum_{t=S+1}^T\frac{D_t-D_{t+1}}{\gamma_t}
  \le\frac{c_m}{2}\Big((S+1+K)D_{S+1}+\sum_{t=S+1}^TD_t\Big).
$$

By \eqref{prop:window:rates:distance:eq1}, $D_{S+1}\le D_*^\tau+\frac12\mcE_{\rm opt}^{\rm last}(S+1)^2$, and, by \eqref{prop:window:rates:distance:eq2}, $\frac1{\theta T}\sum_{t=S+1}^TD_t\le D_*^\tau+\frac12\mcE_{\rm opt}^{\rm win}(T)^2$. Writing $\frac{S+1+K}{\theta T}=\frac{1-\theta}{\theta}+\frac{1+K}{\theta T}$ and using $D_*^\tau\le\frac12\mcE_{\rm stat}(\alpha)^2$, the average 
$\frac{1}{\theta T}\sum_{t=S+1}^T\frac{D_t-D_{t+1}}{\gamma_t}$ is at most
$$
  \frac{c_m}{4\theta}\mcE_{\rm stat}(\alpha)^2
  +\frac{c_m}{4}\Big(\frac{(1+K)\,\mcE_{\rm stat}(\alpha)^2+(S+1+K)\,\mcE_{\rm opt}^{\rm last}(S+1)^2}{\theta T}
  +\mcE_{\rm opt}^{\rm win}(T)^2\Big).
$$

   Since $D\mapsto(\sqrt{2D}+r_{1,*})^2$ is concave and nondecreasing, Jensen's inequality and~\eqref{prop:window:rates:distance:eq2} imply $\frac{1}{\theta T}\sum_{t=S+1}^T2\big(\sqrt{2D_t}+r_{1,*}\big)^2\le2\big(\mcE_{\rm stat}(\alpha)+\mcE_{\rm opt}^{\rm win}(T)\big)^2\le4\mcE_{\rm stat}(\alpha)^2+4\mcE_{\rm opt}^{\rm win}(T)^2$. 

Finally, $\frac{M^2}{2}\sum_{t=S+1}^T\gamma_t\le\frac{M^2}{c_m}\log\frac{1}{1-\theta}$ by Lemma~\ref{lemma:series:bounds}. Collecting terms and using the convexity of $F$ gives~\eqref{prop:window:rates:gap}.

\emph{Proof of (iii), windowed average.} From \eqref{eq:rates:signature}, 
\[
  \|\bmmu_{\bar w^T}-\bmmu\|_2
  \le r_{1,*}+\sqrt{\frac{\varepsilon_{n,k}}{1-\varepsilon_{n,k}}(F(\bar w^T)-F^*)},
\]
and~\eqref{prop:window:rates:gap} bounds $F(\bar w^T)-F^*$. This proves item (iii). 

For $w^1:=\unit_{[n]}/n$, we can take
$
  \Omega:=2M_{\mbX}^2\varepsilon_{n,k}^2+\frac{\alpha k}{2n(n-k)}.
$
Indeed, for every $w\in\Delta_{n,k}$, $\|\bmmu_w-\bmmu_{w^1}\|_2\le M_{\mbX}\|w-w^1\|_1\le2M_{\mbX}\varepsilon_{n,k}$ by Lemma~\ref{lemma:capped:tv}, and $\|w-w^1\|_2^2\le\frac{k}{n(n-k)}$ by Lemma~\ref{lemma:capped:diameter}; hence $D_{\mbX,\alpha}(w\|w^1)=\frac12\|\bmmu_w-\bmmu_{w^1}\|_2^2+\frac\alpha2\|w-w^1\|_2^2\le\Omega$, in particular for $w=w^*$.
\end{proof}

\section{Multiplicative weights update on the capped simplex}\label{s:mwu}

Again, we keep the setting and notation of Section~\ref{s:landscape}, in particular~\eqref{def:w:*:r:1}, \eqref{def:r:1:*}, \eqref{def:R:alpha:r:1:alpha} and~\eqref{eq:tau}.

In this section we analyze the MWU algorithm~\eqref{eq:mwu} stated in Section \ref{s:main:results} with score $s^t$, step sizes $\gamma_t>0$ and constant ${M_\infty}$ as defined there; in particular $\frac{\gamma_t}{{M_\infty^2}}\|s^t\|_\infty\le\frac12$ for all $t$. We will again use  properties established in Section \ref{s:landscape} relative to the mean-induced divergence: the spectral signature property for $F$ (Proposition \ref{prop:IQG}) and, indirectly, the inexact weakly star-convex oracle $\hat F$ (see inequality \eqref{prop:weakly:convex:oracle:eq1} in the proof of Proposition~\ref{prop:weakly:convex:oracle}). Unless otherwise stated, the initialization is arbitrary in $\Delta_{n,k}\cap\ri(\Delta_n)$, and the iterates remain in this set.

\begin{proposition}[MWU: ergodic rates]\label{prop:mwu:ergodic}
Let $k<n/3$, so that $m:=\frac{n-k}{k}>2$, let $\Omega\ge D_{\KL}(w^*\|w^1)$, $T\ge\Omega$ and run~\eqref{eq:mwu} with the constant step size 
\begin{align}
  \gamma_t\equiv\sqrt{\Omega/T},~~~(t\in[T]).
\end{align}
Define the \emph{statistical} and \emph{optimization errors}
\begin{align}
  \mcE_{\rm stat}^{\rm MW}
  &:=\frac{\sqrt m+\sqrt2}{\sqrt m-\sqrt2}\,r_{1,*}
+\sqrt{\frac{2}{b_m}}\cdot\frac{2\sqrt m}{\sqrt m-\sqrt2}\,r_{1,*} +\sqrt{\frac{\tau}{b_m}}\label{prop:mwu:ergodic:errors:stat:MWU}\\
  \mcE_{\rm opt}^{\rm MW}(T) &:=\Big(\frac{2{M_\infty^2}}{b_m}\Big)^{1/2}
  \left(\frac{\Omega}{T}\right)^{1/4}.
  \label{prop:mwu:ergodic:errors:opt:MWU}
\end{align}
Let $\bar w^T:=\frac1T\sum_{t=1}^Tw^t\in\Delta_{n,k}$. Then:
\begin{enumerate}
\item[\rm (i)] \textbf{Ergodic average:}
$
\|\bmmu_{\bar w^T}-\bmmu_{w^*}\|_2\le r_{1,*}+\mcE_{\rm stat}^{\rm MW}+\mcE_{\rm opt}^{\rm MW}(T)
$
and
\begin{align}
  F(\bar w^T)-F^* &\le\frac1T\sum_{t=1}^T\big(F(w^t)-F^*\big)
  \le\underbrace{
    \Big(1+\frac{2}{b_m}\Big)\frac{16\,m\,r_{1,*}^2}{(\sqrt m-\sqrt2)^2}
    +\Big(1+\frac{4}{b_m}\Big)\tau
  }_{\text{statistical}} \\
  &+\underbrace{2\Big(1+\frac{4}{b_m}\Big){M_\infty^2}\sqrt{\frac{\Omega}{T}}}_{\text{optimization}}.
  \label{prop:mwu:ergodic:dist:gap}
\end{align}
\item[\rm (ii)] \textbf{Distance to mean:}
$
  \|\bmmu_{\bar w^T}-\bmmu\|_2\le\mcE_{\rm stat}^{\rm MW}+\mcE_{\rm opt}^{\rm MW}(T).
$
\end{enumerate}
In particular, for $w^1=\unit_{[n]}/n$, one can take $\Omega=\log\frac{1}{1-\varepsilon_{n,k}}$.
\end{proposition}
\begin{proof}
Let $\gamma:=\sqrt{\Omega/T}$. From the assumptions on $T\ge\Omega$ and ${M_\infty^2}$, $\gamma\frac{\|s^t\|_\infty}{{M_\infty^2}}\le\frac12$ for all $t$. The regret bound of~\cite[Theorem~2.4]{2012Arora:Hazan:Kale} implies\footnote{Although~\cite{2012Arora:Hazan:Kale} state it for a constant step size, their proof is a one-step inequality that uses only $\gamma_t\frac{\|s^t\|_\infty}{{M_\infty^2}}\le\frac12$.}
\begin{align}
    \frac{1}{T}\sum_{t=1}^T\langle s^t,w^t\rangle
    \le\frac{1}{T}\sum_{t=1}^T\langle s^t,w^*\rangle
      +\frac{\gamma}{T}\sum_{t=1}^T\langle|s^t|,w^*\rangle
      +\frac{{M_\infty^2}\,D_{\KL}(w^*\|w^1)}{\gamma T}.
    \label{MWU:regret}
\end{align}
Write $x_t:=\|\bmmu_{w^t}-\bmmu\|_2$, $d_t:=\|\bmmu_{w^t}-\bmmu_{w^*}\|_2$ and
$u_t:=\sqrt{\frac{\varepsilon_{n,k}}{1-\varepsilon_{n,k}}\big(F(w^t)-F^*\big)}$.

By definition of $\bv^t$,
$\langle s^t,w^t\rangle=\langle\bv^t,(\bfSigma_{w^t}-\hat\bfSigma)\bv^t\rangle=\hat F(w^t)$,
and $\hat F(w^t)\ge F(w^t)-x_t^2-\lambda_1(\hat\bfSigma-\bfSigma)$ by~\eqref{prop:weakly:convex:oracle:eq1}.

By the bias--variance decomposition,
\begin{align}
        \frac{1}{T}\sum_{t=1}^T\langle s^t,w^*\rangle
    &=\frac{1}{T}\sum_{t=1}^T\langle\bv^t,[\bfSigma_{w^*}(\bmmu_{w^t})-\hat\bfSigma]\bv^t\rangle\\
    &=\frac{1}{T}\sum_{t=1}^T\langle\bv^t,[\bfSigma_{w^*}(\bmmu)-\bfSigma]\bv^t\rangle
    +\frac{1}{T}\sum_{t=1}^T\langle\bv^t,[\bfSigma-\hat\bfSigma]\bv^t\rangle\\
    &~~~-\frac{1}{T}\sum_{t=1}^T\langle\bv^t,\bmmu_{w^*}-\bmmu\rangle^2
    +\frac{1}{T}\sum_{t=1}^T\langle\bv^t,\bmmu_{w^t}-\bmmu_{w^*}\rangle^2\\
    &\le F(w^*)+\lambda_1(\bfSigma-\hat\bfSigma)+\frac{1}{T}\sum_{t=1}^Td_t^2.
\end{align}

Moreover, by~\eqref{eq:spectral:centrality:at:bmmu:w:*} with $w=w^t$ and the line
following it, $d_t\le2r_{1,*}+u_t$; since the first two terms on the right-hand
side of~\eqref{eq:spectral:centrality:at:bmmu:w:*} bound $x_t$ by item~(a) of
Definition~\ref{def:favorable}, the same argument gives $x_t\le r_{1,*}+u_t$.

From ${M_\infty^2}\ge\sup_{t\in\mbN}\|s^t\|_\infty$, $\|w^*\|_1=1$, $\Omega\ge D_{\KL}(w^*\|w^1)$ and $\gamma=\sqrt{\Omega/T}$,
\begin{align}
  \frac{\gamma}{T}\sum_{t=1}^T\langle|s^t|,w^*\rangle
      +\frac{{M_\infty^2}\,D_{\KL}(w^*\|w^1)}{\gamma T}
  \le {M_\infty^2}\Big(\gamma+\frac{\Omega}{\gamma T}\Big)=2{M_\infty^2}\sqrt{\frac{\Omega}{T}}=b_m\,\mcE_{\rm opt}^{\rm MW}(T)^2.
\end{align}

Collecting the previous bounds in~\eqref{MWU:regret}, the left-hand side is at least $\frac1T\sum_t\big(F(w^t)-x_t^2\big)-\lambda_1(\hat\bfSigma-\bfSigma)$, and the right-hand side is at most $F^*+\lambda_1(\bfSigma-\hat\bfSigma)+\frac1T\sum_td_t^2+b_m\mcE_{\rm opt}^{\rm MW}(T)^2$. Rearranging,
\begin{align}
  \frac1T\sum_{t=1}^T\big(F(w^t)-F^*\big)
  \le\frac1T\sum_{t=1}^T\big(x_t^2+d_t^2\big)+b_m\mcE_{\rm opt}^{\rm MW}(T)^2+\tau.
  \label{prop:mwu:ergodic:pf:key}
\end{align}

We now express both sides through $u_t$. On the left, $F(w^t)-F^*=mu_t^2$, since
$\frac{1-\varepsilon_{n,k}}{\varepsilon_{n,k}}=m$. On the right,
$x_t\le r_{1,*}+u_t$ and $d_t\le2r_{1,*}+u_t$, so
$x_t^2+d_t^2\le2(u_t+2r_{1,*})^2$. Thus~\eqref{prop:mwu:ergodic:pf:key} becomes
\begin{align}
  \frac1T\sum_{t=1}^T\phi(u_t)\le b_m\mcE_{\rm opt}^{\rm MW}(T)^2+\tau,
  ~~\text{where}~~
  \phi(u):=mu^2-2(u+2r_{1,*})^2.
  \label{prop:mwu:ergodic:pf:phi}
\end{align}
Since $\phi(u)=\Phi_{m,0,2r_{1,*}}(u^2/2)$, whose radius in
Lemma~\ref{lemma:Phi:G} is 
$
  D_*^{\rm MW} := 4(\frac{r_{1,*}}{\sqrt{\frac{n-k}{k}}-\sqrt2})^2,
$
Lemma~\ref{lemma:Phi:G}(ii) gives
$
  \phi(u)\ge b_m\big(u-\sqrt{2D_*^{\rm MW}}\big)_+^2
  -2\big(\sqrt{2D_*^{\rm MW}}+2r_{1,*}\big)^2,~~u\ge0.
$
Inserting this into~\eqref{prop:mwu:ergodic:pf:phi} gives
\begin{align}
  \frac1T\sum_{t=1}^T\big(u_t-\sqrt{2D_*^{\rm MW}}\big)_+^2
  \le\mcE_{\rm opt}^{\rm MW}(T)^2+\frac{\tau+2\big(\sqrt{2D_*^{\rm MW}}+2r_{1,*}\big)^2}{b_m}.
  \label{prop:mwu:ergodic:pf:excess}
\end{align}

By the Cauchy--Schwarz inequality, \eqref{prop:mwu:ergodic:pf:excess} and $\sqrt{a+b+c}\le\sqrt a+\sqrt b+\sqrt c$,
\begin{align}
  \frac1T\sum_{t=1}^Tu_t
  \le\sqrt{2D_*^{\rm MW}}+\frac1T\sum_{t=1}^T\big(u_t-\sqrt{2D_*^{\rm MW}}\big)_+
  \le\mcE_{\rm stat}^{\rm MW}-r_{1,*}+\mcE_{\rm opt}^{\rm MW}(T).
\end{align}

The maps $w\mapsto\|\bmmu_w-\bmmu\|_2$ and $w\mapsto\|\bmmu_w-\bmmu_{w^*}\|_2$ are convex; hence $\|\bmmu_{\bar w^T}-\bmmu\|_2\le\frac1T\sum_tx_t$ and $\|\bmmu_{\bar w^T}-\bmmu_{w^*}\|_2\le\frac1T\sum_td_t$. With $x_t\le r_{1,*}+u_t$ and $d_t\le2r_{1,*}+u_t$, this gives 
$
  \|\bmmu_{\bar w^T}-\bmmu\|_2\le\mcE_{\rm stat}^{\rm MW}+\mcE_{\rm opt}^{\rm MW}(T)
$
and 
$
  \|\bmmu_{\bar w^T}-\bmmu_{w^*}\|_2\le r_{1,*}+\mcE_{\rm stat}^{\rm MW}+\mcE_{\rm opt}^{\rm MW}(T),
$
which are~(ii) and the first claim in~(i). 

For the optimality gap, $F$ is convex, so $F(\bar w^T)-F^*\le\frac1T\sum_t(F(w^t)-F^*)$, and the right-hand side of~\eqref{prop:mwu:ergodic:pf:key} is at most $\frac2T\sum_t(u_t+2r_{1,*})^2+b_m\mcE_{\rm opt}^{\rm MW}(T)^2+\tau$. Splitting $u_t+2r_{1,*}\le\big(\sqrt{2D_*^{\rm MW}}+2r_{1,*}\big) + \big(u_t-\sqrt{2D_*^{\rm MW}}\big)_+$ and using $(a+b)^2\le2a^2+2b^2$ together with~\eqref{prop:mwu:ergodic:pf:excess},
\begin{align}
  F(\bar w^T)-F^*
  &\le b_m\mcE_{\rm opt}^{\rm MW}(T)^2+\tau+4\big(\sqrt{2D_*^{\rm MW}}+2r_{1,*}\big)^2 \\
  &+4\mcE_{\rm opt}^{\rm MW}(T)^2+\frac{4\tau+8\big(\sqrt{2D_*^{\rm MW}}+2r_{1,*}\big)^2}{b_m},
\end{align}
and $(b_m+4)\mcE_{\rm opt}^{\rm MW}(T)^2=2\big(1+\frac{4}{b_m}\big){M_\infty^2}\sqrt{\Omega/T}$ gives~\eqref{prop:mwu:ergodic:dist:gap}.

We may take $\Omega=\log\frac{1}{1-\varepsilon_{n,k}}$ by Lemma~\ref{lemma:capped:kl}.
\end{proof}

\section{Statistical and optimization rates}\label{s:stat:opt:rates}

This section combines Sections~\ref{s:stability:spectral:signature}-\ref{s:mwu} into explicit statistical and optimization errors for the algorithms above, under Assumption~\ref{assump:finite:covariance:contamination} and the following setup. Let $\delta\in(0,1)$ and $\eta\in(0,1]$ and suppose that  
$
    (1+\eta)\varepsilon+\frac{\log(3/\delta)+2}{n}<\frac13.    
$
Choose $k:=\lceil(1+\eta)\varepsilon n\rceil+\lceil\log(3/\delta)\rceil$ as in \eqref{eq:choice:k}. Then $\varepsilon n\le k < n/3$ and $\varepsilon_{n,k}\le(1+\eta)\varepsilon+\frac{\log(3/\delta)+2}{n}<\frac13$, so that $m:=\frac{n-k}{k}>2$.

By Theorem~\ref{thm:favorability:subgaussian:rate}, on the event ${\rm E}(\delta)$ of Corollary~\ref{cor:favorable:event}, the sample $\tilde\bx_{1:n}$ is $(r_1,\varepsilon_{n,k})$-favorable with $r_1=O(r)$, where $r=r(n,\delta,\varepsilon)$ is the rate~\eqref{eq:subgaussian:rate:alpha:v2}. Hence \eqref{def:w:*:r:1}--\eqref{def:r:1:*} hold for every minimizer $w^*$ of~\eqref{eq:landscape:problem}, and, with $r_{1,\alpha}$, $R_\alpha$ and $\tau$ as in~\eqref{def:R:alpha:r:1:alpha}--\eqref{eq:tau}: $F$ satisfies $\IQG_{h_{\mbX,\alpha}}(m,r_{1,\alpha})$ on $\Delta_{n,k}$ for every $\alpha\ge0$ (Proposition~\ref{prop:IQG}); $\partial\hat F$ is a $(4,r_{1,*},\tau)$-weakly star-convex inexact oracle for $F$ at $w^*$ relative to $h_{\mbX}$ (Proposition~\ref{prop:weakly:convex:oracle}); and $\partial\hat F$ is a Lipschitz inexact oracle relative to the reference function of each method, with the constants $M$ of Table~\ref{tab:setting} (Lemma~\ref{lemma:lipschitz}). Since $k<n/3$, $r_{1,*}\le(1+\sqrt3)\,r_1$, $r_{1,\alpha}=2r_{1,*}+R_\alpha$, and $R_\alpha\le\sqrt\alpha\,\frac{\sqrt{2k}}{n-k}$ (Lemma~\ref{lemma:capped:diameter}). We also write $a_m:=\sqrt m-\sqrt2>0$, so that $c_m=2\sqrt2\,a_m$ and ${b}_m=\sqrt m\,a_m$ in~\eqref{eq:md:constants}.

All methods start at $w^1=\unit_{[n]}/n$ and use the parameters of Table~\ref{tab:setting}. Combining Propositions~\ref{prop:ergodic:rates}, \ref{prop:window:rates} and~\ref{prop:mwu:ergodic} with these bounds yields the following corollaries, up to absolute constants.

\begin{table}[h]
\centering
\caption{Parameters of each method. $\Omega\ge D(w^*\|w^1)$; $M$ is the relative Lipschitz constant; for MWU, ${M_\infty^2}$ plays this role and the analysis uses the spectral signature in place of the oracle properties; $\alpha>0$ only for mean-induced mirror descent.}
\label{tab:setting}
\small
\setlength{\tabcolsep}{4pt}
\begin{tabular}{@{}lllll@{}}
\toprule
Method & Reference & $M$ (or ${M_\infty^2}$) & $\Omega$ & Result\\
\midrule
Entropic MD & $h_{\KL}$ & $\frac52M_{\mbX}^2$ & $\log\frac{1}{1-\varepsilon_{n,k}}\lesssim\varepsilon_{n,k}$ & Prop.~\ref{prop:ergodic:rates}(a)\\
Euclidean MD & $h_{\E}$ & $\frac52\sqrt n\,M_{\mbX}^2$ & $\frac{k}{2n(n-k)}$ & Prop.~\ref{prop:ergodic:rates}(b)\\
Mean-induced MD & $h_{\mbX,\alpha}$ & $M_{\mbX}\sqrt{4+\frac{nM_{\mbX}^2}{4\alpha}}$ & $2M_{\mbX}^2\varepsilon_{n,k}^2+\frac{\alpha k}{2n(n-k)}$ & Prop.~\ref{prop:window:rates}\\
MWU & $D_{\KL}$ & ${M_\infty^2}=2\max\{4M_{\mbX}^2,\|\hat\bfSigma\|\}$ & $\log\frac{1}{1-\varepsilon_{n,k}}\lesssim\varepsilon_{n,k}$ & Prop.~\ref{prop:mwu:ergodic}\\
\bottomrule
\end{tabular}
\end{table}

\begin{corollary}[Entropic and Euclidean mirror descent]\label{cor:rates:md}
  On the event ${\rm E}(\delta)$ of probability at least $1-\delta$, after $T$ iterations with step size $\gamma_t\equiv\frac{1}{M}\sqrt{2\Omega/T}$, the ergodic average $\bar w^T=\frac{1}{T}\sum_{t=1}^Tw^t\in\Delta_{n,k}$ satisfies
  \begin{align}
  \lambda_1\big(\bfSigma_{\bar w^T}(\bmmu)-\bfSigma\big)
  &\lesssim\underbrace{\big(1+a_m^{-2}\big)\frac{r_1^2}{\varepsilon_{n,k}}+\big(1+a_m^{-1}\big)\tau}_{\text{statistical}}
  +\underbrace{\big(1+a_m^{-1}\big)M_{\mbX}^2\sqrt{\frac{\varepsilon_{n,k}}{T}}}_{\text{optimization}},\\
  \|\bmmu_{\bar w^T}-\bmmu\|_2
  &\lesssim\underbrace{\big(1+a_m^{-1}\big)r_1+\big(1+a_m^{-1/2}\big)\sqrt{\tau\varepsilon_{n,k}}}_{\text{statistical}}
  +\underbrace{\big(1+a_m^{-1}\big)^{1/2}M_{\mbX}\,\varepsilon_{n,k}^{3/4}\,T^{-1/4}}_{\text{optimization}}.
\end{align}  
\end{corollary}

\begin{corollary}[Mean-induced mirror descent]\label{cor:rates:mean:induced}
  Let $K:=\sqrt{2m}-1$. On the event ${\rm E}(\delta)$ of probability at least $1-\delta$, after $T$ iterations with step size $\gamma_t:=\frac{2}{c_m(t+K)}$, the last iterate satisfies
  \begin{align}
  \|\bmmu_{w^T}-\bmmu\|_2
  \lesssim\underbrace{\big(1+a_m^{-1}\big)\left(r_1+\sqrt{\alpha}\frac{\sqrt{2k}}{n-k}\right)+\Big(\frac{\tau}{a_m}\Big)^{1/2}}_{\text{statistical}}
  +\underbrace{\frac{K\sqrt\Omega}{T+K}+\frac{M}{a_m\sqrt{T+K}}}_{\text{optimization}}.
\end{align}
For the windowed average with parameters $\theta\in(0,1)$, $S=(1-\theta)T$, 
\begin{align}
  \lambda_1\big(\bfSigma_{\bar w^T}(\bmmu)-\bfSigma\big)
  &\lesssim\underbrace{\Big(1+a_m^{-2}+\frac{1}{a_m\theta}\Big)\frac{\left(r_1+\sqrt{\alpha}\frac{\sqrt{2k}}{n-k}\right)^2}{\varepsilon_{n,k}}
  +\Big(\frac1\theta+\frac1{a_m}\Big)\tau}_{\text{statistical}}
  +\underbrace{\mcE_{\rm opt}^{\Sigma}(T)}_{\text{optimization}},\\
  \|\bmmu_{\bar w^T}-\bmmu\|_2
  &\lesssim\underbrace{\Big(1+a_m^{-1}+\frac{1}{\sqrt{a_m\theta}}\Big)\left(r_1+\sqrt{\alpha}\frac{\sqrt{2k}}{n-k}\right)
  +\Big(\frac1\theta+\frac1{a_m}\Big)^{1/2}\sqrt{\tau\varepsilon_{n,k}}}_{\text{statistical}}\\
  &+\underbrace{\sqrt{\varepsilon_{n,k}\,\mcE_{\rm opt}^{\Sigma}(T)}}_{\text{optimization}},
\end{align}
  where, for $\mcE_{\rm stat}(\alpha)$ and $\mcE_{\rm opt}^{\Sigma}(T)$ as in~\eqref{prop:ergodic:rates:stat:error} and~\eqref{prop:window:rates:errors:sigma:opt}, 
\begin{align}
  \mcE_{\rm opt}^{\Sigma}(T) &\lesssim \frac{1}{\theta T}\Big[a_mK\,\mcE_{\rm stat}(\alpha)^2
  +(1+a_m)\frac{K^2\Omega}{S+K}+\left(1+\log\frac1{1-\theta}\right)\big(a_m^{-1}+a_m^{-2}\big)M^2\Big],\\
  \mcE_{\rm stat}(\alpha)^2 &\lesssim\left(1+a_m^{-1}\right)^2\left(r_1^2+\alpha\,\frac{2k}{(n-k)^2}\right)+\frac{\tau}{a_m}.
\end{align}
\end{corollary}

In the mean-induced geometry, the statistical error contains the additional term $\sqrt\alpha\,\nicefrac{\sqrt{2k}}{n-k}$, which is traded off against the Lipschitz constant $M=M_{\mbX}\sqrt{4+nM_{\mbX}^2/(4\alpha)}$, and hence against the optimization error. The choice $\alpha\asymp\nu_2^2/\varepsilon_{n,k}$ preserves the optimal statistical rate: since $n-k\ge\nicefrac{2n}{3}$,
$
  \sqrt\alpha\,\nicefrac{\sqrt{2k}}{n-k}
  \le\nicefrac{3}{\sqrt2}\sqrt{\nicefrac{\alpha\,\varepsilon_{n,k}}{n}}
  \lesssim\nicefrac{\nu_2}{\sqrt n}\le\nu_2\sqrt{\nicefrac{\log(3/\delta)}{n}}.
$

\begin{corollary}[Multiplicative weights update]\label{cor:rates:mwu}
  On the event ${\rm E}(\delta)$ of probability at least $1-\delta$, after $T$ iterations with step size $\gamma_t\equiv\sqrt{\Omega/T}$, the ergodic average $\bar w^T=\frac{1}{T}\sum_{t=1}^Tw^t\in\Delta_{n,k}$ satisfies
  \begin{align}
  \lambda_1\big(\bfSigma_{\bar w^T}(\bmmu)-\bfSigma\big)
  &\lesssim\underbrace{\big(1+a_m^{-3}\big)\frac{r_1^2}{\varepsilon_{n,k}}+\big(1+a_m^{-1}\big)\tau}_{\text{statistical}}
  +\underbrace{\big(1+a_m^{-1}\big){M_\infty^2}\sqrt{\frac{\varepsilon_{n,k}}{T}}}_{\text{optimization}},\\
  \|\bmmu_{\bar w^T}-\bmmu\|_2
  &\lesssim\underbrace{\big(1+a_m^{-3/2}\big)r_1+\big(1+a_m^{-1/2}\big)\sqrt{\tau\varepsilon_{n,k}}}_{\text{statistical}}\\
  &+\underbrace{\big(1+a_m^{-1}\big)^{1/2}M_\infty\,\varepsilon_{n,k}^{3/4}\,T^{-1/4}}_{\text{optimization}}.
\end{align}  
\end{corollary}

Up to factors depending only on $\theta$ and on $a_m^{-1}$ (bounded when $\varepsilon_{n,k}$ is bounded away from $\frac13$), and for $\alpha\lesssim\nu_2^2/\varepsilon_{n,k}$ in the mean-induced geometry, the statistical errors of all four algorithms are of order $r(n,\delta,\varepsilon)$, the optimal sub-Gaussian rate~\eqref{eq:subgaussian:rate:alpha:v2}, plus a surrogate term: $\sqrt{\varepsilon_{n,k}\,\tau}$ for the averaged iterates and $\sqrt\tau\,\varepsilon_{n,k}^{1/4}$ for the last iterate of mean-induced mirror descent, where $\varepsilon_{n,k}\asymp\varepsilon+\nicefrac{\log(3/\delta)}{n}$. 

The statements in Theorem \ref{thm:informal:rates} and Table~\ref{tab:rates} in the introduction follow from the Corollaries \ref{cor:rates:md}, \ref{cor:rates:mean:induced} and \ref{cor:rates:mwu} and observations of this section with $\theta=\frac12$. Lower-order terms are omitted.

\begin{remark}[Approximate top eigenvector computation]\label{rem:approx:eigenvector:details}
Recall the setting of Remark \ref{rem:approx:eigenvector} in Section \ref{s:main:results} where, we assume that, for all $t\in\mathbb{N}$, $\hat\bfSigma$ is such that $\bfSigma_{w^t}-\hat\bfSigma\succeq\mathbf{0}$ and $\bv^t\in\mbS^{d-1}$ is only a $\chi$-approximate top eigenvector. We keep the definitions of $s^t$, $c^t$ and $g^t:=s^t+c^t\unit_{[n]}$ of Section~\ref{s:main:results}, now with this $\bv^t$.

Tracing the proofs of Proposition~\ref{prop:weakly:convex:oracle} and of the MWU regret argument, one obtains\footnote{The relative Lipschitz bounds (Lemma~\ref{lemma:lipschitz}), the spectral signature and the MWU regret bound hold for every unit vector and are unaffected.}
$$
  \langle g^t,w^*-w^t\rangle \le -(1-\chi)(F(w^t)-F^*)+(2-\chi)(\sqrt{2D_{\mbX}(w^*\|w^t)}+r_{1,*})^2+\tau_\chi.
$$
The quantity 
\begin{align}
  \tau_\chi := \tau+\chi(\lambda_d(\bfSigma-\hat\bfSigma))_++\chi\,r_1^2/\varepsilon_{n,k}
\end{align}
bounds $\tau+\chi\lambda_d(\bfSigma-\hat\bfSigma)+\chi F^*$ (we use $F^*\le r_1^2/\varepsilon_{n,k}$). Combining with the quadratic growth of Proposition~\ref{prop:IQG} and $D_{\mbX}\le D_{\mbX,\alpha}$, one obtains
$$
  \langle g^t,w^*-w^t\rangle\le-(1-\nicefrac\chi2)\Phi_{m_\chi,r_{1,\alpha},r_{1,*}}(D_{\mbX,\alpha}(w^*\|w^t))+\tau_\chi, 
$$
with $m_\chi:=\frac{2(1-\chi)}{2-\chi}m$ and $m=\frac{1-\varepsilon_{n,k}}{\varepsilon_{n,k}}$. Hence all convergence proofs for mirror descent and MWU go through with $\Phi_{m,r,\zeta}$ replaced by $(1-\frac\chi2)\Phi_{m_\chi,r,\zeta}$ and $\tau$ by $\tau_\chi$. In particular, the previous threshold $\varepsilon_{n,k}<\frac{1}{3}$ becomes $m_\chi>2$, i.e., $\varepsilon_{n,k}<\frac{1-\chi}{3-2\chi}$ and the constant $a_m$ is replaced by $a_{m_\chi}$, which is bounded away from zero for this threshold. In addition, only the decreasing step sizes of mean-induced mirror descent change, becoming $\gamma_t=\frac{2}{(1-\chi/2)\,c_{m_\chi}(t+K)}$ with $K\ge\sqrt{2m_\chi}-1$.
\end{remark}

\section{Simulations}\label{s:simulations}

\emph{Coordinatewise spike contamination.} We draw $n=500$ samples from $\calN(\mathbf{0},\bfI_d)$, $d\in\{5,10\}$,
and replace $k=\mathrm{round}(\varepsilon n)$ randomly chosen rows by the spike $3\be_1$, for
$\varepsilon\in\{0.1,0.2,0.3,0.4,0.45,0.495\}$; the point $\varepsilon=0$ in the plots is the uncontaminated sample
mean. Values $\varepsilon\ge\frac13$ lie outside the regime $\varepsilon_{n,k}<\frac13$ covered by our theory. All
methods use $\hat\bfSigma=\mathbf{0}$ and run for $T=1000$ iterations. Mean-induced MD uses the raw score $g^t$ and
$\alpha\in\{10,n/k\}$, solves its proximal subproblem with FISTA~\cite{2009Beck:Teboulle}, and follows the step-size
schedule of Corollary~\ref{cor:rates:mean:induced}; for $m\le2$, i.e.\ $\varepsilon\ge\frac13$, where this schedule
is not defined, we replace $c_m$ by $2\sqrt2\sqrt m$. PSM, Hedge and MWU (score $s^t$) use uniform normalization,
with constant $\frac52M_{\mbX^2}$ for Hedge and MWU and $\frac52\sqrt n\,M_{\mbX^2}$ for PSM. We report the last
iterate averaged over $10$ seeds; the windowed average over $t\ge T/2$ differs from it by at most $0.02$. The plots
also show the contaminated sample mean (dotted), the mean of the surviving clean points or the top eigenvalue of the
clean set (dash-dotted), and the floor $\sqrt{d/(n-k)}$ (dashed).

\begin{figure}[h]
  \centering
  \includegraphics[width=\textwidth]{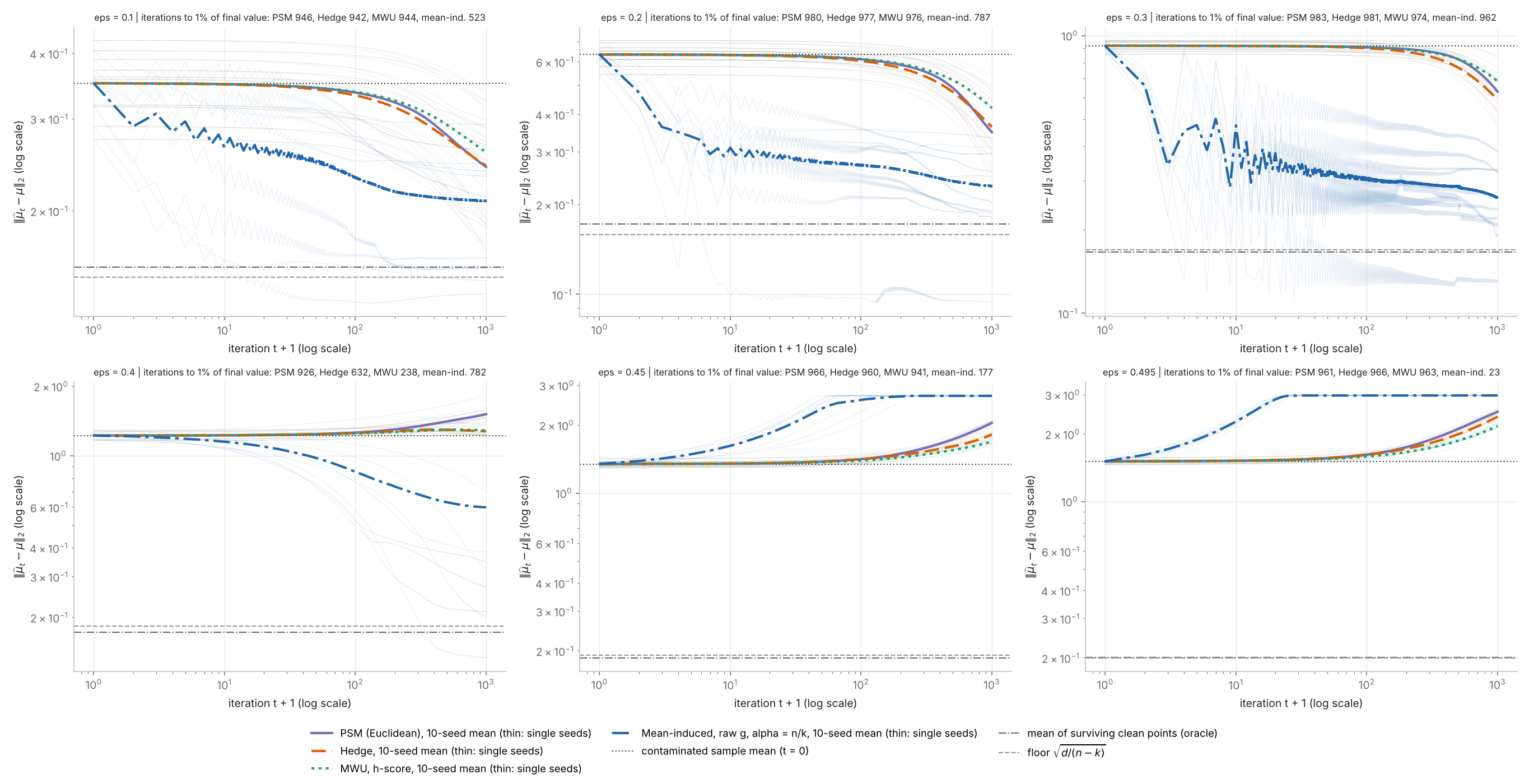}
  \caption{Estimation error $\|\bmmu_{w^t}-\bmmu\|_2$ against $t$ for $d=10$ and each $\varepsilon>0$, with the step
  sizes of Corollaries~\ref{cor:rates:md} and~\ref{cor:rates:mwu} for PSM, Hedge and MWU and of
  Corollary~\ref{cor:rates:mean:induced} for mean-induced MD ($\alpha=n/k$). Thick lines: average over $10$ seeds;
  thin lines: individual seeds.}
  \label{fig:spike:error:iter:d10:const}
\end{figure}

Figure~\ref{fig:spike:error:iter:d10:const} compares mean-induced MD ($\alpha=n/k$) with PSM, Hedge and MWU run with
the step sizes of Corollaries~\ref{cor:rates:md} and~\ref{cor:rates:mwu}, with $\Omega$ as in
Table~\ref{tab:setting}. As expected, mean-induced MD converges much faster: after $T=1000$ iterations the other three
methods are still far from converged. The remaining figures compare mean-induced MD ($\alpha=10$) with the other three methods
run with the more favorable step-decay schedule of~\cite{2019Ge:Kakade:Kidambi:Netrapalli}, $\gamma_0=\frac12$ halved
every $100$ iterations, which is not covered by our theory.

\begin{figure}[h]
  \centering
  \includegraphics[width=\textwidth]{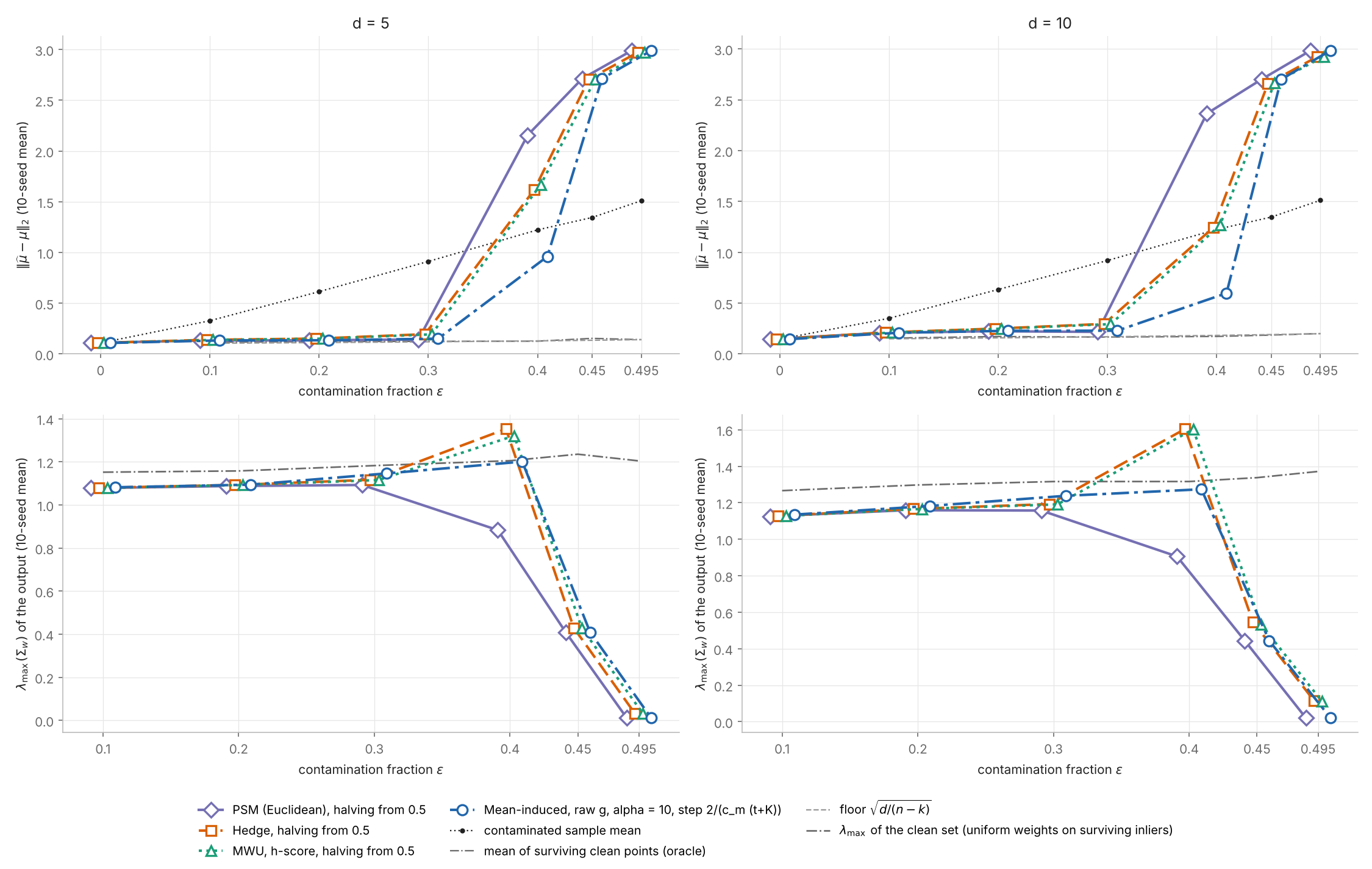}
  \caption{Estimation error $\|\bmmu_w-\bmmu\|_2$ (top) and $\lambda_1(\bfSigma_w)$ (bottom) of the last iterate
  against $\varepsilon$, for $d=5$ (left) and $d=10$ (right). Markers are offset horizontally for readability.}
  \label{fig:spike:eps}
\end{figure}

\begin{figure}[h]
  \centering
  \includegraphics[width=0.85\textwidth]{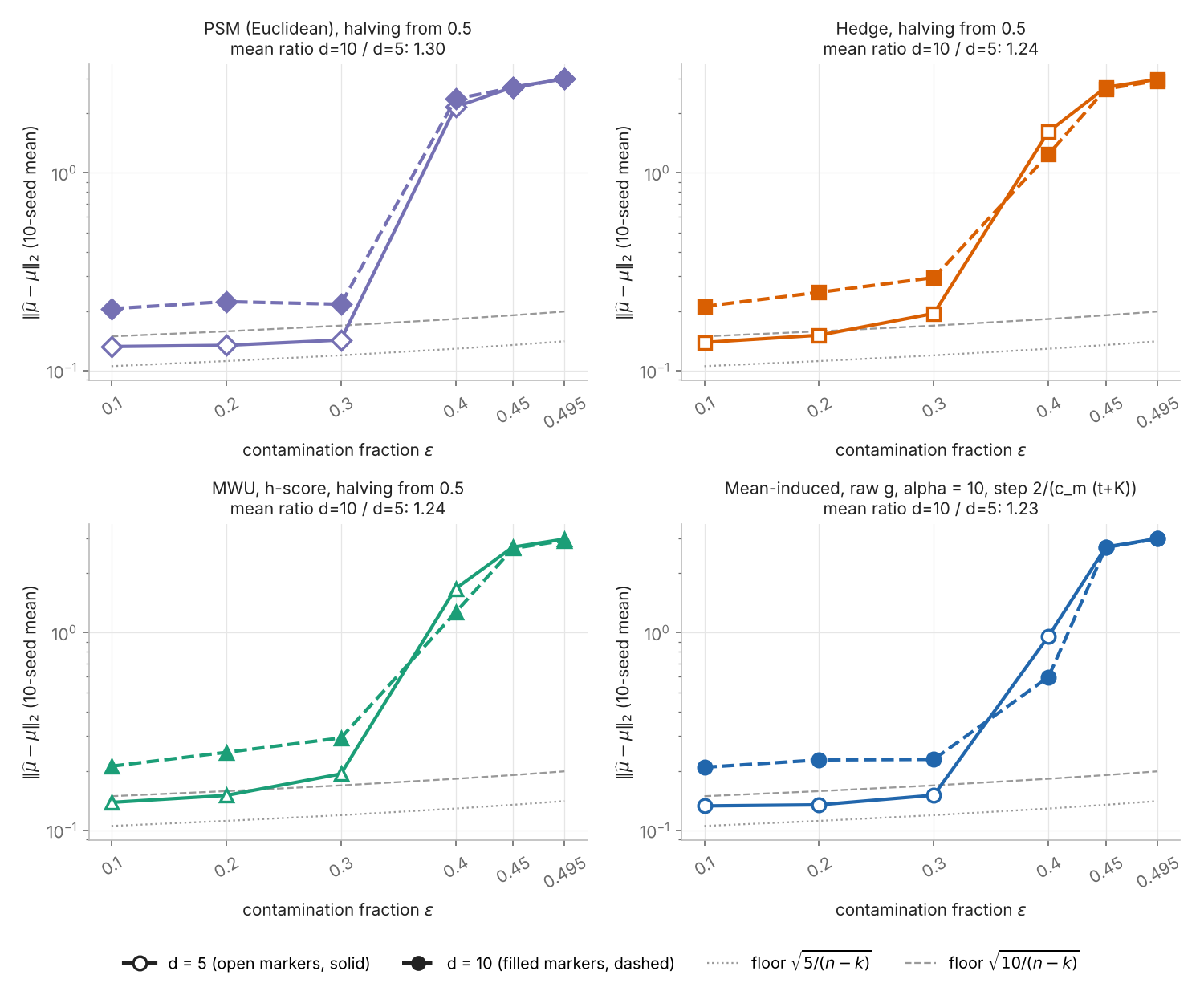}
  \caption{Effect of the dimension: estimation error against $\varepsilon$ (log scale) for $d=5$ (open markers, solid
  lines) and $d=10$ (filled markers, dashed lines), with the floors $\sqrt{5/(n-k)}$ (dotted) and
  $\sqrt{10/(n-k)}$ (dashed), whose ratio is $\sqrt2\approx1.41$. Panel titles give the ratio of the $d=10$ to the
  $d=5$ error, averaged over $\varepsilon$.}
  \label{fig:spike:dimension}
\end{figure}

\begin{figure}[h]
  \centering
  \includegraphics[width=\textwidth]{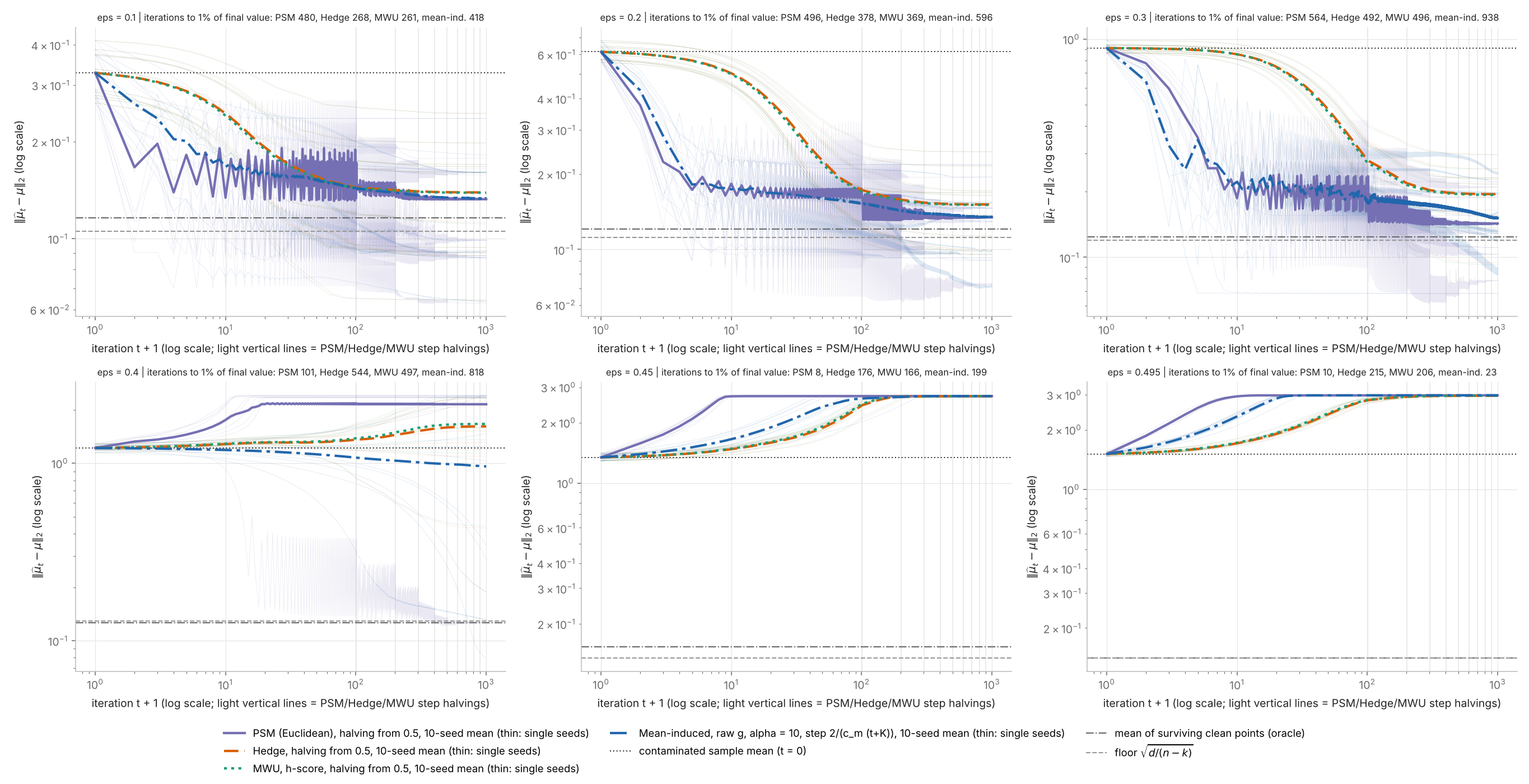}
  \caption{Estimation error $\|\bmmu_{w^t}-\bmmu\|_2$ against $t$ for $d=5$ and each $\varepsilon>0$. Thick lines:
  average over $10$ seeds; thin lines: individual seeds. Panel titles give the number of iterations each method
  needs to come within $1\%$ of its final value.}
  \label{fig:spike:error:iter:d5}
\end{figure}

Figures~\ref{fig:spike:eps} and~\ref{fig:spike:dimension} plot the errors against $\varepsilon$, the latter showing
the effect of $d$. For $\varepsilon\le0.3$, all four methods come close to the oracle. Euclidean and mean-induced MD
perform best and almost identically. Entropic MD and MWU are practically
indistinguishable from each other and fall behind as $\varepsilon$ grows. All four methods break down between $\varepsilon=0.3$ and $0.4$, close to the threshold
$\frac13$. Mean-induced MD degrades most gracefully: at $\varepsilon=0.4$ it is the only method whose average error
remains below that of the sample mean. For $\varepsilon\ge0.45$ (and already at $\varepsilon=0.4$ for PSM), the output has a smaller top eigenvalue than the clean set but a large error (Figure~\ref{fig:spike:eps}, bottom): the optimization succeeds, but a small $\lambda_1(\bfSigma_w)$ no longer certifies an accurate mean, consistent with the requirement $\varepsilon_{n,k}<\frac13$. Figure~\ref{fig:spike:error:iter:d5} plots the error against $t$ for $d=5$.
Within ten iterations, Euclidean and mean-induced MD bring the error down to $0.17$--$0.25$, whereas entropic MD and
MWU are still at $0.24$--$0.84$ and need about $100$ iterations to catch up.

\appendix

\section{Remaining proofs}\label{s:appendix}

\subsection{Proof of Lemma \ref{lemma:capped:tv}}\label{app:capped:tv}

We prove the bound on the total variation. Let $w\in\Delta_{n,k}$. By the identity $\frac1n=\frac{1-\varepsilon_{n,k}}{n-k}$, $w_i-\frac{1}{n}\le w_i-(1-\varepsilon_{n,k})w_i=\varepsilon_{n,k} w_i$ for every $i$. Since $\varepsilon_{n,k} w_i\ge0$, taking positive parts preserves this inequality, so $(w_i-\frac{1}{n})_+\le\varepsilon_{n,k} w_i$. Summing over $i$ and using $\sum_iw_i=1$ gives
$
  \mathrm{TV}(w,\nicefrac{\unit_{[n]}}{n})=\sum_{i=1}^n(w_i-\frac{1}{n})_+\le\varepsilon_{n,k}\sum_{i=1}^nw_i=\varepsilon_{n,k}.
$

Next, we prove the condition on equality on the total variation. Let $w=\frac{\unit_S}{n-k}$ with $|S|=n-k$. Since $w$ and $\frac{\unit_{[n]}}{n}$ both have total mass one,
$\sum_i(w_i-\frac{1}{n})=0$, implying
$
  \mathrm{TV}(w,\nicefrac{\unit_{[n]}}{n})=\sum_{i=1}^n(w_i-\frac{1}{n})_+=\sum_{i=1}^n(\frac{1}{n}-w_i)_+.
$
We have $w_i=\frac1{n-k}$ for $i\in S$ and $w_i=0$ for $i\notin S$. Using the second expression,
$
    \mathrm{TV}(w,\nicefrac{\unit_{[n]}}{n})=\sum_{i\notin S}\frac1n=\frac{n-|S|}{n}=\frac kn=\varepsilon_{n,k}.
$

\subsection{Proof of Lemma \ref{lemma:weight:dif}}\label{app:weight:dif}

Let $\delta_i:=w_i-\frac1n$ and write $\delta=\delta_+-\delta_-$ with
$\delta_\pm\ge0$. Since $\sum_i\delta_i=0$, we have
$\sum_i(\delta_+)_i=\sum_i(\delta_-)_i=\mathrm{TV}(w,\unit_{[n]}/n)\le\varepsilon_{n,k}$
by Lemma~\ref{lemma:capped:tv}. Hence
$
    \sum_i\delta_ia_i=\sum_i(\delta_+)_ia_i-\sum_i(\delta_-)_ia_i
    \le\mathrm{TV}(w,\nicefrac{\unit_{[n]}}{n})(\max_ia_i-\min_ia_i),
$
which gives both inequalities. If $a_i\ge0$, drop the nonpositive term
$-\sum_i(\delta_-)_ia_i$ to get $\sum_i\delta_ia_i\le\mathrm{TV}(w,\unit_{[n]}/n)\max_ia_i$. 

\subsection{Proof of Lemma~\ref{lemma:capped:diameter}}\label{app:capped:diameter}

\emph{Diameter.} By Lemma~\ref{lemma:capped:extreme},
$\Delta_{n,k}\times\Delta_{n,k}$ is the convex hull of the pairs of vertices, so the convex function $(w,w')\mapsto\|w-w'\|_2^2$ attains its maximum at a pair $w=\unit_S/(n-k)$, $w'=\unit_{S'}/(n-k)$ with $|S|=|S'|=n-k$. Then $w-w'$ equals $\pm\frac1{n-k}$ on $S\triangle S'$ and $0$ elsewhere, so 
$
  \|w-w'\|_2^2=\frac{|S\triangle S'|}{(n-k)^2}=\frac{2(n-k)-2|S\cap S'|}{(n-k)^2}.
$
Since $|S\cap S'|\ge\max\{0,n-2k\}$, with equality when $S'$ contains as much of $[n]\setminus S$ as possible, the maximum equals
$\frac{2(n-k)-2\max\{0,n-2k\}}{(n-k)^2}=\frac{2\min\{k,n-k\}}{(n-k)^2}$.

\emph{Distance to the uniform distribution.} Since $\langle w,\unit_{[n]}/n\rangle=\|\unit_{[n]}/n\|_2^2=\frac1n$ for $w\in\Delta_n$, we have $\|w-\unit_{[n]}/n\|_2^2=\|w\|_2^2-\frac1n$. Moreover, $\|w\|_2^2\le\max_iw_i\sum_iw_i\le\frac{1}{n-k}$ for $w\in\Delta_{n,k}$, with
equality at the vertices. Hence the maximum equals $\frac{1}{n-k}-\frac1n=\frac{k}{n(n-k)}$.

\subsection{Proof of Lemma \ref{lemma:calM:d:convex:hull}}\label{app:calM:d:convex:hull}

Each $\bu^{\otimes2}$ with $\bu\in\mbS^{d-1}$ is positive semidefinite with
$\tr(\bu^{\otimes2})=\|\bu\|_2^2=1$, and $\calM_d$ is convex, so it contains the
convex hull of these points. Conversely, by the spectral theorem, write $\bfM\in\calM_d$ as $\bfM=\sum_{j=1}^d\lambda_j\bv_j^{\otimes2}$ with orthonormal eigenvectors $\bv_j\in\mbS^{d-1}$ and eigenvalues $\lambda_j\ge0$ satisfying
$\sum_j\lambda_j=\tr(\bfM)=1$; this is a convex combination. For the consequence,
$\langle\bfA,\bfM\rangle=\sum_j\lambda_j\bv_j^\top\bfA\bv_j\le\lambda_{\max}(\bfA)$,
with equality when $\bfM=\bu^{\otimes2}$ for a top unit eigenvector $\bu$ of $\bfA$; finally, $\lambda_{\max}(\bfA)=\|\bfA\|$ since $\bfA\succeq0$. 

\subsection{Proof of Lemma \ref{lemma:trim:trunc}}\label{app:trim:trunc}

First note that
$
  \min_{w\in\Delta_{n,k}}\sfT_{2}(\bz_{1:n}\,|\,w,\bfM)
  =\frac{1}{n-k}\sum_{i\in S}\langle\bz_i^{\otimes2},\bfM\rangle,
$
where $S=S(\bz_{1:n},\bfM)\subseteq[n]$ is the set as stated in  Lemma~\ref{lemma:min:weight:trimming}. Now, we have that on the
event $\Count_2(\bz_{1:n},B,k)$, $\|\bfC_{\bfM}\bz_\ell\|_2>B$, or equivalently
$\langle\bz_\ell^{\otimes2},\bfM\rangle>B^2$, can only hold for at most $k$ indices
$\ell\in[n]$. Therefore, for all $i\in S$, 
$\langle\bz_i^{\otimes2},\bfM\rangle\le B^2$ and 
$\langle\bz_i^{\otimes2},\bfM\rangle=\langle\bz_i^{\otimes2},\bfM\rangle\wedge B^2$. As a result,
$
  \sfT_{2,k}(\bz_{1:n}\,|\,\bfM)
  =\frac{1}{n-k}\sum_{i\in S}\langle\bz_i^{\otimes2},\bfM\rangle
  =\frac{1}{n-k}\sum_{i\in S}\big(\langle\bz_i^{\otimes2},\bfM\rangle\wedge B^2\big).
$
Since $t\in[0,\infty)\mapsto t\wedge B^2$ is nondecreasing, the average in the
right-hand side of the preceding display is taken over the indices achieving the
$n-k$ smallest values of $\langle\bz_i^{\otimes2},\bfM\rangle\wedge B^2$. Taking an
average over all indices can only give a larger result, so
$
  \frac{1}{n-k}\sum_{i\in S}\big(\langle\bz_i^{\otimes2},\bfM\rangle\wedge B^2\big)
  \le\frac{1}{n}\sum_{i=1}^n\big(\langle\bz_i^{\otimes2},\bfM\rangle\wedge B^2\big)
  =\sfT_{2,B}(\bz_{1:n}\,|\,\bfM).
$
On the other hand, since
$\langle\bz_i^{\otimes2},\bfM\rangle\wedge B^2\in[0,B^2]$ and $|[n]\setminus S|=k$,
$
  \frac{1}{n-k}\sum_{i\in S}\big(\langle\bz_i^{\otimes2},\bfM\rangle\wedge B^2\big)
  \ge\frac{1}{n}\sum_{i\in S}\big(\langle\bz_i^{\otimes2},\bfM\rangle\wedge B^2\big)
  \ge\sfT_{2,B}(\bz_{1:n}\,|\,\bfM)-\frac{B^2k}{n}.
$

\subsection{Proof of Lemma \ref{lemma:CM:moment}}\label{app:CM:moment}

Write $\bfM=\sum_{j=1}^d\lambda_j\bu_j^{\otimes2}$ with $\lambda_j\ge0$, $\sum_j\lambda_j=\tr(\bfM)=1$ and $\bu_j\in\mbS^{d-1}$. Then $\|\bfC_{\bfM}\bh\|_2^2=\langle\bh,\bfM\bh\rangle=\sum_j\lambda_j\langle\bu_j,\bh\rangle^2$, a convex combination. The claim follows from
$\langle\bu_j,\bh\rangle^2\le\|\bh\|_2^2$. For the second, since $x\mapsto x^{p/2}$ is convex on $[0,\infty)$ for $p\ge2$, Jensen's inequality gives
$
    \|\bfC_{\bfM}\bxi\|_2^p
    =(\sum_j\lambda_j\langle\bu_j,\bxi\rangle^2)^{p/2}
    \le\sum_j\lambda_j|\langle\bu_j,\bxi\rangle|^p,
$
and taking expectations, $\esp\|\bfC_{\bfM}\bxi\| _2^p\le\sum_j\lambda_j\nu_p^p=\nu_p^p$.

\subsection{Proof of Proposition \ref{prop:multvect}}\label{app:multvect}

When $\bg\sim\Gamma$, $\langle\bg,\bfC_{\bfM}\bh\rangle$ is a one-dimensional
normal random variable with mean $0$ and variance
$\|\bfC_{\bfM}\bh\|_2^2=\langle\bh,\bfM\bh\rangle$. Therefore,
$\sqrt{\nicefrac{2}{\pi}}\,\|\bfC_{\bfM}\bh\|_2=\Gamma|\langle\bg,\bfC_{\bfM}\bh\rangle|$.
Now let $\bar\bg:=\bfC_{\bfM}^\top\bg/\|\bfC_{\bfM}^\top\bg\|_2$ if
$\bfC_{\bfM}^\top\bg\neq\bzero$; otherwise, let $\bar\bg$ be any fixed vector in
$\mbS^{d-1}$. In both cases,
$\langle\bg,\bfC_{\bfM}\bh\rangle=\langle\bfC_{\bfM}^\top\bg,\bh\rangle
=\|\bfC_{\bfM}^\top\bg\|_2\,\langle\bar\bg,\bh\rangle$. Since $\tr(\bfM)=1$, we have
$\bfC_{\bfM}\neq\bzero$, so $\Gamma\|\bfC_{\bfM}^\top\bg\|_2>0$, and by Jensen's
inequality,
$\Gamma\|\bfC_{\bfM}^\top\bg\|_2\le(\Gamma\|\bfC_{\bfM}^\top\bg\|_2^2)^{1/2}
=\tr(\bfM)^{1/2}=1$. Then
\begin{align}
  \Gamma|\langle\bg,\bfC_{\bfM}\bh\rangle|
  &=\int_{\re^d}\|\bfC_{\bfM}^\top\bg\|_2\,|\langle\bar\bg,\bh\rangle|\,\Gamma(d\bg)\\
  &=\big(\Gamma\|\bfC_{\bfM}^\top\bg\|_2\big)\int_{\re^d}
  \frac{\|\bfC_{\bfM}^\top\bg\|_2}{\Gamma\|\bfC_{\bfM}^\top\bg\|_2}\,
  |\langle\bar\bg,\bh\rangle|\,\Gamma(d\bg)\\
  &=\big(\Gamma\|\bfC_{\bfM}^\top\bg\|_2\big)\int_{\mbS^{d-1}}
  |\langle\bu,\bh\rangle|\,\bar{\Gamma}_{\bfC_{\bfM}}(d\bu),
\end{align}
where $\bar{\Gamma}_{\bfC_{\bfM}}$ is defined by
$
  \bar{\Gamma}_{\bfC_{\bfM}}(A):=\int_{\re^d}
  \frac{\|\bfC_{\bfM}^\top\bg\|_2}{\Gamma\|\bfC_{\bfM}^\top\bg\|_2}\,
  \unit_{\{\bar\bg\in A\}}\,\Gamma(d\bg)
$
for all Borel sets $A\subseteq\mbS^{d-1}$. This is a probability measure, since its density $\|\bfC_{\bfM}^\top\bg\|_2/\Gamma\|\bfC_{\bfM}^\top\bg\|_2$ is nonnegative with $\Gamma$-mean one, and $\bar\bg\in\mbS^{d-1}$ always. This proves the identity.

\subsection{Proof of Lemma \ref{lemma:bias}}\label{app:bias}

Fix $\bv\in\mbS^{d-1}$. Since
${\sf clip}_B(\langle\bv,\bxi\rangle)+{\sf rem}_B(\langle\bv,\bxi\rangle)=\langle\bv,\bxi\rangle$
has mean zero,
$\esp[{\sf clip}_B(\langle\bv,\bxi\rangle)]=-\esp[{\sf rem}_B(\langle\bv,\bxi\rangle)]$.
We use the facts that $|{\sf rem}_B(t)|=(|t|-B)_+\le|t|^p/B^{p-1}$ for $t\in\re$ and
$p\ge1$, and $(s-B^2)_+\le s^{p/2}/B^{p-2}$ for $s\ge0$ and $p\ge2$. The first gives
$
  \big|\esp[{\sf rem}_B(\langle\bv,\bxi\rangle)]\big|
  \le\frac{\esp|\langle\bv,\bxi\rangle|^p}{B^{p-1}}
  \le\frac{\nu_p^p}{B^{p-1}}.
$

For $\bfM\in\calM_d$, we have
$\langle\bfSigma,\bfM\rangle=\esp\langle\bxi^{\otimes2},\bfM\rangle
=\esp\|\bfC_{\bfM}\bxi\|_2^2$, and ${\sf clip}_B^2(t)=t^2\wedge B^2$ for $t\ge0$. The
second fact, with $s=\|\bfC_{\bfM}\bxi\|_2^2$, gives
$
  0\le\|\bfC_{\bfM}\bxi\|_2^2-{\sf clip}_B^2(\|\bfC_{\bfM}\bxi\|_2)
  =\big(\|\bfC_{\bfM}\bxi\|_2^2-B^2\big)_+
  \le\frac{\|\bfC_{\bfM}\bxi\|_2^p}{B^{p-2}} .
$
Taking expectations and applying Lemma~\ref{lemma:CM:moment},
$
  0\le\langle\bfSigma,\bfM\rangle-\esp\big[{\sf clip}_B^2(\|\bfC_{\bfM}\bxi\|_2)\big]
  \le\frac{\esp\|\bfC_{\bfM}\bxi\|_2^p}{B^{p-2}}
  \le\frac{\nu_p^p}{B^{p-2}}. 
$

\subsection{Proof of Corollary \ref{cor:favorable:event}}\label{app:favorable:event}

Since $s\in\mbN$ and $s\ge\log(3/\delta)$, Lemma~\ref{lemma:counting} gives the tail bound $\prob\{\Count_2(\bxi_{1:n},B,s)\}\ge1-e^{-s}\ge1-\delta/3$, and Lemma~\ref{lemma:truncated:processes} gives $\prob({\rm E}_1(\delta/3)),\prob({\rm E}_2(\delta/3))\ge1-\delta/3$. By the union bound, $\prob({\rm E}(\delta))\ge1-\delta$.

Now fix the event ${\rm E}(\delta)$. By Assumption~\ref{assump:finite:covariance:contamination}, $\tilde\bxi_i=\bxi_i$
for all but at most $\lfloor\varepsilon n\rfloor$ indices, so, for every
$\bfM\in\calM_d$, $\#\{i:\|\bfC_{\bfM}\tilde\bxi_i\|_2>B\}\le s+\lfloor\varepsilon n\rfloor=k$; this is item~(i) of Definition~\ref{def:favorable:counting}. By Lemma~\ref{lemma:truncated:processes} with $\delta/3$ in place of $\delta$, and since $\frac{28}{3}\le\frac{224\sqrt{2\pi}}{3}$, $\frac{4B}{3n}\le\frac{8B}{3n}$ and $\log(3/\delta)\le s$, we have $S_1\le \bar r_1(\delta)$ and $S_2\le B\bar r_1(\delta)$; these are items~(ii) and~(iii).

\subsection{Proof of Proposition \ref{prop:favorable:stable}}\label{app:favorable:stable}

By Proposition~\ref{prop:counting:favorable}, $\tilde\bx_{1:n}$ is $(r_1,\varepsilon_{n,k})$-favorable; by item~(b) of Definition~\ref{def:favorable}, the minimum defining $w^*$ is at most $r_1^2/\varepsilon_{n,k}$, so
$
  \lambda_1(\bfSigma_{w^*}(\bmmu)-\bfSigma)\le r_1^2/\varepsilon_{n,k}.
$
By item (i) in Definition \ref{def:favorable:counting} and Assumption \ref{assump:finite:covariance:contamination}, $\Count_2(\tilde\bxi_{1:n},B,k)$ holds. Moreover, $\Count_2(\tilde\bxi_{1:n},B,k)$ implies $\Count_2(\tilde\bxi_{1:n},B,2k)$, items~(ii)--(iii) of Definition~\ref{def:favorable:counting} do not involve $k$, and $\varepsilon n\le2k<n$; hence $\tilde\bx_{1:n}$ is also $(B,\bar r_1,\varepsilon_{n,2k})$-favorable and, by Proposition~\ref{prop:counting:favorable}, $(R_1,\varepsilon_{n,2k})$-favorable. In particular,
\begin{enumerate}
  \item for all $w\in\Delta_{n,2k}$,
  $
    \|\bmmu_w-\bmmu\|_2\le R_1 + \sqrt{\frac{\varepsilon_{n,2k}}{1-\varepsilon_{n,2k}}
    \Big(\lambda_1(\bfSigma_w(\bmmu)-\bfSigma)+\frac{R_1^2}{\varepsilon_{n,2k}}\Big)};
  $
  \item $\min_{w\in\Delta_{n,2k}}\lambda_d\big(\bfSigma_w(\bmmu)-\bfSigma\big)\ge-R_1^2/\varepsilon_{n,2k}$.
\end{enumerate}

Let $w\in\Delta_{n,k}(w^*)$, i.e., $(1-\varepsilon_{n,k})w_i\le w_i^*$ for all $i\in[n]$. Then $w_i\le\frac{1}{(1-\varepsilon_{n,k})(n-k)}=\frac{n}{(n-k)^2}\le\frac{1}{n-2k}$, since $(n-k)^2\ge n(n-2k)$; thus $w\in\Delta_{n,2k}$.

We first bound the deviation from $\bfSigma$. Let $\eta_i:=w_i^*-(1-\varepsilon_{n,k})w_i\ge0$. Then $\bfSigma_w(\bmmu)=\frac{1}{1-\varepsilon_{n,k}}\big(\bfSigma_{w^*}(\bmmu)-\sum_i\eta_i\tilde\bxi_i^{\otimes2}\big)$, and since $\sum_i\eta_i\tilde\bxi_i^{\otimes2}\succeq0$,
$
  \bfSigma_w(\bmmu)-\bfSigma\preceq\frac{1}{1-\varepsilon_{n,k}}
  \big(\bfSigma_{w^*}(\bmmu)-\bfSigma+\varepsilon_{n,k}\bfSigma\big)
  \preceq R_2\,\bfI_d,
$
where we defined 
$
  R_2:=\frac{r_1^2/\varepsilon_{n,k}+\varepsilon_{n,k}\|\bfSigma\|}{1-\varepsilon_{n,k}}. 
$

Since $w\in\Delta_{n,2k}$, item~2 gives $\bfSigma_w(\bmmu)-\bfSigma\succeq-\frac{R_1^2}{\varepsilon_{n,2k}}\bfI_d$. Hence $\|\bfSigma_w(\bmmu)-\bfSigma\|\le\max\{R_2,R_1^2/\varepsilon_{n,2k}\}$.

We now bound the deviation from $\bmmu$. Since $\lambda_1(\bfSigma_w(\bmmu)-\bfSigma)\le R_2$, item~1 gives
$
  \|\bmmu_w-\bmmu\|_2\le R_1+\sqrt{\frac{\varepsilon_{n,2k}}{1-\varepsilon_{n,2k}}
  \Big(R_2+\frac{R_1^2}{\varepsilon_{n,2k}}\Big)}=R,
$
where the equality uses $\varepsilon_{n,2k}=2\varepsilon_{n,k}$.

Finally, $R\ge R_1$ gives $R_1^2/\varepsilon_{n,2k}\le R^2/(2\varepsilon_{n,k})$, and $R^2\ge\frac{\varepsilon_{n,2k}}{1-\varepsilon_{n,2k}}R_2\ge2\varepsilon_{n,k}R_2$ gives $R_2\le R^2/(2\varepsilon_{n,k})$. Hence, for every $w\in\Delta_{n,k}(w^*)$, $\|\bmmu_w-\bmmu\|_2\le R$ and $\|\bfSigma_w(\bmmu)-\bfSigma\|\le R^2/\varepsilon_{n,k}$, which is $(R,\varepsilon_{n,k})$-stability with respect to $w^*$.

\subsection{Proof of Lemma \ref{lemma:Phi:G}}\label{app:Phi:G}

We use the notation $D_*=D_{*,m}(r,\zeta)$ only within this proof. 

\emph{Proof of (i).} Write $s:=\sqrt{2D}$ and $s_*:=\sqrt{2D_*}
=\frac{\sqrt m\,r+\sqrt2\,\zeta}{\sqrt m-\sqrt2}$. For $s\le r$,
$\Phi_{m,r,\zeta}(D)=-4D-4\zeta\sqrt{2D}-2\zeta^2$; for $s\ge r$,
$\Phi_{m,r,\zeta}(D)=2(m-2)D-2(mr+2\zeta)\sqrt{2D}+mr^2-2\zeta^2$. Both
expressions are convex in $D$, since $D\mapsto-\sqrt D$ is convex, and both give $\frac{d\Phi_{m,r,\zeta}}{ds}=-4(r+\zeta)$ at $s=r$. Hence $\Phi_{m,r,\zeta}$ is continuously differentiable on $(0,\infty)$ with nondecreasing derivative, so it is convex on $[0,\infty)$. Since $m>2$, we have $s_*>r$ and $\sqrt m\,(s_*-r)=\sqrt2\,(s_*+\zeta)$, so $\Phi_{m,r,\zeta}(D_*)=0$ and
\begin{align}
    \Phi_{m,r,\zeta}'(D_*)
    =\frac{2m(s_*-r)-4(s_*+\zeta)}{s_*}
    =\frac{2\sqrt2\,(\sqrt m-\sqrt2)(s_*+\zeta)}{s_*}.
\end{align}
Convexity gives the tangent-line inequality at $D_*$.

\emph{Proof of (ii).} Write $s:=\sqrt{2D}$ and $s_*:=\sqrt{2D_*}$, so that
$\Phi_{m,r,\zeta}(D)=m(s-r)_+^2-2(s+\zeta)^2$ and $(\sqrt m-\sqrt2)\,s_*=\sqrt m\,r+\sqrt2\,\zeta$; in particular, $s_*\ge r$. For $s\ge s_*$, we have $(s-r)_+=s-r$ and
\begin{align}
  \Phi_{m,r,\zeta}(D)
  =\big(\sqrt m(s-r)-\sqrt2(s+\zeta)\big)\big(\sqrt m(s-r)+\sqrt2(s+\zeta)\big),
\end{align}
where the first factor equals $(\sqrt m-\sqrt2)(s-s_*)$ and the second is at least $\sqrt m(s-r)\ge\sqrt m(s-s_*)$; hence $\Phi_{m,r,\zeta}(D)\ge\sqrt m(\sqrt m-\sqrt2)(s-s_*)^2$. For $s\le s_*$, $\Phi_{m,r,\zeta}(D)\ge-2(s+\zeta)^2\ge-2(s_*+\zeta)^2$. Both cases give~(ii).

\subsection{Proof of Lemma \ref{lemma:recursion:slow}}\label{app:recursion:slow}

Let $\gamma:=\frac1M\sqrt{2\Omega/T}$, and let $c_{m,\zeta}:=c_m(1+\zeta/\sqrt{2D_*})\ge c_m>0$, so that Lemma~\ref{lemma:Phi:G}(i) gives
$\Phi_{m,r,\zeta}(\bar D_t)\ge c_{m,\zeta}(\bar D_t-D_*)$. Hence \eqref{lemma:recursion:slow:eq0} yields
$
    c_{m,\zeta}\,\gamma\,(\bar D_t-D_*)\le D_t-D_{t+1}+\gamma\tau+\frac{\gamma^2M^2}{2}.
$
Summing over $t\in[T]$, using $D_{T+1}\ge0$ and $D_1\le\Omega$, and dividing by
$c_{m,\zeta}\gamma T$,
$
    \frac1T\sum_{t=1}^T\bar D_t-D_*
    \le\frac{1}{c_{m,\zeta}}(\frac{\Omega}{\gamma T}+\frac{\gamma M^2}{2}+\tau)
    \le \frac{1}{c_m}(\frac{\Omega}{\gamma T}+\frac{\gamma M^2}{2}+\tau)
    = \frac{M}{2(\sqrt m-\sqrt2)}\sqrt{\frac{\Omega}{T}}+\frac{\tau}{c_m},
$
where the second inequality uses $c_{m,\zeta}\ge c_m$ and the nonnegativity of
the term in parentheses, and the equality uses $\frac{\Omega}{\gamma T}=\frac{\gamma M^2}{2}=M\sqrt{\frac{\Omega}{2T}}$.

\subsection{Proof of Lemma \ref{lemma:recursion:G:plus}}\label{app:recursion:G:plus}

We first prove the following auxiliary lemma.
\begin{lemma}\label{lemma:contraction}
Let $A_0\ge0$ and $K>1$, and let $(H_t)_{t\ge1}$ be a nonnegative sequence such
that, for all $t\in\mbN$,
$
    H_{t+1}\le(1-\frac{2}{t+K})H_t+\frac{A_0}{(t+K)^2}.
$
Then, for all $t\in\mbN$,
$
    H_t\le(\frac{1+K}{t+K})^2H_1
    +A_0(1+\frac{1}{1+K})^2\frac{t-1}{(t+K)^2}.
$
\end{lemma}
\begin{proof}
For $t=1$ the claim reads $H_1\le H_1$. Fix $s\in\mbN$. Since $K>1$, we have $s+K>2$, so $0<1-\frac{2}{s+K}<1$. Moreover,
$(1+\frac{1}{s+K})^2(1-\frac{2}{s+K})\le1$, since, by
$\log(1+x)\le x$ for $x>-1$, $2\log(1+\frac{1}{s+K})+\log(1-\frac{2}{s+K}) \le\frac{2}{s+K}-\frac{2}{s+K}=0$. Multiplying the assumed recursion at time $s$ by $(s+1+K)^2=(s+K)^2(1+\frac{1}{s+K})^2$ and using $H_s\ge0$ and $1+\frac{1}{s+K}\le1+\frac{1}{1+K}$,
$
    (s+1+K)^2H_{s+1}\le(s+K)^2H_s+A_0(1+\frac{1}{1+K})^2.
$
Summing over $s=1,\ldots,t-1$, the terms telescope, giving
$(t+K)^2H_t\le(1+K)^2H_1+A_0(1+\frac{1}{1+K})^2(t-1)$. Dividing by $(t+K)^2$ proves the claim.
\end{proof}

\begin{proof}[Proof of Lemma \ref{lemma:recursion:G:plus}]

By Lemma~\ref{lemma:Phi:G}(i),
$\Phi_{m,r,\zeta}(D)\ge c_{m,\zeta}(D-D_*)$, where
$c_{m,\zeta}:=c_m(1+\zeta/s_*)$ and $s_*:=\sqrt{2D_*}$. Since
$s_*\ge\frac{\sqrt2\,\zeta}{\sqrt m-\sqrt2}$, we have
$c_m\le c_{m,\zeta}\le c_m\sqrt{m/2}$, and the choice of $K$ gives
$\gamma_tc_{m,\zeta}\le\frac{\sqrt{2m}}{t+K}\le1$. Hence,
by~\eqref{lemma:recursion:G:plus:eq0}, writing $D_t-D_*=(D_t-D_*^\tau)+\frac{\tau}{c_m}$,
\begin{align}
    D_{t+1}-D_*^\tau
    &\le(1-\gamma_tc_{m,\zeta})(D_t-D_*^\tau)
    +\gamma_t\tau\left(1-\nicefrac{c_{m,\zeta}}{c_m}\right)+\frac{\gamma_t^2M^2}{2}\\
    &\le(1-\gamma_tc_{m,\zeta})(D_t-D_*^\tau)+\frac{\gamma_t^2M^2}{2},
\end{align}
since $c_{m,\zeta}\ge c_m$. Let $H_t:=(D_t-D_*^\tau)_+$. If $D_t\le D_*^\tau$, the first term is nonpositive, so
$H_{t+1}\le\frac{\gamma_t^2M^2}{2}$. If $D_t>D_*^\tau$, then, since
$\gamma_tc_{m,\zeta}\ge\gamma_tc_m$,
$H_{t+1}\le(1-\gamma_tc_m)H_t+\frac{\gamma_t^2M^2}{2}$. 
  
In both cases,
$
    H_{t+1}\le(1-\frac{2}{t+K})H_t+\frac{2M^2}{c_m^2(t+K)^2}.
$
Since $K\ge\sqrt{2m}-1>1$, Lemma~\ref{lemma:contraction} applies with $A_0:=\frac{2M^2}{c_m^2}$. Finally, $\frac{t-1}{(t+K)^2}\le\frac{1}{t+K}$ and, since $1+K\ge\sqrt{2m}>2$, $(1+\frac{1}{1+K})^2\le\frac94$, so $A_0(1+\frac{1}{1+K})^2\frac{t-1}{(t+K)^2}\le\frac{9M^2}{2c_m^2(t+K)}=\frac{A}{t+K}$.

\end{proof}

\subsection{Proof of Lemma \ref{lemma:series:bounds}}\label{app:series:bounds}

For a nonincreasing $\phi$ on $[S,T]$, $\sum_{t=S+1}^T\phi(t)\le\int_S^T\phi(x)\,dx$.
With $\phi(x)=(x+K)^{-2}$ this gives
$\sum_{t=S+1}^T(t+K)^{-2}\le\frac{1}{S+K}$, and with $\phi(x)=(x+K)^{-1}$ it gives
the second bound; both functions are positive and nonincreasing on $[S,T]$, since
$S+K>0$. Finally, for $K\ge0$, $\frac{T+K}{S+K}\le\frac TS=\frac{1}{1-\theta}$.

\subsection{Proof of Lemma \ref{lemma:capped:kl}}\label{ss:capped:KL}

By definition, and with the convention $0\log0:=0$, $D_{\KL}(w\|\unit_{[n]}/n)=\sum_{i=1}^nw_i\log(nw_i)$. For every $i$ with $w_i>0$,
the cap $w_i\le\frac{1}{n-k}$ gives $\log(nw_i)\le\log\frac{n}{n-k}=\log\frac{1}{1-\varepsilon_{n,k}}$; multiplying by $w_i$, summing over $i$, and using $\sum_iw_i=1$ gives the first inequality. At a vertex $w=\unit_S/(n-k)$, every nonzero coordinate equals $\frac{1}{n-k}$, so
equality holds. Finally, $\log y\ge1-\frac1y$ for $y>0$, applied with
$y=1-\varepsilon_{n,k}$, gives the second inequality.

\bibliographystyle{plain}
\bibliography{robust_mean_estimation_landscape_MP}

\end{document}